\documentclass[11pt,letterpaper]{article}
\pdfoutput=1
\usepackage{jheppub}
\usepackage{amsmath,amssymb,amsfonts,bm,mathrsfs}
\usepackage{graphicx,wrapfig}
\usepackage{multirow}
\usepackage{verbatim}
\usepackage{appendix}
\usepackage{rotating}

\graphicspath{{./figures/}}

\numberwithin{equation}{section}
\newtheorem{thm}{Theorem}
\newtheorem{prop}{Proposition}
\newtheorem{lemma}{Lemma}
\newtheorem{defn}{Definition}

\newcommand{\id}{{1\hspace{-.14cm}1}}
\newcommand{\im}{{\rm im}}
\renewcommand{\ker}{{\rm ker}}
\renewcommand{\aa}{{\mathscr A}}
\newcommand{\pp}{{\mathfrak p}}

\newcommand{\bea}{\begin{eqnarray}}
\newcommand{\eea}{\end{eqnarray}}
\newcommand{\be}{\begin{equation}}
\newcommand{\ee}{\end{equation}}
\newcommand{\bse}{\begin{subequations}}
\newcommand{\ese}{\end{subequations}}
\newcommand{\mb}{\mathbf}

\newcommand{\wt}{\widetilde}

\newcommand{\ol}{\overline}

\newcommand{\eg}{\emph{e.g.}}
\newcommand{\ie}{\emph{i.e.}}
\newcommand{\cf}{\emph{cf.}}

\newcommand{\Z}{{\mathbb Z}}

\newcommand{\C}{{\mathbb C}}
\newcommand{\Q}{{\mathbb Q}}
\newcommand{\PP}{{\mathbb P}}

\newcommand{\cp}{{\mathbb{CP}}}

\newcommand{\bs}{\backslash}
\newcommand{\pd}{\partial}

\newcommand{\CA}{\mathcal{A}}

\newcommand{\CC}{\mathcal{C}}

\newcommand{\CL}{\mathcal{L}}

\newcommand{\CN}{\mathcal{N}}

\newcommand{\CP}{\mathcal{P}}

\newcommand{\CR}{\mathcal{R}}
\newcommand{\CS}{\mathcal{S}}

\newcommand{\CW}{\mathcal{W}}
\newcommand{\CX}{\mathcal{X}}

\title{A Spectral Perspective on Neumann-Zagier}

\author[1]{Tudor Dimofte}
\author[2]{Roland van der Veen}

\affiliation[1]{Institute for Advanced Study, Einstein Dr., Princeton, NJ 08540, USA}
\affiliation[2]{Korteweg-de Vries Institute for Mathematics, University of Amsterdam P.O. Box 94248
1090 GE Amsterdam The Netherlands}

\abstract{We provide a new topological interpretation of the symplectic properties of gluing equations for triangulations of hyperbolic 3-manifolds, first discovered by Neumann and Zagier.
We also extend the symplectic properties to more general gluings of $PGL(2,\C)$ flat connections on the boundaries of 3-manifolds with topological ideal triangulations, proving that gluing is a $K_2$ symplectic reduction of $PGL(2,\C)$ moduli spaces.
Recently, such symplectic properties have been central in constructing quantum $PGL(2,\C)$ invariants of 3-manifolds.
Our methods adapt the spectral network construction of Gaiotto-Moore-Neitzke to relate framed flat $PGL(2,\C)$ connections on the boundary $\CC$ of a 3-manifold to flat $GL(1,\C)$ connections on a double branched cover $\Sigma\to \CC$ of the boundary. Then moduli spaces of both $PGL(2,\C)$ connections on $\CC$ and $GL(1,\C)$ connections on $\Sigma$ gain coordinates labelled by the first homology of $\Sigma$, and inherit symplectic properties from the intersection form on homology.}

\begin{document}

\maketitle

%%%%%%%%%%%%%%%%%%%%%%%%%%%%%%%%%%%%%%%%%%%%%%%%%%%%%%%%%%%%%%%%%%%%%%%%%%

\section{Introduction and motivation}
\label{sec:intro}

Systems with non-abelian gauge symmetry can sometimes be analyzed very effectively using related systems with abelian gauge symmetry. A famous example involves the use of (abelian) Seiberg-Witten theory in four dimensions to compute (non-abelian) Donaldson invariants \cite{Witten-Don}. Donaldson and Seiberg-Witten theory are smoothly connected in physics: they are two limits of the same four-dimensional quantum field theory.

In the present paper we use the same basic philosophy to study moduli spaces of flat (non-abelian) $PGL(2,\C)$ connections in two and three dimensions, by means of closely related --- but conceptually much simpler --- moduli spaces of (abelian) $GL(1,\C)$ flat connections. Again, these two types of moduli spaces occur naturally in the same physical systems. Moduli spaces of flat connections on a two-dimensional surface describe vacua of four-dimensional $\CN=2$ supersymmetric theories of ``class $\CS$'' \cite{Gaiotto-dualities, GMN} (further developed in, \eg\, \cite{GMNII, GMNIV, GMN-spectral}); while moduli spaces of flat connections on a 3-manifold describe vacua of three-dimensional $\CN=2$ theories of ``class $\CR$'' \cite{DGH, Yamazaki-3d,  DG-Sdual, DGG}. Just as 3-manifolds can have 2-dimensional boundaries, the 3d theories of class $\CR$ (labelled by 3-manifolds) describe boundary conditions for the 4d theories of class $\CS$ (labelled by surfaces).

The main question that we address, and hope to shed light on, is the symplectic nature of Thurston's gluing equations \cite{thurston-1980} for ideal triangulations of hyperbolic manifolds, and their generalizations.
The gluing equations for cusped hyperbolic manifolds were first shown to have symplectic properties by Neumann and Zagier \cite{NZ}. The symplectic properties immediately implied a formula for the variation of the volume of a hyperbolic 3-manifold as cusps are deformed. They have since been used to show that A-polynomials of hyperbolic 3-manifolds $M$ are $K_2$-Lagrangian submanifolds in natural $K_2$-symplectic spaces associated to $\pd M$ \cite{Dunfield-Mahler, Champ-hypA}, and that the hyperbolic structures on 3-manifolds can be systematically quantized \cite{hikami-2006, DGLZ, Dimofte-QRS, KashAnd, DG-quantumNZ, Gar-index}.%
\footnote{Ideas about quantization of A-polynomials go back to \cite{gukov-2003, garoufalidis-2004}. Alternative methods of quantization include skein calculus (\eg\ \cite{Frohman-Gelca}) and topological recursion (\eg\ \cite{DF, GS-quant, BorotEynard}), which are expected to be equivalent to triangulation constructions.} %
They also played a crucial role in the construction of 3d $\CN=2$ quantum field theories associated to 3-manifolds \cite{DGG, DGG-index}, which (in principle) provide a categorification of hyperbolic invariants along the lines of \cite{Wfiveknots}.

The gluing equations and their symplectic properties have been generalized in many ways since the work of Neumann and Zagier. Neumann \cite{Neumann-combinatorics} showed that they held for topological ideal triangulations (not necessarily of hyperbolic manifolds). It later became clear that in the topological setting the gluing equations naturally describe a gluing of framed flat $PGL(2,\C)$ connections (\eg, \cite{Zickert-rep, DGG-Kdec}), which include hyperbolic metrics. Symplectic properties of gluing equations were conjectured in \cite{DGV-hybrid} for 3-manifolds partially glued from ideal tetrahedra, with a proposed proof in \cite{Guilloux-PGL}.%
\footnote{This situation is closely related to Bonahon's formula for the deformation of the volume of hyperbolic 3-manifolds with geodesic boundary \cite{Bonahon-vol}.} %
Symplectic properties were also conjectured for generalized $PGL(K)$ gluing equations in \cite{BFG-sl3, DGG-Kdec}, with recent proposed proofs in \cite{GZ-gluing, Guilloux-PGL}.

Unfortunately, so far, all proofs of symplectic properties of gluing equations have involved subtle combinatorics, and have been relatively unintuitive. (The impressive works of \cite{Neumann-combinatorics, GZ-gluing} are testament.) We seek to remedy this situation with an elementary topological construction.

First, we observe that symplectic properties of gluing equations have to do with framed flat connections on \emph{boundaries} rather than interiors of 3-manifolds --- for example, boundaries of ideal tetrahedra and the torus boundary of a fully-glued cusped 3-manifold.
Thus, for the most part, we are dealing with an intrinsically two-dimensional problem.
Then we borrow (and extend) a construction of Gaiotto, Moore, and Neitzke
\cite{GMN-spectral, GMN-snakes} who showed, in the context of 4d $\CN=2$ theories of class $\CS$, that the moduli space of framed flat $GL(K,\C)$ connections with certain singularities on a surface $\CC$ is (roughly) symplectomorphic to a space of flat $GL(1,\C)$ connections on a $K$-fold branched cover $\Sigma\overset\pi\to\CC$\,.\,%
\footnote{Such a spectral-cover construction of local systems also features (independently) in yet-unpublished work of Goncharov and Kontsevich \cite{GK-spectral}, where coordinates on the relevant $GL(1,\C)$ and $GL(K,\C)$ moduli spaces are promoted to fully non-commutative variables.} %
This correspondence was called a non-abelianization map. In \cite{GMN-spectral}, the cover $\Sigma\overset\pi\to\CC$ is a spectral cover, and the non-abelianization map was defined using the data of a related ``spectral network'' on $\CC$ --- hence the title of our paper. For our purposes, we will treat spectral networks (and the non-abelianization maps they induce) as purely topological objects.

The space of $GL(1,\C)$ connections on $\Sigma$ is extremely simple. It has coordinates $x_\gamma\in \C^*$ labelled by cycles $\gamma\in H_1(\Sigma,\Z)$, and a Poisson bracket given by the intersection form in homology, $\{x_\gamma,x_{\gamma'}\} = \langle\gamma,\gamma'\rangle x_\gamma x_{\gamma'}$.
Physically, $H_1(\Sigma)$ is just the electric-magnetic charge lattice of the 4d $\CN=2$ gauge theory labelled by a ``UV curve'' $\CC$ and a Seiberg-Witten curve $\Sigma$.
Via non-abelianization, the space of framed flat $GL(K,\C)$ connections on $\CC$ inherits the coordinates $x_\gamma$ and their simple Poisson bracket, which coincides with the inverse of the Atiyah-Bott symplectic form \cite{AtiyahBott-YM}.
The non-abelianization map can further be modified to provide a symplectomorphism between $PGL(K,\C)$ connections on $\CC$ and a projectivized space of $GL(1,\C)$ connections on $\Sigma$, whose coordinates are labelled by elements of odd homology $\gamma\in H_1^-(\Sigma,\Z) := \ker\,\big[\pi_*:H_1(\Sigma)\to H_1(\CC)\big]$.

In special cases, the $x_\gamma$ coincide with Fock-Goncharov cluster coordinates on spaces of framed flat connections, which complexify Thurston's shear coordinates (or, dually, Penner's length coordinates) in Teichm\"uller theory. As pointed out in \cite{GMN-snakes, HN-FN} and as we discover here, they can also be much more general. Even in the $PGL(2,\C)$ case, the coordinates $x_\gamma$ include complexified Fenchel-Nielsen coordinates on boundaries of 3-manifolds, of the type discussed in \cite{NRS, DG-Sdual, DGV-hybrid, Kabaya-pants, BFG-sl3}.

We claim that symplectic (and in fact $K_2$) properties of gluing equation are an obvious consequence of the topological fact that if a 3-manifold $M$ is glued out of ideal tetrahedra $\{\Delta_i\}_{i=1}^N$ (or, in fact, other 3-manifolds) then the odd homologies of appropriate double-covers of boundaries $\Sigma\overset\pi\to \CC=\pd M$ and $\Sigma_{\Delta_i}\overset\pi\to\CC_\Delta=\pd \Delta$ are related by a lattice symplectic reduction,
\be H_1^-(\Sigma,\Z) \simeq  \oplus_{i=1}^N H_1^-(\Sigma_{\Delta_i},\Z)/\!/G \qquad\text{(modulo $K$-torsion)}\,, \label{HredK} \ee
where $G$ is a certain isotropic subgroup of gluing cycles. 
Setting up all the right structure and definitions needed to understand \eqref{HredK} and the non-abelianization map that relates \eqref{HredK} to a statement about gluing equations is a little tricky. In this paper, we will provide the necessary definitions in the case $K=2$, \ie\ for spaces of framed flat $PGL(2,\C)$ connections. Once the definitions are in place, all proofs are elementary.

We now describe our constructions and main results in a little more detail.

\subsection{Symplectic structures from homology}

\begin{wrapfigure}{r}{2in}
\centering
\includegraphics[width=1.8in]{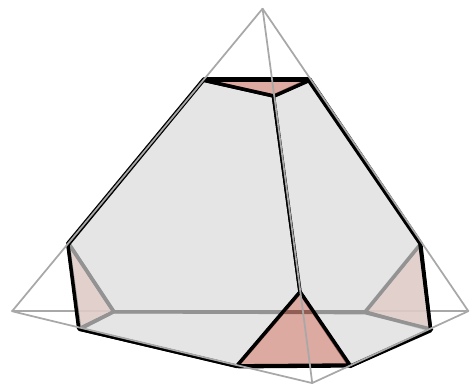}
\caption{Truncated tetrahedron}
\label{fig:trunc}
\end{wrapfigure}

To generalize the notion of an ideal triangulation, we work with a class of ``framed'' 3-manifolds (Section \ref{sec:framed}). They are oriented manifolds $M$ that can be constructed by gluing together pairs of big, hexagonal faces of truncated tetrahedra (Figure \ref{fig:trunc}). We call this a triangulation $\mb t$ of $M$. A framed 3-manifold has its boundary $\CC=\pd M$ split into several parts $\CC=\CC_{\rm big}\cup \CC_{\rm small}\cup \CC_{\rm def}$.
The ``big'' boundary $\CC_{\rm big}$ is tiled by unglued hexagonal tetrahedron faces (we call this tiling a 2d ideal triangulation $\mb t_{2d}$), while the ``small'' boundary $\CC_{\rm small}$ is tiled by the small triangular faces of truncated tetrahedra. If only interiors of some tetrahedron faces are glued there may be also be a ``defect'' boundary $\CC_{\rm def}$, consisting of annuli around unglued edges.

One example of a framed 3-manifold is the tetrahedron $\Delta$ itself. Its big boundary is a 4-holed sphere and it small boundary contains four discs that fill in the holes. Another example is a cusped hyperbolic manifold, such as a knot complement $M=S^3\bs K$. Its small boundary consists of a torus $T^2$ at each cusp, and its big boundary is empty. (An ideal hyperbolic triangulation of $M$ induces a triangulation $\mb t$ as a framed 3-manifold, with $\pd M$ tiled by truncated vertices of tetrahedra.) A closed hyperbolic 3-manifold with a spun triangulation  \cite{thurston-1980} (\cf\ \cite{LTY-spinning}) is a framed 3-manifold whose boundary only contains small spheres, at the vertices of the spun triangulation. Taking either the cusped or closed hyperbolic examples and deleting all (big) edges of the triangulation $\mb t$ produces framed 3-manifolds with with defects, $\CC_{\rm def} \neq \oslash$.
See also Figures \ref{fig:admM}--\ref{fig:admMD} on page \pageref{fig:admM}.

Given a framed 3-manifold $M$, there exists a canonical two-fold branched cover of its boundary $\Sigma\overset\pi\to \CC$. The cover can be constructed by placing a branch point in every face of a triangulation $\mb t_{2d}$ of $\CC_{\rm big}$, and branch cuts along a trivalent graph dual to the triangulation $\mb t_{2d}$ (Figure \ref{fig:branch}, page \pageref{fig:branch}), as well as along the noncontractible cycles of $\CC_{\rm def}$. It turns out that the topological type of the cover is independent of the choice of triangulation used to define it (Lemma \ref{lemma:cover}).

We can define odd homology of the cover as $H_1^-(\Sigma) := \ker\,\big[\pi_*:H_1(\Sigma)\to H_1(\CC)\big]$, working implicitly with $\Z$ coefficients. The odd homology is a nondegenerate symplectic lattice, with skew-symmetric product $\langle *,*\rangle$ given by the usual intersection form. If $M$ is such that $\CC_{\rm small}$ contains only discs and annuli (say) and  $\CC_{\rm def}$ is empty, then a short calculation shows that
\be \label{dim-intro} {\rm rank}\, H_1^-(\Sigma) = 6\,{\rm genus}(\CC)- 6 + 2\,(\text{\# small discs in $\CC_{\rm small}$})\,. \ee
(See \eqref{rankH} for a more general formula.)
We also introduce the twisted homology group $\wt H_1^-(\Sigma)$, a $\Z_2$ extension of $H_1^-(\Sigma)$, which is defined as the (odd) homology of the unit tangent bundle $T_1\Sigma$ with a $\Z_2$ reduction of the fiber class $u$ (Section \ref{sec:twistH}). For the boundary of a framed 3-manifold without defects, there is a natural splitting $\wt H_1^-(\Sigma)\simeq H_1^-(\Sigma)\oplus \Z_2$ (Lemma \ref{lemma:gen}). It is induced by a surjective map
\be \label{path-intro} \wt h:\, \mb P \to\hspace{-.3cm}\to \wt H_1^-(\Sigma) \ee
from a certain group $\mb P$ of paths on $\CC_{\rm small}$ (Section \ref{sec:path}), whose image is a copy of $H_1^-(\Sigma)$.

One may recognize \eqref{dim-intro} as dimension of Teichm\"uller space of a punctured surface $\CC^*$, formed by puncturing $\CC$ once on each small disc. More relevantly for us, it is the complex dimension of the space $\CX[\CC]$ of framed flat $PGL(2,\C)$ connections on $\CC^*$, with unipotent holonomy around the punctures. This space is defined fully in Section \ref{sec:PGL}, following \cite{FG-Teich, DGG-Kdec, DGV-hybrid}; the ``framing'' of a flat connection consists of an extra choice of invariant flag (\ie\ an eigenline of the $PGL(2,\C)$ holonomy) on every component of $\CC_{\rm small}$.
We in fact show (Propositions \ref{prop:coords}--\ref{prop:PB}): \medskip

\noindent{\it Suppose $\pi_1(\CC_{\rm small})$ is abelian. Given any triangulation $\mb t_{\rm 2d}$ of $\CC_{\rm big}$, there is an algebraically open subset $\CP[\CC;\mb t_{\rm 2d}]\subset \CX[\CC]$ and a map
\be x :\, \CP[\CC;\mb t_{\rm 2d}]\times H_1^-(\Sigma) \to \C^*\,, \label{PH-intro} \ee
that's a homomorphism on the second factor (\ie\ $x_{\gamma+\gamma'} = x_\gamma x_{\gamma'}$ for $\gamma,\gamma'\in H_1^-(\Sigma)$) and nondegenerate in the sense that any basis $\{\gamma_i\}_{i=1}^d$ of $H_1^-(\Sigma)$ provides global coordinates $(x_{\gamma_i})\in (\C^*)^{d}$ on $\CP[\CC;\mb t_{\rm 2d}]$.
The map \eqref{PH-intro} may be extended to twisted homology $\wt H_1^-(\Sigma)$, with the convention that the fiber class $u$ maps to $x_u\equiv -1$.
Moreover, there is a non-degenerate holomorphic symplectic structure on $\CP[\CC]$, which agrees with the Atiyah-Bott structure on the space of ordinary (un-framed) flat connections, whose Poisson brackets are}
\be \{x_\gamma,x_\gamma'\} = \langle \gamma,\gamma'\rangle x_\gamma x_{\gamma'}\,. \label{PB-intro} \ee
%\medskip

One proof of these statements follows by labeling both cycles $\gamma\in H_1^-(\Sigma)$ and coordinates on $\CP[\CC]$ by paths $\pp\in \mb P$, using \eqref{path-intro}, and simply computing Poisson brackets. We will follow this approach in Section \ref{sec:PGL}. More fundamentally, the statements follow from a non-abelianization map
\be \label{NA-intro} \Phi[\mb t_{\rm 2d}]:\, \wt \CX_{\rm ab}^-[\Sigma]\big|_{R} \to \CP[\CC;\mb t_{\rm 2d}]\,, \ee
defined using spectral networks in Section \ref{sec:NA}. Here $\wt \CX_{\rm ab}^-[\Sigma]$ is a moduli space of (twisted and projectivized) flat $GL(1,\C)$ connections on $\Sigma$, \ie\ flat $GL(1,\C)$ connections on the unit tangent bundle $T_1\Sigma$ with fiber holonomy $-1$, modulo a certain projective identification. The space $\wt \CX_{\rm ab}^-[\Sigma]$ is parametrized by the holonomies of flat connections along cycles $\gamma\in \wt H_1^-(\Sigma)$, giving an obvious nondegenerate homomorphism
\be x:\, \wt\CX_{\rm ab}^-[\Sigma]\times \wt H_1^-(\Sigma)\to \C^*\,, \ee
and the Atiyah-Bott Poisson bracket among functions $x_\gamma$ is given trivially by \eqref{PB-intro}. We show (Proposition \ref{prop:NA}) that, subject to some mild restrictions `$R$' on the domain, the non-abelianization map $\Phi[\mb t_{\rm 2d}]$ is 1-1 and a symplectomorphism. Therefore, $\CP[\CC,\mb t_{\rm 2d}]$ inherits the coordinates $x_\gamma$ and their simple Poisson bracket.

The holomorphic symplectic form $\omega$ on $\wt \CX_{\rm ab}^-[\Sigma]$ has an avatar $\hat \omega$ in the K-theory group $K_2(F^*)\otimes \Q$, where $F$ is the field of functions on $\wt \CX_{\rm ab}^-[\Sigma]$. It can be written $\hat\omega =\frac12\sum_{ij}(\epsilon^{-1})^{ij}x_i\wedge x_j$ where $\{x_i\}$ are coordinates associated to a basis $\{\gamma_i\}$ of $H_1^-(\Sigma)$ and $\langle \gamma_i,\gamma_j\rangle=:\epsilon_{ij}$. It follows from the fact that $\Phi[\mb t_{\rm 2d}]$ preserves $x_\gamma$ functions that the non-abelianization map is in fact a $K_2$ symplectomorphism, inducing a $K_2$ avatar of the holomorphic symplectic form on $\CP[\CC;\mb t_{\rm 2d}]$. Such avatars (and their motivic versions) were first introduced in \cite{FG-Teich}.

\subsection{Gluing}

Now, suppose that a framed 3-manifold $M'$ is glued together by identifying pairs of hexagonal faces in the big-boundary triangulation $\mb t_{\rm 2d}$ of a framed 3-manifold $M$. For example, $M'$ could be a knot complement, and $M$ could be a disjoint collection of truncated tetrahedra. We assume (largely for simplicity) that neither $M$ nor $M'$ have defects.

We separate the gluing procedure into two steps. First, by gluing only the interiors of pairs of faces of $M$ we form a framed 3-manifold $M_0$ that \emph{does} have defects along some edges of its triangulation. Then we fill in the defects to recover $M'$,
\be M\; \overset{\text{glue interiors of faces}}\leadsto \;M_0\; \overset{\text{fill in edges}}\leadsto \;M'\,. \ee
Each step of the gluing is compatible with the canonical covers $\Sigma,\,\Sigma_0,\,\Sigma'$ of the respective boundaries $\CC,\,\CC_0,\,\CC'$.
It is then an easy exercise to show that, up to 2-torsion, the homology $\wt H_1^-(\Sigma')$ is a lattice symplectic reduction of $\wt H_1^-(\Sigma)$. More precisely (Proposition \ref{prop:glue}), there is an injection of finite (2-torsion) cokernel
\begin{subequations} \label{qg-intro}
\be \wt g : \wt H_1^-(\Sigma_0) \hookrightarrow \wt H_1^-(\Sigma)\,,\ee
which preserves the intersection form,
and there is a distinguished subgroup $\wt G\subset \wt H_1^-(\Sigma_0)$ of ``gluing cycles'' (cycles killed by filling in the defects) such that $\langle \wt G,\wt G\rangle=0$ (\ie\ $\wt G$ is isotropic) fitting into the exact sequence (Lemma \ref{lemma:gen})
\be 0 \to \wt G \overset{\tilde i}\to \wt K \overset{\tilde q}\to \wt H_1^-(\Sigma') \to 0\,, \ee
\end{subequations}
where $\wt K$ is a finite-index subgroup of the complement $\wt K\subset \wt K':= \ker\,\langle \wt G,*\rangle\big|\raisebox{-.1cm}{$\wt H_1^-(\Sigma_0)$}$, with $\wt K'/\wt K=$ 2-torsion. Thus $\wt H_1^-(\Sigma') = \wt K/\wt G \simeq \wt H_1^-(\Sigma_0)/\!/\wt G$ (modulo 2-torsion), and more generally there's a finite-index sublattice $\wt H\subset \wt H_1^-(\Sigma)$ such that
\be \wt H_1^-(\Sigma') = \wt K/\wt G \simeq \wt g(\wt K)/\wt g(\wt G) = \wt H/\!/g(\wt G)\,. \ee
(These maps and equivalences hold for un-twisted homology as well.)

The gluing $M\leadsto M_0\leadsto M'$ also induces a gluing of $PGL(2,\C)$ and $GL(1,\C)$ moduli spaces (Sections \ref{sec:gluePGL} and \ref{sec:glueab})
\be \label{glueP-intro}
\begin{array}{rccc}
g_{PGL(2)} :& \, \CP[\CC;\mb t_{\rm 2d}]\big|_{x_{g(\tilde G)}=1,\, R} &\to\hspace{-.3cm}\to& \CP[\CC';\mb t_{\rm 2d}']\,, \\[.1cm]
g_{GL(1)}: &\, \wt\CX_{\rm ab}^-[\Sigma]\big|_{x_{g(\tilde G)}=1} &\to\hspace{-.3cm}\to&  \wt\CX_{\rm ab}^-[\Sigma']\,,
\end{array}
\ee
where `$R$' denotes some mild (open) extra restrictions. In both cases, the gluing maps turn out to be controlled by the gluing equations
\be \label{eq-intro} \boxed{x_{\tilde g(\gamma)} = x_{\tilde q(\gamma)} \qquad \forall\; \gamma\in \wt K\subset  \wt H_1^-(\Sigma_0)}\,. \ee
Here the RHS contains all coordinates on $\CP[\CC']$ (say), since $\wt q$ is surjective; and they are identified with functions on $\CP[\CC]$ on the LHS. If $\mu\in \wt G\subset \wt K$, then $\wt q(\mu)=0$, so the gluing equations simply say $x_{\tilde g(\mu)}=1$, matching the restriction on the domain in \eqref{glueP-intro}.
Our main results then follow quite quickly: \medskip

\noindent {\bf Theorem \ref{thm:H}} (p. \pageref{thm:H}) {\it The $PGL(2)$ gluing map is the symplectic reduction of a finite quotient
\be g_{PGL(2)}\,:\; (\CP[\CC]\big|_{R}/Z)\big/\!\!\big/(\C^*)^{{\rm rank}(\tilde G)}  = (\CP[\CC]\big|_{R}/Z)\big|_{\tilde g(\tilde G)=1} \big/(\C^*)^{{\rm rank}(\tilde G)} \,\overset{\sim}\to\, \CP[\CC']\,,\ee
where `$R$' is a mild (open) restriction, $Z \simeq \wt H_1^-(\Sigma)/\wt H$ is a finite group (at most 4-torsion), and the group action $(\C^*)^{{\rm rank}(\tilde G)}$ for symplectic reduction is generated by using $x_{\tilde g(\mu)}$ ($\mu\in \tilde G$) as moment maps with respect to the holomorphic symplectic structure. More so, $g_{PGL(2)}$ is a $K_2$ symplectic reduction with respect to $K_2$ avatars $\hat\omega$, $\hat\omega'$ of the symplectic forms on $\CP[\CC]$, $\CP[\CC']$.} \medskip

\noindent An analogue of Theorem \ref{thm:H} holds (rather trivially) for the $GL(1)$ gluing map, and and in fact induces the $PGL(2)$ result via non-abelianization:
\medskip

\noindent {\bf Theorem \ref{thm:NA}} (p. \pageref{thm:NA}) 
{\it
Gluing and non-abelianization maps fit into a commutative diagram
\be \begin{array}{cccl} \label{NAglue-intro}
 \wt \CX_{\rm ab}^-[\Sigma]\big|_{x_{\tilde g(\tilde G)}=1} &\overset{g_{GL(1)}}{\to\hspace{-.3cm}\to} & \wt \CX_{\rm ab}^-[\Sigma'] &= (\wt \CX_{\rm ab}^-[\Sigma]/Z) \big/\!\!\big/ (\C^*)^{{\rm rank}\,\tilde G} \\[.2cm]
 \hspace{-.38in}\Phi[\mb t_{\rm 2d}]\downarrow & & 
  \hspace{-.38in}\Phi[\mb t_{\rm 2d}']\downarrow \\[.2cm]
 \CP[\CC;\mb t_{\rm 2d}]\big|_{x_{\tilde g(\tilde G)}=1} & \overset{g_{PGL(2)}}{\to\hspace{-.3cm}\to} &  \CP[\CC';\mb t_{\rm 2d}'] &
  = (\CP[\CC;\mb t_{\rm 2d}]/Z)\big/\!\!\big/ (\C^*)^{{\rm rank}\,\tilde G}\,,
\end{array}\ee
with appropriate open restrictions on the domains (as in \eqref{NA-intro}, \eqref{glueP-intro}); the vertical maps are 1-1 $K_2$ symplectomorphisms and the horizontal maps are $K_2$ symplectic reduction.
}

\subsection{Example: knot complement}
\label{sec:intro-eg}

To finish the introduction, we illustrate in some detail how the abstract formalism described above applies to a simple example, the ideal triangulation of a knot complement. Let $M'=S^3\bs K$ be the knot complement, viewed as a framed 3-manifold with small torus boundary $\CC'=\CC'_{\rm small}\simeq T^2$, and let $M=\sqcup_{i=1}^N \Delta_i$ be the disjoint union of truncated tetrahedra from which $M'$ is glued. Topologically its boundary is a union of spheres, $\CC = \sqcup_{i=1}^N \CC_{\Delta i}$,\, $\CC_{\Delta i}\simeq S^2$.

\begin{wrapfigure}{r}{2in}
\centering
\includegraphics[width=1.5in]{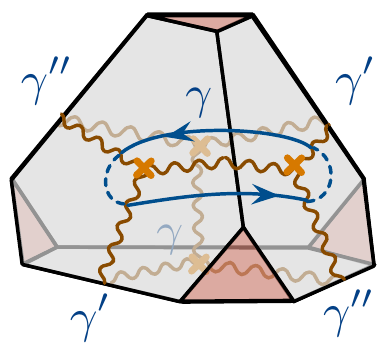}
\caption{Six edge cycles generating $\wt H_1^-(\Sigma_\Delta)$.}
\label{fig:tetH-intro}
\end{wrapfigure}

The canonical cover of the boundary of a tetrahedron $\Sigma_\Delta\overset\pi\to \CC_\Delta$ is branched over four points, one on each face, as illustrated in Figure \ref{fig:tetH-intro}. Thus it has the topology of a torus, $\Sigma_\Delta\simeq T^2$. The first homology (twisted or untwisted) is entirely odd. A convenient basis of generators is given by cycles $\gamma,\gamma',\gamma'$ encircling pairs of branch points, and thus naturally associated to edges of the tetrahedron. We draw $\gamma$ in Figure \ref{fig:tetH-intro}, with the convention that solid paths lie on the top sheet of the cover $\Sigma_\Delta$ and dashed lines lie on the bottom sheet. Such smooth curves also have canonical lifts to the unit tangent bundle $T_1\Sigma_\Delta$ (given by their tangent vectors), and thus represent cycles in twisted homology as well. The cycles $\gamma,\gamma',\gamma''$ are equal on opposite edges and sum to zero. Therefore,
\be \begin{array}{rl}
 H_1^-(\Sigma_\Delta) &= \langle \gamma,\gamma',\gamma''\,|\,\gamma+\gamma'+\gamma''=-1\rangle \simeq \Z^2\,,\\
 \wt H_1^-(\Sigma_\Delta) &= \langle \gamma,\gamma',\gamma'',u\,|\,\gamma+\gamma'+\gamma''=-1,\, 2u=0\rangle \simeq H_1^-(\Sigma_\Delta)\oplus \Z_2\,,
\end{array}
\ee
where $u$ generates the extra fiber class in twisted homology. The intersection product is
\be \langle \gamma,\gamma'\rangle = \langle \gamma',\gamma''\rangle = \langle \gamma'',\gamma\rangle = 1\,,\qquad \langle u,*\rangle = 0\,.\ee

Correspondingly, the space of framed flat $PGL(2)$ connections on the boundary of a tetrahedron (\cf\ \cite{Dimofte-QRS, DGG-Kdec}) is
\be \CP[\CC_\Delta] = \{ x_\gamma,x_{\gamma'},x_{\gamma''}\,|\, x_\gamma x_{\gamma'} x_{\gamma''}=1 \} \simeq (\C^*)^2\,, \ee
with holomorphic symplectic form, $K_2$ form, and Poisson brackets
\be \omega = \frac{dx_\gamma}{x_\gamma}\wedge \frac{dx_{\gamma'}}{x_{\gamma'}}\,,\quad \hat \omega = x_\gamma\wedge x_{\gamma'}\,;\qquad \{\log x_\gamma,\log x_{\gamma'}\} = \{\log x_{\gamma'},\log x_{\gamma''}\}= \{\log x_{\gamma''},\log x_\gamma\}= 1\,.\ee
The coordinates on $\CP[\CC_\Delta]$ (which coincide in this case with Fock-Goncharov coordinates \cite{FG-Teich}, or a complexification of Thurston's shear coordinates) are labelled, as promised, by cycles $\gamma$; and their Poisson brackets are induced by the intersection form on $\wt H_1^-(\Sigma)$. 

To connect with hyperbolic geometry, we note that standard shape parameters are related to $x$ coordinates as
\be z = -x_\gamma  = x_{\gamma+u}\,,\qquad z' = x_{\gamma'+u}\,,\qquad z''=x_{\gamma''+u}\,, \ee
so that $zz'z''=x_u = -1$. Of course, there is a second standard relation among the shapes in hyperbolic geometry, namely $z+z'{}^{-1}-1=0$. This second relation describes framed flat connections that extend from the boundary to the \emph{interior} of a tetrahedron, and cuts out a $K_2$ Lagrangian submanifold $\CL_\Delta\subset \CP[\CC_\Delta]$. We will not need this second relation.

For the full collection $M$ of tetrahedra, we have $H_1^-(\Sigma) = \oplus_{i=1}^N H_1^-(\Sigma_{\Delta i}) \simeq \Z^{2N}$ and $\wt H_1^-(\Sigma) \simeq H_1^-(\Sigma)\oplus \Z_2$ (in twisted homology, fiber classes on disconnected components are all identified, hence a single $\Z_2$ extension). 
Correspondingly, $\CP[\CC] = \prod_{i=1}^N \CP[\CC_{\Delta i}] \simeq (\C^*)^{2N}$. We denote the generators of $\wt H_1^-(\Sigma)$ as $\gamma_i,\gamma_i',\gamma_i''$ and the associated functions on $\CP[\CC]$ as $z_i=-x_{\gamma i},\;z_i'=-x_{\gamma i'}$, etc., with index $i$ for the $i$-th tetrahedron.

\begin{figure}[htb]
\centering
\includegraphics[width=6in]{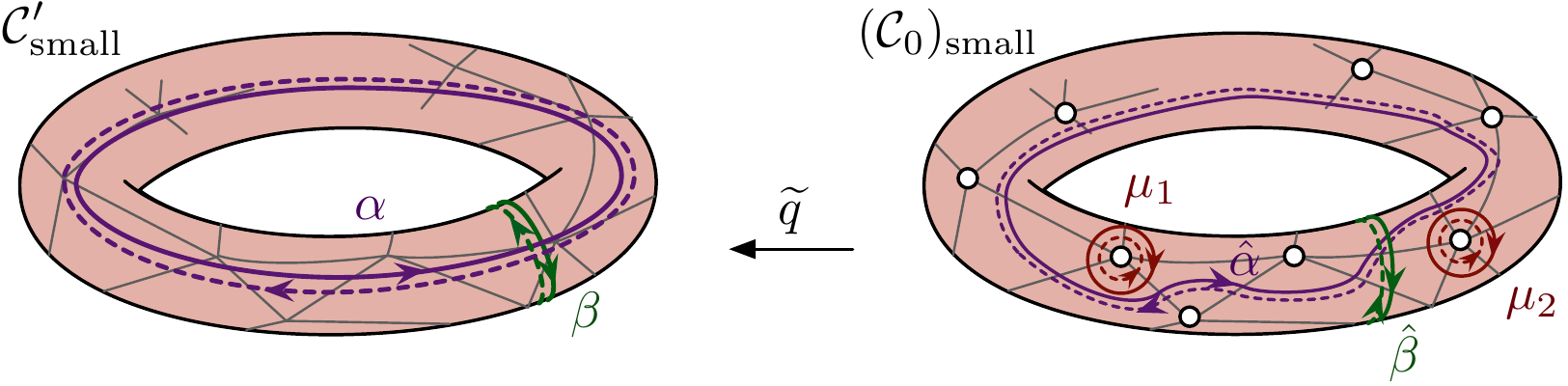}
\caption{The small boundaries of $M'$ and $M_0$, with their (unbranched) canonical covers $\Sigma',\,\Sigma_0$, and curves representing generators of $\wt H_1^-(\Sigma')$ and $\wt H_1^-(\Sigma_0)$.}
\label{fig:intro-cycles}
\end{figure}

For the glued-up knot complement $M'$, the canonical cover of the boundary $\Sigma'\overset\pi\to\CC'\simeq T^2$ is unbranched and disconnected, \ie\ $\Sigma' \simeq T^2\sqcup T^2$. Therefore,
\be \begin{array}{rl}
 H_1^-(\Sigma') &= \Z\langle \alpha,\beta\rangle \simeq \Z^2\,, \\
 \wt H_1^-(\Sigma') &= \Z\langle \alpha,\beta,u\,|\, 2u=0\rangle = H_1^-(\Sigma')\oplus \Z_2\,,
\end{array}
\ee
where $\alpha$ and $\beta$ are ``odd double lifts'' of (say) the meridian and longitude cycles on the torus boundary of $M'$, as on the left of Figure \ref{fig:intro-cycles}. The intersection product is
\be \langle \alpha,\beta \rangle  = 2\,,\ee
with $\langle u,*\rangle\equiv 0$ (always). Correspondingly, there is a complex symplectic space \cite{NZ,gukov-2003}
\be \CP[\CC'] = \{x_\alpha,x_\beta\} \simeq (\C^*)^2\,,\qquad x_\alpha=\ell^2\,,\; x_\beta = m^2\,,\ee
of framed flat $PGL(2)$ connections on the torus boundary of the knot complement, parametrized by the squares of meridian and longitude eigenvalues $x_\alpha = \ell^2$ and $x_\beta=m^2$. Now
\be \omega' = \frac12\frac{dx_\alpha}{x_\alpha}\wedge \frac{dx_\beta}{x_\beta}\,,\quad \hat\omega' = \frac12 x_\alpha\wedge x_\beta\,;\qquad
\{\log x_\alpha,\log x_\beta\} = 2\,.
\ee

In order to relate the boundaries of the tetrahedra $M$ and the knot complement $M'$, we form an intermediate framed 3-manifold $M_0$ with defects: $M_0$ is obtained by gluing tetrahedra (only) along the interiors of their big faces, or equivalently by removing the edges from the triangulation $\mb t$ of $M'$. If there are $N$ tetrahedra there are also $N$ edges in the triangulation. Therefore, the boundary $\CC_0=\pd M_0$ is topologically a surface of genus $N+1$. Its small part $(\CC_0)_{\rm small}$ is a torus with $2N$ holes (formed from the torus $\CC'_{\rm small}$ with endpoints of edges removed), and its defect part $(\CC_0)_{\rm def}$ consists of $N$ annuli. (In this case $(\CC_0)_{\rm big}=\oslash$.)
The canonical cover $\Sigma_0\overset\pi\to \CC_0$ has no branch points, but is characterized by a branch cut along the non-contractible cycle of each annulus in $(\CC_0)_{\rm def}$. It has genus $2N+1$, whence
\be {\rm rank}\,H_1^-(\Sigma_0) = {\rm rank}\,H_1(\Sigma_0)-{\rm rank}\,H_1(\CC_0) = 2(2N+1)-2(N+1) = 2N\,; \ee
thus $H_1^-(\Sigma_0)\simeq \Z^{2N}$ and $\wt H_1^-(\Sigma_0)$ is a $\Z_2$ extension thereof.
The finite-cokernel injection $\wt g$ in \eqref{qg-intro} maps $\wt H_1(\Sigma_0)$ into $\wt H_1(\Sigma)\simeq \Z^{2N}\oplus \Z_2$, preserving the fiber class $u$.

Inside $\wt H_1^-(\Sigma_0)$ lies the subgroup $\wt G$ of gluing cycles, generated by ``odd double lifts'' $\mu_j$ of curves surrounding endpoints of defects $I_j$ on $(\CC_0)_{\rm small}$, as on the right of Figure \ref{fig:intro-cycles}. (The two cycles at the two endpoints of a defect $I_j$ are equivalent, so they may unambiguously be called $\mu_j$.) It is easy to see that the only relation among the $\mu_j$ is $\sum_{j=1}^N \mu_j=0$, whence
\be \wt G \simeq \Z^{N-1}\,.\ee
Moreover, since $\langle \mu_j,\mu_{j'}\rangle=0$, $\wt G$ is an isotropic subgroup of $\wt H_1^-(\Sigma_0)$. The ``complement'' $\wt K:=\ker\,\langle \wt G,*\rangle|_{\tilde H_1^-(\Sigma_0)}$ includes $\wt G$ itself and is a subgroup of rank $N+1$,
\be \wt G \subset \wt K \simeq \Z^{N+1}\oplus\Z_2\,. \ee
The two additional generators of $\wt K$ (besides the fiber class $u$) are represented by ``odd double lifts'' $\hat\alpha,\hat \beta$ of any curves on $(\CC_0)_{\rm small}$ that map to $\alpha,\beta$ once the holes in $(\CC_0)_{\rm small}$ are filled in. The map that fills in the holes is $\wt q$ from (\ref{qg-intro}b); it sends $(\hat\alpha,\hat \beta,\mu_j)$ to $(\alpha,\beta,0)$. Thus
\be \wt H_1^-(\Sigma') \simeq \wt K/\wt G \simeq \wt H_1^-(\Sigma_0)/\!/\wt G\,.\ee

Now, Thurston's gluing equations state that the product of shapes $z_i,z_i',z_i''$ around any edge $I_j$ in the triangulation of $M'$ is trivial, and that the product around longitude and meridian paths (with exponents $\pm 1$) equals the squares $\ell^2$, $m^2$ of longitude and meridian eigenvalues (\cf\ \cite{NZ}). In terms of homology, it turns out that the injection $\wt g$ from (\ref{qg-intro}a) also sends each $\mu_j$ to a sum of tetrahedron cycles $\gamma_i,\gamma_i',\gamma_i''$ around edge $I_j$, and sends $\hat\alpha,\hat \beta$ to sums and differences of parameters around longitude and meridian paths. We will show in Section \ref{sec:Thurs} that, together with some signs coming from the fiber class $u$, Thurston's gluing equations are precisely and succinctly written as
\be x_{\tilde q(\gamma)} = x_{\tilde g(\gamma)}\,,\qquad \forall\; \gamma\in \wt K\,. \ee
When $\gamma=\mu_j$ for some $j$, this is an edge equation; and when $\gamma=\hat\alpha$ or $\hat \beta$ these are cusp equations. It follows (Theorem \ref{thm:H}) that $\CP[\CC']=(\CP[\CC]/Z)/\!/(\C^*)^{N-1} = (\CP[\CC]/Z)|_{\tilde g(\mu_j)=1}/(\C^*)^{N-1}$ is the symplectic reduction of a finite quotient. Explicitly, the quotient action is $x_\gamma\sim (t^\mu)^{\langle \gamma,\tilde g(\mu)\rangle}x_\gamma$ for all $\gamma\in \wt H_1^-(\Sigma)$ and $\mu\in \wt G$, where (say) if $\mu = \sum_{j=1}^{N-1} a_j\mu_j$ in any basis for $\wt G$ then $t^\mu \in \C^*$ is defined as $t_1^{a_1}\cdots t_{N-1}^{a_{N-1}}$ ($t_j\in \C^*$).

There is also a slightly weaker algebraic version of the symplectic properties. Choosing $\{\gamma_i,\gamma_i'\}_{i=1}^N$ as a basis for $H_1^-(\Sigma)$ and $z_i=-x_{\gamma i},\; z_i'=-x_{\gamma i'}$ as the corresponding coordinates on $\CP[\CC]$, the gluing equations take the form
\be  \hspace{.5in} \begin{array}{rll} \ell^2 &= x_\alpha &= \pm z^A z'{}^{A'} \\
 m^2 &= x_\beta &= \pm z^B z'{}^{B'} \\
 1&= x_{\tilde q(\mu_j)} &= \pm z^{C_j} z'{}^{C_j'}\qquad (\forall\; 1\leq j\leq N)
 \end{array}
\ee
for $N$-dimensional integer vectors $A,A',B,B'$ and $N\times N$ matrices $C,C'$. The $(N+1)\times 2N$ matrix
\be  \mb g = \begin{pmatrix} A & A' \\ B & B' \\ C & C' \end{pmatrix} \ee
is nothing but the matrix of the injection $\wt g$ restricted to $\wt K$ (ignoring the fiber class $u$). The fact that $\wt g$ preserves the intersection form implies
\be \mb g\, J_{2N}\, \mb g^{T}  = (2J_2) \oplus 0_{N\times N}\,,\qquad J_{2n}:=\begin{pmatrix} 0 & I_n \\ -I_n & 0 \end{pmatrix}\,; \ee
while because ${\rm rank}( \wt K)=N+1$ we must have ${\rm rank}(\mb g)=N+1$. This was precisely the form of the symplectic properties initially discovered by Neumann and Zagier \cite{NZ}.

\section{Preliminaries}
\label{sec:prelim}

We begin by defining the odd homology of a cover and its twisted version. Then we formally introduce the ``framed'' 3-manifolds that support framed flat connections and construct canonical covers of their boundaries. Section \ref{sec:glue-gen} contains the basic but central result that gluing of framed 3-manifolds induces a lattice symplectic reduction on the odd homology of their boundaries.

\subsection{Odd homology}
\label{sec:oddH}

Let $\CC$ be a closed, oriented surface or disjoint union thereof, and let $\pi:\Sigma\to \CC$ be an oriented double cover, branched over a finite (possibly empty) collection of isolated points $p\in \mathfrak b$. We can think of $\mathfrak b$, the branching locus, as a subset of either $\CC$ or $\Sigma$. We always assume that the branching is simple, with ramification index $2$ at any $p\in \mathfrak b$.

We define the ``odd'' homology of the cover (with $\Z$ coefficients always assumed) as
\be H_\bullet^-(\Sigma) := \ker\,\big[\pi_*:H_\bullet(\Sigma)\to H_\bullet(\CC)\big]\,. \label{def-1} \ee
Note that, since $\Sigma$ is oriented, $H_\bullet(\Sigma)$ is torsion-free, whence $H_\bullet^-(\Sigma)$ is torsion-free as well. It is also convenient to introduce the deck transformation $\sigma:\Sigma\to\Sigma$, an orientation-preserving involution. Letting $\sigma_*:H_\bullet(\Sigma)\overset\sim\to H_\bullet(\Sigma)$ denote the induced automorphism on homology and $P_\pm:=\id\pm \sigma_*$ the associated (quasi-)projections, we may equivalently define
\be H_\bullet^-(\Sigma) := \ker\,P_+ = \ker\, \pi_*\,,\qquad H_\bullet^+(\Sigma):=\ker\, P_-\,. \label{def-2}\ee
The equivalence of \eqref{def-1} and \eqref{def-2} follows from the existence of an injection $\ell^+:H_\bullet(\CC)\hookrightarrow H_\bullet(\Sigma)$ that sends a cycle to the sum of its pre-images%
\footnote{This is sometimes called a ``transfer homomorphism,'' \cf\ \cite[Chapter G.3]{Hatcher}},
obeying
\be \pi_*\circ\ell^+ = 2\id_{H_\bullet(\CC)}\,,\qquad \ell^+\circ \pi_* = P_+\,.\label{l+} \ee
We call $\ell^+$ an ``even double lift.'' A more careful description of $\ell^+$ appears in Appendix \ref{app:odd}, together with some other simple results about odd homology.

The first homology group $H_1(\Sigma)$ has a non-degenerate skew-symmetric intersection form $\langle\;,\;\rangle: \bigwedge^2 H_1(\Sigma,\Z)\to \Z$, which is preserved by $\sigma_*$ (since $\sigma$ is orientation-preserving). It follows (Lemma \ref{lemma:int}, Appendix \ref{app:odd}) that the intersection form is non-degenerate on $H_1^-(\Sigma)$ as well.

It follows from injectivity of $\ell^+$ that
\be {\rm rank}\, H_1^-(\Sigma) = {\rm rank}\,H_1(\Sigma)-{\rm rank}\,H_1(\CC)\,. \label{rank-1}\ee
Combined with the Riemann-Hurwitz formula $\chi(\Sigma)=2\chi(\CC)-\#(\mathfrak b)$ (see \eqref{RH}), we also obtain
\be {\rm rank}\, H_1^-(\Sigma) = - \chi(\CC) + \#(\mathfrak b)\,. \label{rank-2}\ee

\subsection{Twisted homology}
\label{sec:twistH}

We will often need a $\Z_2$ extension of odd homology.
Let $S$ be any closed, oriented, connected surface and let $T_1S$ denote its unit tangent bundle. The homology $H_1(T_1S)$ is an extension%
\footnote{Showing this is a classic application of the Gysin sequence.} %
of $H_1(S)$ by $\Z/(\chi(S)\Z)$, where $\Z/(\chi(S)\Z)$ is generated by the class $u$ of the unit-tangent fiber above any point of $S$, which satisfies $\chi(S)\, u = 0$. We define the twisted homology $\wt H_1(S) := H_1(T_1S)/\Z\langle 2u\rangle$ to be the quotient of $H_1(T_1S)$ by the subgroup generated by $2u$; this makes sense since $\chi(S) = 2-2g(S)$ is always even. Thus $\wt H_1(S)$ is a $\Z_2$ extension,
\be \label{twistseq}
\begin{array}{ccl} 0&\to& \hspace{.1in} \Z_2 \hspace{.2in} \hookrightarrow \wt H_1(S) \overset{p}\to H_1(S) \to 0\,, \\
  &&\hspace{-.2in}=\Z\langle u\rangle/\Z\langle 2u\rangle
  \end{array}
\ee
where $p$ just sets $u\mapsto 0$.
More generally, if $S = \sqcup_{i=1}^n S_i$ has multiple connected components $S_i$, we define $\wt H_1(S) := \big[\oplus_{i=1}^n \wt H_1(S_i)\big]/\Z\langle u_i-u_j\rangle$ by identifying all the fiber classes $u_i\sim u_j$ ($1\leq i,j\leq n$). Then $\wt H_1(S)$ is still a single $\Z_2$ extension of $H_1(S)$, as in \eqref{twistseq}.

\begin{figure}[htb]
\centering
\includegraphics[width=3.5in]{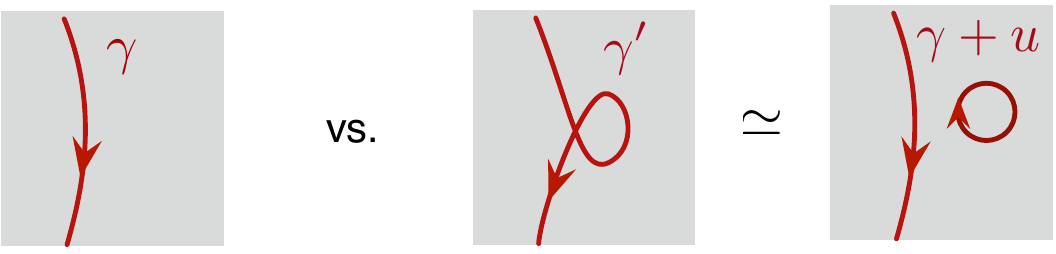}
\caption{Two smooth curves on a surface $S$ representing $[\gamma]$ and $[\gamma']=[\gamma]+u$ in $\wt H_1(S)$.}
\label{fig:u}
\end{figure}

We can represent all classes $[\gamma]\in \wt H_1^-(S)$ by drawing smooth curves $\gamma$ on $S$, and using their unit tangent vectors to lift them canonically to curves on $T_1S$. In particular, the fiber class $u$ is represented by any small contractible loop, as in Figure \ref{fig:u}.

Now suppose $\Sigma\overset\pi\to\CC$ is an oriented branched double cover as above. The deck transformation $\sigma$ has nonvanishing derivative $d\sigma$, and thus extends to an involution on $T_1\Sigma$, which in turn induces an involution $\sigma_*$ on twisted homology $\wt H_1(\Sigma)$. Indeed, since the deck transformation is orientation-preserving, we have $\sigma_* u= u$. Letting $P_\pm = \id\pm \sigma_*$, we define
\be \wt H_1^\pm(\Sigma) := \ker\, P_\mp\big|\raisebox{-.2cm}{$\wt H_1(\Sigma)$}\,.\ee
Note that $P_+u=P_-u=0$, so $u\in \wt H_1^\pm(\Sigma)$. Indeed, $\wt H_1^-(\Sigma)$ (say) is just the preimage of $H_1^-(\Sigma)\subset H_1(\Sigma)$ under the map $p$ in \eqref{twistseq}, so $\wt H_1^-(\Sigma)$ is a $\Z_2$ extension of $H_1^-(\Sigma)$,
\be \label{twistseq-}
 0 \to \Z_2 \hookrightarrow \wt H_1^-(\Sigma) \overset p\to H_1^-(\Sigma)\to 0\,. \ee

The intersection product on $\wt H_1^-(\Sigma)$ is simply pulled back from $H_1^-(\Sigma)$, \ie\  $\langle \gamma,\gamma'\rangle_{\tilde H_1^-} := \langle p(\gamma),p(\gamma')\rangle_{H_1^-}$.  Since the product on $H_1^-(\Sigma)$ is nondegenerate, it follows that $\langle \gamma,\gamma'\rangle=0$ $\forall\, \gamma'\in \wt H_1^-(\Sigma)$ if and only if $\gamma=0$ or $\gamma=u$.

\subsection{Boundaries of framed 3-manifolds}
\label{sec:framed}

Framed 3-manifolds, as introduced in \cite{DGG-Kdec,DGV-hybrid}, are a class of 3-manifolds with extra structure on which it is natural to define moduli spaces of framed flat connections. They include both knot complements and single tetrahedra. In the present section, we look at framed 3-manifolds topologically, construct canonical branched covers of their boundaries, and study how the homology of these covers behaves under gluing. We use the most liberal possible definition:

\begin{defn} \label{def:framed}
A \emph{framed 3-manifold} $M$ is an oriented 3-manifold with boundary (or disjoint union of such), together with a splitting of its boundary $\pd M=\pd M_{\rm big}\cup \pd M_{\rm small}$ into ``big'' and ``small'' parts, such that

a) $\pd M_{\rm big}$ consists of surfaces of any genus with at least one hole and negative Euler character;

b) $\pd M_{\rm small}$ consists of surfaces of any genus with any number of holes;

c) the two parts of the boundary attach along circles, $\pd M_{\rm big}\cap \pd M_{\rm small}=\sqcup_i S^1_i$\,.
\end{defn}

\noindent Often (as in \cite{DGG-Kdec,DGV-hybrid}) it is useful to impose additional restrictions on the small boundary. For example, requiring $\pd M_{\rm small}$ to have abelian fundamental group ensures that framing data does not restrict the choice of a flat connection on $\pd M$; and requiring $\pd M_{\rm small}$ to admit a Euclidean structure ensures the existence of additional data needed for quantization (\eg\ combinatorial flattenings of \cite{Neumann-combinatorics}, \cf\ \cite{Dimofte-QRS, DG-quantumNZ}). For the moment, we don't need these restrictions.

The basic example of a framed 3-manifold is a truncated tetrahedron (Figure \ref{fig:trunc}), whose big boundary is a 4-holed sphere (formed from four big hexagonal faces) and whose small boundary is four small triangles (at the truncated vertices).
We define a \emph{triangulation} (or tiling) $\mb t$ of a framed 3-manifold to be a decomposition into truncated tetrahedra that are glued only along big hexagonal faces.%
\footnote{Such a triangulation always exists, though we don't need this fact here. For example, to construct a triangulation for an arbitrary framed 3-manifold $M$, one can first take two copies $M,\ol M$ with opposite orientation and glue them along the big boundary to form $N=M\cup_{\rm big} \ol M$. The new $N$ has only closed small boundary components. It admits a topological ideal triangulation (equivalent to a tiling by truncated tetrahedra) by classic theorems of Matveev \cite{Matveev-book} (\cf\ \cite[Prop 1.2]{Tillmann-normal}). Then one can refine the triangulation of $N$ so that $\pd M_{\rm big}\subset N$ is realized by faces of tetrahedra, and cut along $\pd M_{\rm big}$ to recover a triangulation of $M$.} %
Thus, $\mb t$ induces

\noindent \;\;- a tiling of $\pd M_{\rm big}$ by big hexagons, equivalent to an ideal 2d triangulation $\mb t_{2d}$; and

\noindent\;\;- a tiling of $\pd M_{\rm small}$ by small triangles.

\begin{figure}[htb]
\centering
\includegraphics[width=5in]{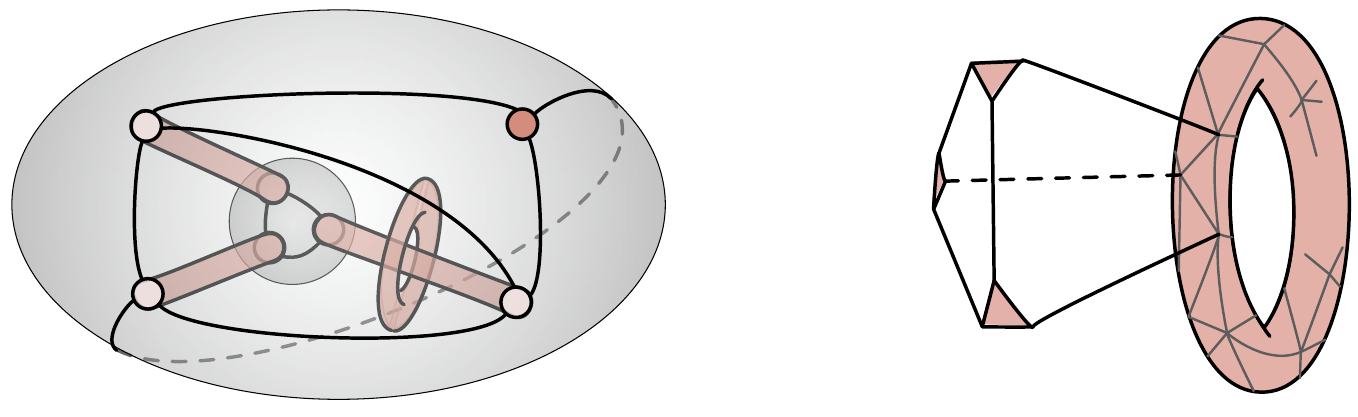}
\caption{Left: a typical framed 3-manifold, with two big boundaries (a 3-holed sphere and a 4-holed sphere) and several small tori, annuli, and discs. A 2d triangulation of the big boundaries is shown. Right: the triangulation of a small torus boundary induced from a 3d triangulation $\mb t$.}
\label{fig:admM}
\end{figure}

For example, any knot complement $M$ with a standard ideal triangulation can be viewed as a framed 3-manifold with a single small torus boundary and a tiling $\mb t$ where all big hexagonal faces have been glued pairwise. The tiling is simply obtained by truncating the ideal triangulation. This is illustrated in Figure \ref{fig:admM}.

We also introduce
\begin{defn} \label{def:framedD}
A \emph{framed 3-manifold with defects} $M$ is an oriented 3-manifold with boundary, together with a splitting $\pd M = \pd M_{\rm big}\cup \pd M_{\rm small}\cup \pd M_{\rm def}$ into three parts, such that $\pd M_{\rm big}$ and $\pd M_{\rm small}$ are as in Def. \ref{def:framed}, while $\pd M_{\rm def}$ is a collection of annuli that attach only to holes on the small boundary:  $\pd M_{\rm big}\cap \pd M_{\rm small} = \sqcup_i S^1_i$, $\pd M_{\rm small}\cap \pd M_{\rm def} = \sqcup_i S^1_j$, $\pd M_{\rm def}\cap\pd M_{\rm big} = \oslash$.
\end{defn}

\begin{figure}[htb]
\centering
\includegraphics[width=5in]{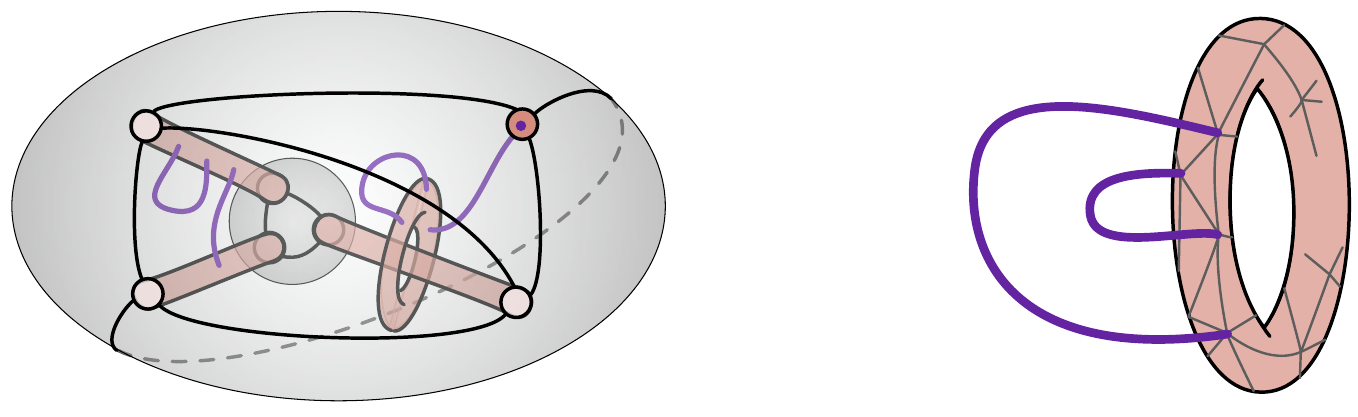}
\caption{Including defects (extra annular boundary components) that begin and end on $\pd M_{\rm small}.$ On the right, the defects are compatible with a 3d triangulation: they begin and end at vertices of the tiling of $\pd M_{\rm small}$.}
\label{fig:admMD}
\end{figure}

\begin{wrapfigure}{r}{2.1in}
\centering
\vspace{-.4cm}
\includegraphics[width=1.9in]{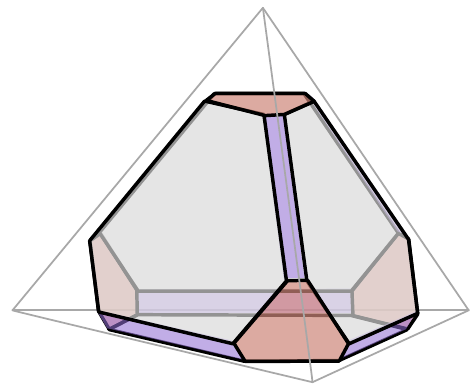}
\caption{Doubly truncated tetrahedron}
\label{fig:ttrunc}
\end{wrapfigure}

We can think about forming a framed 3-manifold with defects in two different ways. We can start with a framed 3-manifold $M$ without defects and a collection of open curves $\{I_i\}$ in the interior of $M$, whose endpoints lie on the small boundary. Excising neighborhoods of the $I_i$ produces a framed 3-manifold with defects $M'=M\bs \cup_i I_i$, as in Figure \ref{fig:admMD}. We say that a triangulation $\mb t$ of $M$ is compatible with the defects if the $I_i$ lie along edges of the triangulation. Then $M'$ inherits a tiling $\mb t'$ by tetrahedra with truncated vertices and some truncated edges (those along the $I_i$), as in Figure \ref{fig:ttrunc}.

Alternatively, we may start with a framed 3-manifold $M$ and a triangulation $\mb t$, and identify the \emph{interiors} (but not the edges) of some pairs of faces of the truncated tetrahedra in $M$. In general this again leads to a framed 3-manifold $M'$ with defects along the new internal edges created by the gluing. It inherits a triangulation $\mb t'$.

In the extreme case, we may build a framed 3-manifold $M'$ with defects entirely from doubly truncated tetrahedra, by gluing their big hexagonal faces together in pairs. Then a defect lies along every edge in the triangulation of $M'$.

\subsection{The canonical cover}
\label{sec:cc}

The boundary $\CC = \pd M$ of any framed 3-manifold, with or without defects, admits a canonical double cover $\Sigma\overset{\pi}\to \CC$.

\begin{figure}[htb]
\centering
\includegraphics[width=4.4in]{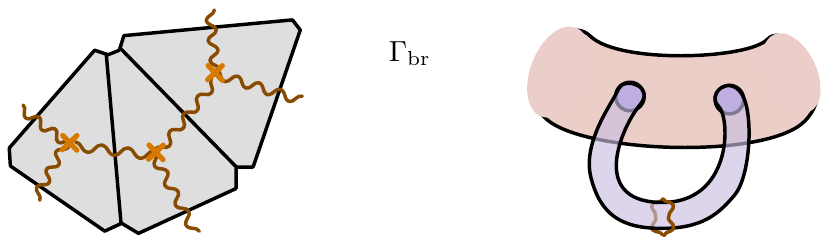}
\caption{The branching graph $\Gamma_{\rm br}$, dual to a big-boundary triangulation (left) and containing closed curves around each defect (right).}
\label{fig:branch}
\end{figure}

To construct it, we first choose a triangulation $\mb t$ of $M$. Let $\mb t_{\rm 2d}$ be the induced tiling of the big boundary $\pd M_{\rm big}$ by hexagons, thought of as a 2d ideal triangulation. Let $\Gamma_{\rm br}$ be the union of a trivalent graph dual to the triangulation $\mb t_{\rm 2d}$ and a closed curve around the girth of every defect of $M$, as in Figure \ref{fig:branch}.
We construct the cover $\Sigma\overset\pi\to \CC=\pd M$ by using $\Gamma_{\rm br}$ as a branching graph --- putting branch cuts on the edges of $\Gamma_{\rm br}$ and branch points at the vertices of $\Gamma_{\rm br}$.
In other words, we take two identical copies $\Sigma^+,\Sigma^-$ of $\CC\bs\Gamma_{\rm br}$, glue $\Sigma^+$ to $\Sigma^-$ (and vice versa) along the edges of $\Gamma_{\rm br}$, and add in one copy of the vertices $\mathfrak b$ of $\Gamma_{\rm br}$ to complete the space to a closed surface $\Sigma$:
\be \Sigma = \Sigma^+\cup\Sigma^-\cup\mathfrak b = (\CC\bs\Gamma_{\rm br})\cup_{\text{edges}(\Gamma_{\rm br})}(\CC\bs\Gamma_{\rm br}) \cup \mathfrak b\qquad \text{(with $\mathfrak b=$vertices$(\Gamma_{\rm br})$)}\,.
\ee

\begin{lemma} \label{lemma:cover}
This definition of the canonical cover $\Sigma\overset{\pi}\to \CC$ does not depend on a choice of triangulation $\mb t$.
\end{lemma}

\begin{figure}[htb]
\centering
\includegraphics[width=5in]{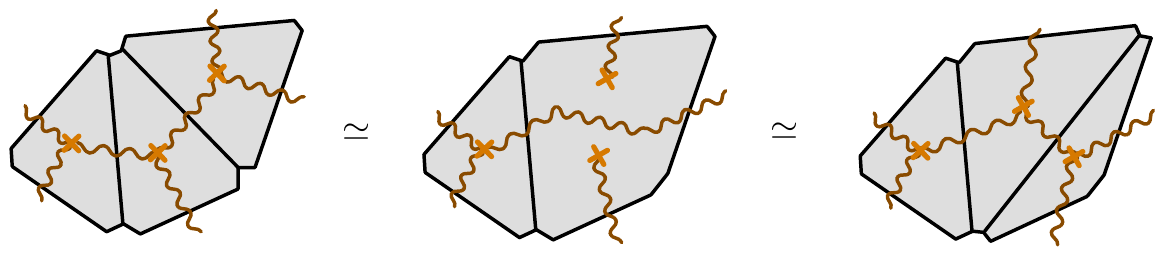}
\caption{Invariance of the canonical cover $\Sigma\overset\pi\to\CC$ under a flip of the big-boundary triangulation.}
\label{fig:branch-inv}
\end{figure}

\noindent\emph{Proof.} The only possible dependence on $\mb t$ comes from the induced ideal triangulation $\mb t_{\rm 2d}$ on the big boundary. Any two triangulations of the big boundary are related by a sequence of flips, \emph{a.k.a.} 2-2 Pachner moves.  Each flip, moreover, changes the cover $\Sigma\overset{\pi}\to \CC$ locally in a way that does not modify its topological type (\ie, corresponds to a homeomorphism of $\Sigma$ and a homotopy of $\pi$). To see this, we simply reconnect branch cuts as in Figure \ref{fig:branch-inv}. (In general, the topological type of a cover $\Sigma\overset\pi\to\CC$ defined from a branching graph only depends on the relative homology class in $H_1(\CC,\mathfrak b;\Z_2)$ induced by the graph.) \;$\square$ \medskip

Another simple but fundamental result is

\begin{lemma} \label{lemma:rank-dim}
Let $M$ be a framed 3-manifold without defects, with boundary $\CC=\pd M$ such that $\pi_1(\CC_{\rm small})$ is abelian. Let $\CC_i$ denote the connected components of $\CC$, with $g_i$ the genus of $\CC_i$ and $n_i$ the number of disc components of $(\CC_i)_{\rm small}$. Let $\Sigma_i\overset\pi\to\CC_i$ be the canonical cover over $\CC_i$. Then
\be \label{rankH}
{\rm rank}\, H_1^-(\Sigma_i) = \begin{cases} 0 & \text{$\CC_i$ is a small sphere} \\
2 & \text{$\CC_i$ is a small torus} \\
 6g_i-6 + 2n_i & \text{otherwise}\,.\end{cases}
\ee
\end{lemma}

\noindent {\it Proof.} This is trivial for small spheres or tori (which themselves are closed, disjoint components of $\CC$). Otherwise, $\CC_i$ contains big boundary with holes connected to small discs and small annuli. We calculate $\chi(\CC_i) = (\text{\# faces of $\mb t_{\rm 2d}$}) + n_i - (\text{\# big edges of $\mb t_{\rm 2d}$}) = -\tfrac12 (\text{\# faces of $\mb t_{\rm 2d}$}) + n_i=-\tfrac12\#(\mathfrak b)+n_i$. On the other hand, $\chi(\CC_i)=2-2g_i$, so $\#(\mathfrak b) = 4g_i-4+2n_i$\,; then the result follows from \eqref{rank-2}. \; $\square$

\subsection{Gluing of canonical covers}
\label{sec:glue-gen}

As a final basic result, we now describe how the odd homology of canonical covers of boundaries changes as framed 3-manifolds are glued.

\begin{figure}[htb]
\centering
\includegraphics[width=5.5in]{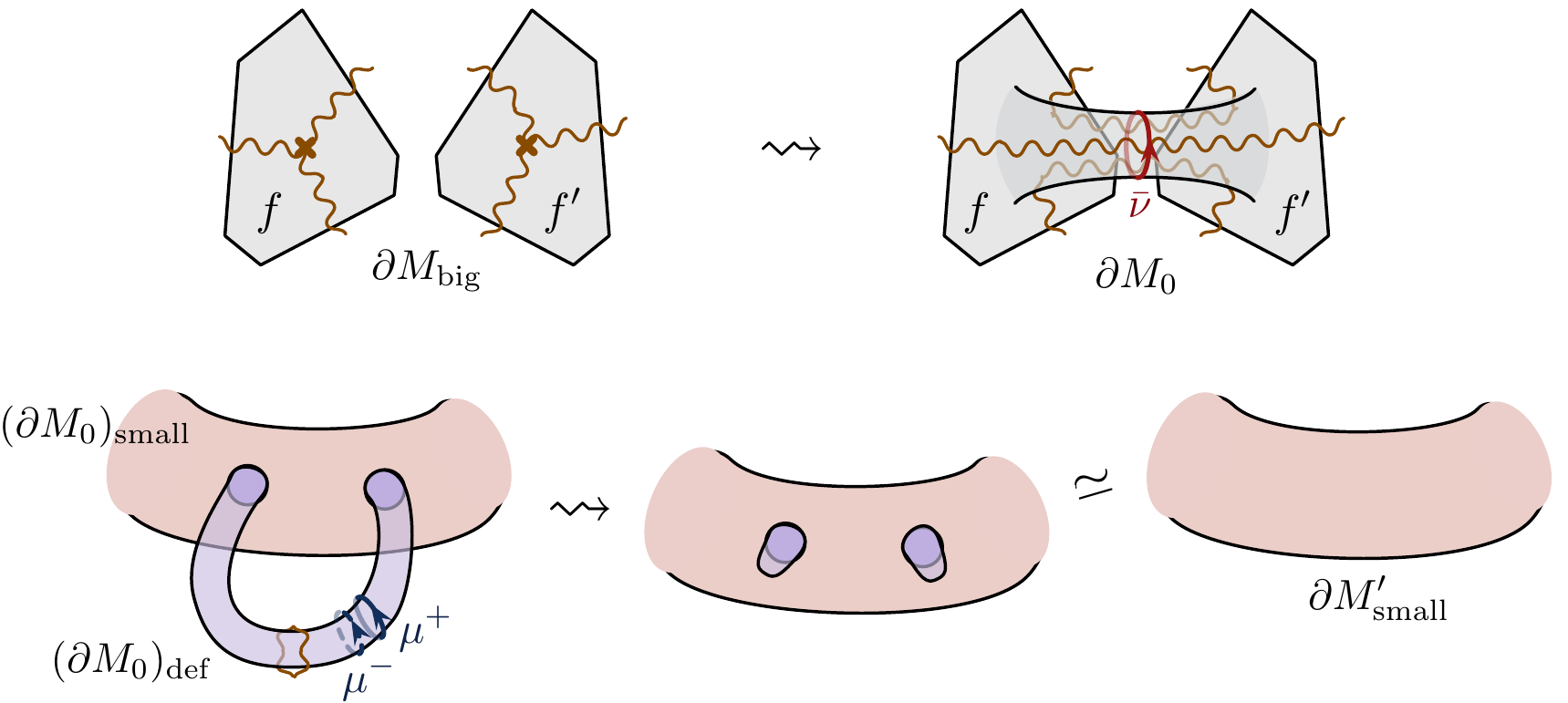}
\caption{Gluing $M\leadsto M_0 \leadsto M'$ in two steps: first, identifying interiors of faces on the big boundary of $M$, then filling in defects of $M_0$.}
\label{fig:SxSg}
\end{figure}

Consider the following setup, used throughout the remainder of the paper.
Let $M$ be a framed 3-manifold with triangulation $\mb t$. ($M$ may have defects and/or multiple components.) Choose several pairs of faces $(f_i,f_i')$ of the induced triangulation on the big boundary $\pd M_{\rm big}$ to glue together, and identify their \emph{interiors} to form a manifold $M_0$\,:
\be M_0 = M|_{f_i^\circ \sim -f_i'{}^\circ}\,. \label{MM1} \ee
In general, $M_0$ will be a framed 3-manifold with defects, and an induced triangulation $\mb t_0$. Let $\{I_j\}$ denote the new defects in $M_0$, corresponding to new internal edges of $\mb t_0$ that were not present in $\mb t$. We may then fill in the defects $I_j$ to form a new 3-manifold $M'$,
\be M' = M_0 \cup(\sqcup_j I_j)\,. \label{MM2} \ee
We illustrate this gluing process in Figure \ref{fig:SxSg}.

A standard example of such gluing involves a disjoint union of truncated tetrahedra $M = \sqcup_{i=1}^N \Delta_i$, glued together to form (say) a knot complement $M'$. In this case $M_0$ is obtained by gluing doubly truncated tetrahedra; the boundary $\pd M_0$ has genus $N+1$ and can be thought of as the small torus boundary of $M'$ with $N$ additional defects drilled out. We will come back to this example in Sections \ref{sec:N}, \ref{sec:Thurs}.

Now, let us denote the boundaries of the three 3-manifolds as $\CC := \pd M$, $\CC_0:=\pd M_0$ and $\CC' := \pd M'$. They have canonical double covers
\be \Sigma\overset\pi\to\CC\,,\qquad \Sigma_0\overset\pi\to\CC_0\,,\qquad
\Sigma'\overset\pi\to\CC'\,.\ee
To relate the homologies of these covers, we introduce oriented closed curves $\ol\mu_j \subset \CC_0$ winding around the girth of each defect $I_j$ (parallel to the branch lines, as in Figure \ref{fig:SxSg}). Denote by $\mu_j^+,\mu_j^-$ the lifts of $\ol\mu_i$ to the two sheets of $\Sigma_0$, and define subgroups
\be \label{defG}
\begin{array}{ll}
G_0 \subset H_1(\Sigma_0) &\quad \text{subgroup generated by all the $\mu_j^\pm$} \\
G' := \ker\,P_+|_{G_0} \subset H_1^-(\Sigma_0)  \\
G := \im\,P_-|_{G_0} \subset H_1^-(\Sigma_0) & \quad\text{subgroup generated by $\mu_j:=\mu_j^+-\mu_j^-$}
\end{array}
\ee
In general, $G'/G$ may be a nontrivial 2-torsion group.

\begin{prop} \label{prop:glue}

Let $K':=\ker\,\langle G',*\rangle|_{H_1^-(\Sigma_0)}$.
Then:

a) There is an injection $g: H_1^-(\Sigma_0) \hookrightarrow H_1^-(\Sigma)$ with finite (2-torsion) cokernel that preserves the intersection form.

b) There is an exact sequence
\be  0 \to  G' \overset{i}\to K' \overset{q}\to H_1^-(\Sigma') \to 0 \label{SESG} \ee
providing an isomorphism $H_1^-(\Sigma') \simeq H_1^-(\Sigma_0)/\!/G' = K'/G'$, preserving the intersection form.

c) Suppose $K\subset H_1^-(\Sigma_0)$ is a subgroup such $G\subset K\subset K'$ and $K'/G'\simeq K/G$.
Then there is a finite-index subgroup $H \subset H_1^-(\Sigma)$ containing $g(G)$ such that
\be H_1^-(\Sigma') \simeq H/\!/g(G) = g(K)/g(G)\,. \label{Hiso} \ee
The quotient $H_1^-(\Sigma)/H$ contains (at most) 2-torsion and 4-torsion.

\end{prop}

\emph{Proof.} For part (a), we observe that relating $\Sigma \leadsto \Sigma_0$ involves the inverse of the basic cut-and-glue operation of Appendix \ref{app:cutglue}. Namely, from the perspective of the boundary $\CC$, gluing the interiors of two faces $(f_i,f_i')$ is equivalent to excising discs in the interiors of these faces (containing branch points) and subsequently identifying the circular boundaries of holes that are created. The gluing can be reversed by cutting the boundary $\CC_0$ along curves $\ol\nu_i$ (Figure \ref{fig:SxSg}), and filling in the resulting holes with discs. This gluing/ungluing lifts to $\Sigma$ and $\Sigma_0$. Applying Lemma \ref{lemma:S-} (p. \pageref{lemma:S-}) to the pair $(\Sigma_0,\Sigma)$, and noting that the group $G_\nu\subset H_1(\Sigma_0)$ generated by lifts $\nu_i$ of the $\ol \nu_i$ is even --- so that $G_\nu^-:=\ker\,P_+|_{G_\nu}=0$ and $\ker\,\langle G_\nu^-,*\rangle|_{H_1^-(\Sigma_0)}=0$ --- we obtain the injection $H_1^-(\Sigma_0)/\!/G_\nu^- =H_1^-(\Sigma_0) \hookrightarrow H_1^-(\Sigma)$. The cokernel is at most 2-torsion and may be nontrivial because the lifts $\nu_i$ are connected (and branch points are destroyed/created in the gluing/ungluing).

For part (b), we observe that relating boundaries $\CC_0\leadsto \CC$ and $\Sigma_0\leadsto \Sigma$ is again an application of a basic cut-and-glue. This time, filling in defects means from the perspective of the boundary to cut $\CC_0$ along the $\ol\mu_i$ and cap off holes with discs. The operation lifts to $\Sigma_0$. Applying Lemma \ref{lemma:S-} and noting that the $\ol\mu_i$ have two distinct lifts (so no branch points are created/destroyed) yields the sequence \eqref{SESG}, with identically vanishing homology.

For part (c), we construct $H$ as follows. Observe that $H_1^-(\Sigma_0)/K'$ is torsion-free (because $\langle G',n\alpha\rangle = 0$ $\Leftrightarrow$ $\langle G',\alpha\rangle = 0$), so we can decompose $H_1^-(\Sigma_0) = K'\oplus K'{}^\perp$ for some (non-unique!) sublattice $K'{}^\perp \simeq H_1^-(\Sigma_0)/K'$. Let $\wt H := K\oplus K'{}^\perp$, and note that $H_1^-(\Sigma_0)/\wt H \simeq K'/K \simeq G'/G$ is 2-torsion. Then we set $H:=g(\wt H) \subset H_1^-(\Sigma)$. Since the cokernel of the injection $g$ is 2-torsion, $H_1^-(\Sigma)/H$ contains at most 2- and 4-torsion.
The isomorphisms \eqref{Hiso} follow from $H_1^-(\Sigma)\simeq K'/G'\simeq K/G$ together with the fact that $g$ is injective (and preserves the intersection form). \;$\square$ \medskip

\noindent {\bf Remark.} We emphasize, following Appendix \ref{app:cutglue}, that $g$ can explicitly be defined by representing a cycle $[\gamma]\in H_1^-(\Sigma_0)$ by a curve (or curves) $\gamma$ away from the gluing regions, then including $\gamma$ as curves in $\Sigma$ and passing to homology. Similarly, $q$ is defined by representing $[\gamma]\in K'\subset H_1^-(\Sigma_0)$ by curve(s) $\gamma$ disjoint from defects, then filling in the defects and passing to homology. \medskip

Gluing can also be extended in a straightforward manner to unit tangent bundles and twisted homology. Namely, the gluing of surfaces $\Sigma\leadsto \Sigma_0$ and $\Sigma_0\leadsto \Sigma'$ extends uniquely to unit tangent bundles ($T_1\Sigma\leadsto T_1\Sigma_0$, $T_1\Sigma_0\leadsto T_1\Sigma'$
) by performing all cut/glue operations smoothly, \ie\ cutting and gluing discs (etc.) along circles without introducing any kinks. Then, for example, the (even) classes $[\nu_i]-u\in \wt H_1(\Sigma_0)$ represented by non-intersecting smooth curves $\nu_i\subset \Sigma_0$ become trivial when included in $\wt H_1(\Sigma)$; and the classes $[\mu_j^\pm]-u\in \wt H_1(\Sigma)$ represented by the $\mu_j^\pm$ become trivial in $\wt H_1(\Sigma')$ once defects are filled in. Let us therefore define
\be \label{defGtw}
\begin{array}{ll@{\quad}l}
\wt G_0 &\subset \wt H_1(\Sigma_0) & \text{subgroup generated by $\mu_j^\pm-u$} \\
\wt G' &:= \ker\,P_+|_{\tilde G_0} \subset \wt H_1^-(\Sigma_0) \\
\wt G  &:= \im\,P_-|_{\tilde G_0} \subset \wt H_1^-(\Sigma_0) & \text{subgroup generated by $\mu_j:= \mu_j^+-\mu_j^-$}\,,
\end{array}
\ee
where we now use `$\mu_j$' and `$\mu_j^\pm$' to denote the classes of the canonical lifts of the gluing curves to $T_1\Sigma_0$.
Let $\wt K' = \ker\,\langle \wt G',*\rangle|\raisebox{-.1cm}{$\wt H_1^-(\Sigma_0)$}$.

We may define a lift $\wt g:\wt H_1^-(\Sigma_0)\to \wt H_1^-(\Sigma)$ of the map $g$ by representing cycles $[\gamma]\in \wt H_1^-(\Sigma_0)$ by curves $\gamma\subset T_1\Sigma_0$ away from the gluing regions, including them as curves in $T_1\Sigma$ and passing to homology. We define a lift $\wt q:K'\to \wt H_1^-(\Sigma')$ of the map $q$ the same way, by representing cycles $[\gamma]\in K'$ by curves $\gamma\subset T_1\Sigma_0$ disjoint from the defects, filling in the defects, and passing to homology. Observe that $\wt g(u)=u$ and $\wt q(u)=u$ (\ie\ these maps preserve the fiber class).  Then \medskip

\noindent {\bf Proposition \ref{prop:glue}'}

{\it a) The map $\wt g$ is an injection of finite cokernel, and $\wt q$ is a surjection with kernel $\wt G'$. They fit into commutative diagrams 
\be \label{gluetw}
\begin{array}{ccccccccc}
 &&&&0&  & 0 \\
 &&&& \downarrow && \downarrow \\
 &&0&&\Z_2& = & \Z_2 \\
 &&\downarrow&& \downarrow && \downarrow \\
 0 & \to & \wt G' & \overset{\tilde i}\hookrightarrow & \wt K' & \overset{\tilde q}\to & \wt H_1^-(\Sigma',\Z) & \to & 0 \\
 &&\hspace{-.27cm}p \downarrow&& \hspace{-.27cm}p \downarrow && \hspace{-.27cm}p \downarrow \\
 0 & \to & G' & \overset{i}\hookrightarrow & K' & \overset{q}\to & H_1^-(\Sigma',\Z) & \to & 0\,, \\
 &&\downarrow&&\downarrow && \downarrow \\
 &&0&& 0 && 0
\end{array}
\hspace{.7in}
\begin{array}{ccc}
 0&&0 \\
 \downarrow &&\downarrow \\
 \Z_2 &=&\Z_2 \\
 \downarrow && \downarrow \\
 \wt H_1^-(\Sigma_0) & \overset{\tilde g}\hookrightarrow & \wt H_1^-(\Sigma) \\
 \hspace{-.27cm}p \downarrow && \hspace{-.27cm}p \downarrow \\
 H_1^-(\Sigma_0) & \overset g\hookrightarrow & H_1^-(\Sigma)\,. \\
 \downarrow && \downarrow \\
 0 &&0
\end{array}
\ee
with all columns and all rows on the left exact.

b) Given any $\wt K\subset \wt H_1^-(\Sigma_0)$ such that $\wt G\subset \wt K\subset \wt K'$ and $\wt K/\wt G=\wt K'/\wt G'$, the quotient $\wt K'/\wt K$ is at most 2-torsion, and there exists $\wt H\subset \wt H_1^-(\Sigma)$ such that $\wt H_1^-(\Sigma)/\wt H$ is finite (at most 2- and 4-torsion), and $\wt H_1^-(\Sigma') = \wt H_1^-(\Sigma)/\!/g(\wt G) = g(\wt K)/g(\wt G)$. If $\wt K$ extends a group $K$ as in Prop \ref{prop:glue}c, then $\wt H$ can be chosen to extend $H$ appearing there. \medskip
}

\noindent \emph{Proof.}  Part (a) is a straightforward diagram chase. On the RHS, the columns are exact by definition of twisted homology. The top square commutes because $\wt g(u)=u$, and the bottom square commutes by comparing definitions of $\wt g$ and $g$. Injectivity of $g$ and $\wt g(u)=u$ imply injectivity of $\wt g$.

On the LHS, the right column is exact by definition, and the middle column is exact because $\wt K'=\langle \wt G',*\rangle|_{\tilde H_1^-(\Sigma_0)}$ contains $u$, hence is a $\Z_2$ extension of $K'=\langle G',*\rangle|_{H_1^-(\Sigma_0)}$. The top square commutes because $\wt q(u)=u$, and $\wt q$ is surjective because $q$ and the $p$'s are surjective. The (bottom-)right square commutes by comparing definitions of $\wt q$ and $q$. Obviously the inclusion $\wt i$ is injective, and the left square commutes by comparing definitions of $\wt G',\wt K'$ and $G',K'$. Exactness of the bottom row at $K'$ implies that $\ker\,\wt q$ is at most a $\Z_2$ extensions of $\im\,\wt i$; but since $\wt q(u)=u$ we must have $\ker\,\wt q=\im\,\wt i$. Hence the middle row is exact, and in particular $\wt G'$ does not contain the fiber class $u$, whence $p:\wt G'\to G'$ is an isomorphism.

Part (b) follows by repeating the proof of Prop.~\ref{prop:glue}c. \; $\square$ \medskip

\section{Paths on the small boundary}
\label{sec:path}

In this section we introduce an algebra $\mb P$ of paths on the boundary of a framed, triangulated 3-manifold. It is inspired by a construction of Neumann in \cite{Neumann-combinatorics} (and earlier \cite{NZ}). This algebra is an extremely useful tool for describing cycles in twisted homology $\wt H_1^-(\Sigma)$, keeping track of these cycles under gluing operations, and ultimately relating the cycles to functions on the non-abelian moduli space $\CX[\CC]$.

\subsection{Paths and generators for $H_1^-(\Sigma)$}
\label{sec:path-gen}

Let $\Sigma\overset\pi\to \CC$ be the (oriented) canonical cover of the boundary of a framed 3-manifold $M$ (with or without defects), with a fixed 2d triangulation $\mb t_{\rm 2d}$ of $\CC_{\rm big}$. Let $\Sigma^+$ and $\Sigma^-$ denote the two sheets of $\Sigma$, constructed from a branching graph induced by $\mb t_{\rm 2d}$ as in Section \ref{sec:cc} (so $\Sigma^+\simeq \Sigma^- \simeq \CC\bs \Gamma_{\rm br}$).

Let $E^\circ$ denote the disjoint union of the interiors all small edges of hexagons in $\mb t_{\rm 2d}$. Note that all the small edges lie on the interface $\pd \CC_{\rm small}\cap \pd \CC_{\rm big}$ between big and small parts of $\CC$. We define $\mb P$ (depending on $\mb t_{\rm 2d}$) as the relative homology group
\be \mb P := H_1(\CC_{\rm small},E^\circ)\,. \label{defP} \ee
Thus, $\mb P$ is generated by oriented paths $\pp$ on $\CC_{\rm small}$\,: either closed paths, or paths whose endpoints lie in the interior of one of the small edges of $\mb t_{\rm 2d}$ (and therefore on the boundary of a distinguished hexagonal face of $\mb t_{\rm 2d}$).

\begin{figure}[htb]
\centering
\includegraphics[width=5.5in]{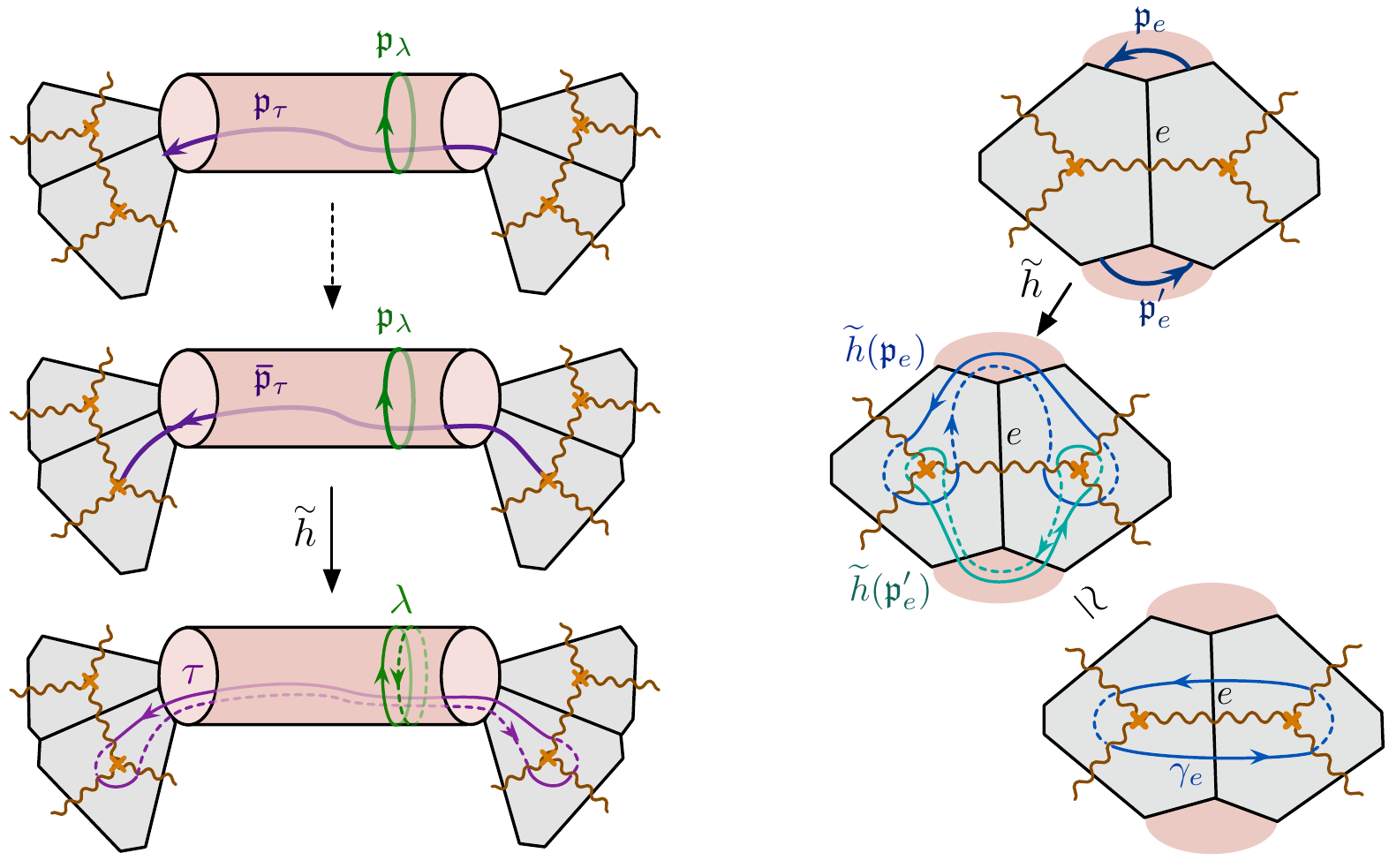}
\caption{Examples of the map $h:\mb P[\mb t_{\rm 2d}]\to H_1^-(\Sigma)$. Left: a closed path $\pp(\lambda)$ surrounding a small annulus and an open path $\pp(\tau)$ traversing it. Right: two paths mapping equivalently to the edge cycle $\gamma_e$. For convenience in visualizing and counting intersection numbers, we slightly deform the curves representing $h(\pp)$ for open $\pp$ away from branch points.}
\label{fig:defh}
\end{figure}

There is a homomorphism
\be \wt h\,:\; \mb P \to \wt H_1^-(\Sigma) \label{defh} \ee
defined as follows.
If $[\pp]\in \mb P$ can be represented by a smooth, closed curve $\pp$, we use the fact that there is no branching along the small boundary define lifts $\pp^\pm$ to $\Sigma^\pm\subset \Sigma$, and then (canonically) to $T_1\Sigma$. We set $\wt h([\pp]) := \ell^-(\pp) = [\pp^+]-[\pp^-]$ to be the class of the ``odd double lift'' of $\pp$. Note that this does not depend on which smooth representative of the curve $\pp$ is used, since changing the representative (as in Figure \ref{fig:u}) can only modify $\wt h([\pp])$ by an even multiple of $u$.

\begin{figure}[htb]
\centering
\includegraphics[width=6in]{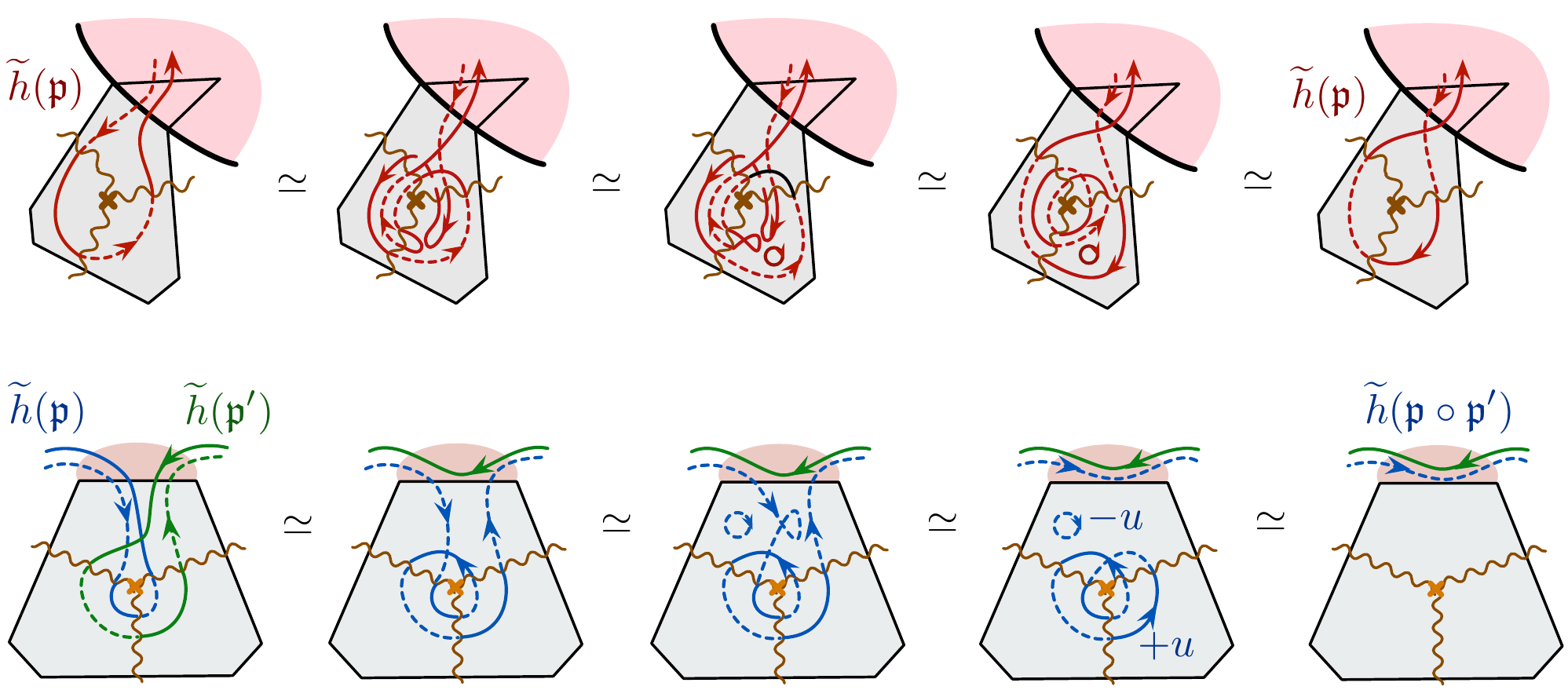}
\caption{Top: verification that it does not matter how $\wt h(\pp)$ is chosen to wind around branch points. Bottom: proof that $\wt h(\pp)+\wt h(\pp') = \wt h(\pp\circ\pp')=\wt h(\pp+\pp')$ in $\wt H_1^-(\Sigma)$, assuring that $\wt h$ is well defined. We represent elements of $\wt H_1(\Sigma)$ by smooth curves on $\Sigma$, keeping track of the fiber class $u$ by means of small loops. In each case, corrections involving the fiber class cancel out at the end.}
\label{fig:defhtw}
\end{figure}

If $[\pp]\in \mb P$ can be represented by a smooth open path $\pp$, then we can extend $\pp$ slightly to a path $\ol \pp$ that starts and ends at the branch points on the faces at the ends of $\pp$ (as on the left of Figure \ref{fig:defh}). We take lifts $\ol\pp^\pm$ to $\Sigma^\pm$, and connect them at the branch points to form a smooth closed curve $\pp^+\circ \pp^-$; we then define $\wt h([\pp]):= \ell^-(\pp)=[\pp^+\circ \pp^-]$. For convenience, we can deform the curve $\pp^+\circ \pp^-$ so that it passes around the branch points, as on the bottom-left of Figure \ref{fig:defh}. An easy exercise verifies that it does not matter whether we pass clockwise or counter-clockwise (top of Figure \ref{fig:defhtw}).

We extend $\wt h$ to all of $\mb P$ by linearity. For the map $\wt h$ to be a well-defined homomorphism, we must check that $\wt h(-\pp)=-\wt h(\pp)$, which follows easily from our definitions, and that
\be \wt h(\pp\circ \pp')=\wt h(\pp)+\wt h(\pp')\,, \label{hsum} \ee
where $\pp\circ\pp'$ is a smooth path homotopic to the concatenation of $\pp$ and $\pp'$ (possibly a cyclic concatenation, in which case $\pp\circ\pp'$ is closed). The proof of \eqref{hsum} is entirely local (depending on behavior at the endpoints of $\pp$ and $\pp'$), and illustrated on the bottom of Figure \ref{fig:defhtw}.

The map $\wt h$ has a nontrivial kernel. If $\pp_e$ and $\pp_e'$ are two paths circling counter-clockwise around the two ends of a big edge $e$ of $\mb t_{\rm 2d}$, as on right of Figure \ref{fig:defh}, then
\be \wt h(\pp_e) = \wt h(\pp_e')\,. \label{he} \ee
Let $\mb P_E=\langle \pp_e-\pp_e'\rangle_{\text{big edges $e$}}$ be the subgroup of $\mb P$ generated by the differences of paths associated to all edges of $\mb t_{2d}$. Clearly $\mb P_E \subset \ker\, \wt h$. In fact,
\begin{lemma} \label{lemma:gen}

If $M$ has no defects, $\ker\, \wt h = \mb P_E$, and $\im\, \wt h \simeq H_1^-(\Sigma)\subset \wt H_1^-(\Sigma)$. In other words,
\be 0\to \mb P_E \hookrightarrow \mb P \overset{\tilde h}\to H_1^-(\Sigma)\to 0 \ee
is exact. Thus $\wt h$ provides a splitting $\wt H_1^-(\Sigma)\simeq H_1^-(\Sigma)\oplus \Z_2$; and we may trivially extend $\wt h$ to a surjection $\wt h:(\mb P\oplus \Z_2)\to\hspace{-.3cm}\to \wt H_1^-(\Sigma)$, defining $\wt h(u)=u$ for the generator $u$ of $\Z_2$\,:
\be 0 \to \mb P_E \hookrightarrow \mb P\oplus \Z_2 \overset{\tilde h}\to \wt H_1^-(\Sigma)\to 0\,. \ee
\end{lemma}

\noindent {\it Proof.} Consider the map $h:= p\circ \wt h:\mb P\to H_1^-(\Sigma)$. Momentarily we will demonstrate that $h$ is surjective with $\ker\, h=\mb P_E$, by producing explicit presentations of $\mb P$ and $H_1^-(\Sigma)$. (These explicit presentations will be used throughout the rest of the paper.) Therefore $\ker\,\wt h \subset \mb P_E$, whence $\ker\,\wt h=\mb P_E$, and the rest of the Lemma follows. \; $\square$

\subsection{Generators of $H_1^-(\Sigma)$ and intersection form}
\label{sec:gen}

Suppose that the triangulated framed 3-manifold $M$ has no defects, and that $\pi_1(\CC_{\rm small}$ is abelian. In other words, the connected components of $\CC_{\rm small}$ are discs, annuli, tori, or spheres. We describe a set of generators for $H_1^-(\Sigma)$ labelled by paths $\pp\in\mb P$. Let $h:=p\circ\wt h:\mb P\to H_1^-(\Sigma)$.

Small spheres are connected components of $\CC$, over which both the path algebra and odd homology are trivial, so we may ignore them. Small tori are also connected components of $\CC$. For each small torus $t$, choose a basis of A and B cycles, let $\pp_\alpha^{(t)},\pp_\beta^{(t)}$ be paths representing these cycles --- generating the path algebra on $t$ --- and set 
\be
\alpha^{(t)} = h(\pp_\alpha^{(t)})\,,\quad  \beta^{(t)} = h(\pp_\beta^{(t)})\,.
\ee
Clearly $\alpha^{(t)},\beta^{(t)}$ generate the part of $H_1^-(\Sigma)$ associated to the small torus.

The remaining connected components of $\CC$ (which are all we consider now) contain big boundary $\CC_{\rm big}$ with holes connected to small discs or annuli. For each small annulus $a$ let $\pp_\lambda^{(a)}$ be a closed path generating $H_1(a)$, and choose an open path $\pp_\tau^{(a)}$ traversing the annulus from one end to the other (Figure \ref{fig:defh}, left). Set
\begin{subequations} \label{defHgens}
\be \lambda^{(a)} = h(\pp_\lambda^{(a)})\,,\quad \tau^{(a)} = h(\pp_\tau^{(a)})\,.\ee
For each edge $e$ of the big-boundary triangulation $\mb t_{\rm 2d}$, let $\pp_e$ and $\pp'_e$ be the open paths on $\CC_{\rm small}$ winding counter-clockwise around the two ends of $e$ (Figure \ref{fig:defh}, right), and set
\be  \gamma_e = h(\pp_e) = h(\pp'_e)\,.\ee
\end{subequations}

\begin{figure}[htb]
\centering
\includegraphics[width=5.8in]{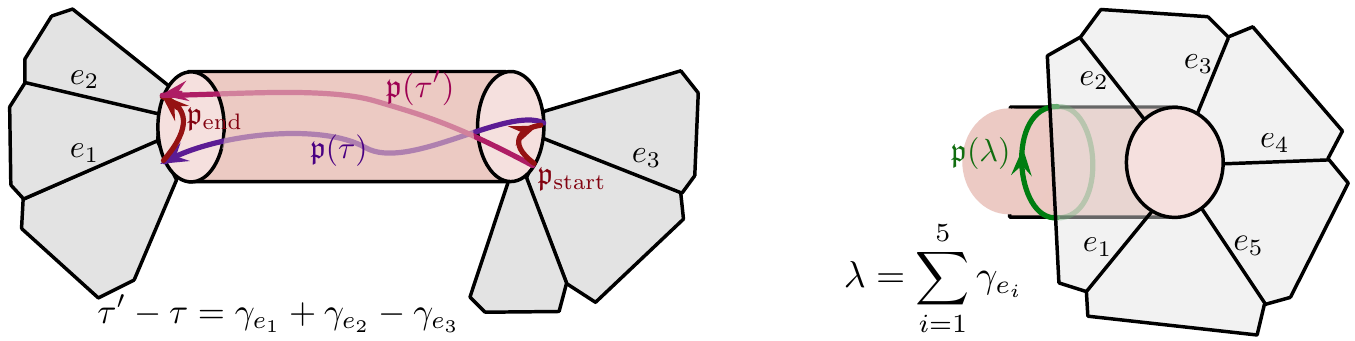}
\caption{Examples of relations among elements $h(\pp)$. On the left, we move the endpoints of $\pp_\tau$ by concatenating with (sums of) edge-paths, $\pp_\tau'=\pp_{\rm end}\circ \pp_\tau \circ \pp_{\rm start}$.}
\label{fig:TLrelsh}
\end{figure}

There are some relations among these elements. It is easy to see the path algebra $\mb P$ for small discs and annuli is generated by $\pp_e,\pp_e',\pp_\lambda^{(a)},\pp_\tau^{(a)}$ subject to the relations that the sum of paths around the boundary of each small disc $d$ vanishes, and the sums of paths around the two boundaries of each annulus $a$ are both equal to $\pp_\lambda^{(a)}$,
\be \label{prels}
\mb P = \Big\langle \{\pp_\lambda^{(a)},\pp_\tau^{(a)}\}_a,\, \{\pp_e,\pp_e'\}_e\,\Big| \sum_{\text{around $\pd\,d$}}(\text{$\pp_e$ or $\pp_e'$}) = 0, \sum_{\text{around $\pd_1a$}} (\text{$\pp_e$ or $\pp_e'$}) = -\hspace{-.2in}\sum_{\text{around $\pd_2a$}} (\text{$\pp_e$ or $\pp_e'$}) = \pp_{\lambda}^{(a)}\Big\rangle\,.
\ee
Also note that any two choices of traversing paths $\pp_\tau^{(a)}$, $\pp_\tau^{(a)}{}'$ are related by adding or subtracting appropriate edge paths $\pp_e,\pp_e'$ (Figure \ref{fig:TLrelsh}, left).

Similarly, since the branching graph $\Gamma_{\rm br}$ (dual to big edges of $\mb t_{\rm 2d}$) together with the extended paths $\ol\pp_\tau^{(a)}$ forms the 1-skeleton of a cell decomposition of $\CC$, lifting to the 1-skeleton of a cell decomposition of $\Sigma$, we find that (ignoring small spheres and tori)
\be \label{Hrels}
H_1^-(\Sigma) = \Big\langle \{\lambda^{(a)},\tau^{(a)}\}_a,\,\{\gamma_e\}_e\,\Big|\, \sum_{\text{$e$ on $\pd d$}}\gamma_e = 0,\, \sum_{\text{$e$ on $\pd_1a$}} \gamma_e = -\sum_{\text{$e$ on $\pd_2a$}} \gamma_e = \lambda^{(a)}\,\Big\rangle\,,
\ee
with each relation coming from a 2-cell. (Such cell decompositions with 0-cells at branch points are discussed in greater detail in Appendix \ref{app:cell}. See Lemma \ref{lemma:cell} there.) Applying $h$ to the presentation \eqref{prels} of $\mb P$ clearly produces the presentation \eqref{Hrels} of $H_1^-(\Sigma)$, with kernel precisely $\mb P_E=\langle \pp_e-\pp_e' \rangle$.

By direct inspection, we also see that

\begin{figure}[htb]
\centering
\includegraphics[width=5.8in]{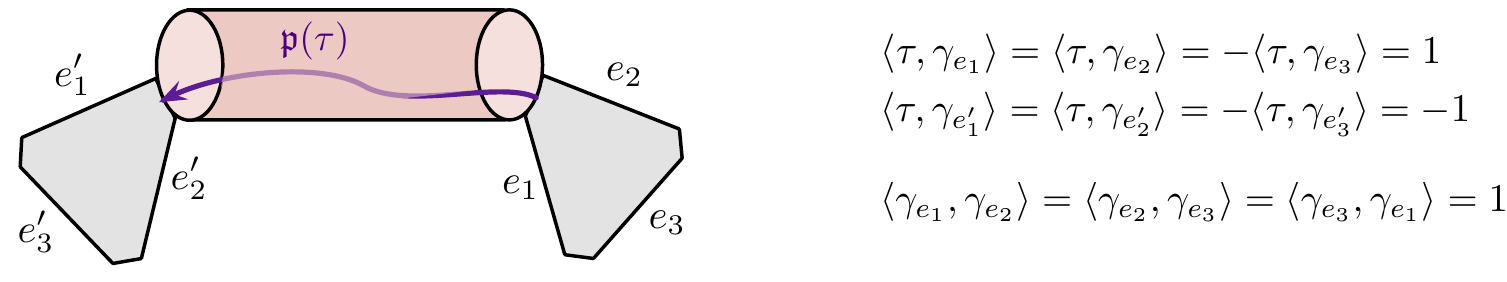}
\caption{Intersection product among $\tau$'s and $\gamma_e$'s.}
\label{fig:Peth}
\end{figure}

\begin{lemma} \label{lemma:intH}
The intersection form among generators of $H_1^-(\Sigma)$ (including small tori and spheres) is
\be \label{int-H}
\begin{array}{l}
 \langle \alpha^{(t)},\beta^{(t')} \rangle = 2 \langle \pp_\alpha^{(t)},\pp_\beta^{(t')}\rangle = 2 \,\delta_{tt'} \\
 \langle \tau^{(a)}, \lambda^{(a')}\rangle = 2 \langle \pp_\tau^{(a)},\pp_\lambda^{(a')} \rangle = \pm 2 \,\delta_{aa'} \\
 \langle \tau^{(a)},\gamma_e \rangle = \pm 1 \quad\text{for the six $e$ on faces adjacent to $\pd \pp_\tau^{(a)}$} \\
 \langle \gamma_e ,\gamma_{e'} \rangle = \text{\# faces shared by $e,e'$}\,,
\end{array}
\ee
where the signs in $\langle \tau,\lambda\rangle$ depends on orientation, the signs in $\langle \tau,\gamma_e\rangle$ are shown in Figure \ref{fig:Peth}, and the faces in $\langle \gamma_e,\gamma_e'\rangle$ are also counted with orientation, $+1$ $(-1)$ if $e'$ occurs counter-clockwise (clockwise) after $e$.
\end{lemma}

We have shown that $h:\mb P\to H_1^-(\Sigma)$ is surjective with kernel $\mb P_E$ when $\pi_1(\CC_{\rm small})$ is abelian. The same basic argument can easily be extended to reach the same conclusion for general $\CC_{\rm small}$ (as long as $M$ has no defects). For closed components of $\CC_{\rm small}$ of any genus $g$, we choose $2g$ A and B cycles to generate $\mb P$; then applying $h$ to these cycles (\ie\ taking odd double lifts) directly produces generators of $H_1^-(\Sigma)$. For components of $\CC$ containing both big and small parts, we choose multiple non-intersecting paths akin to the $\pp_\tau$'s, traversing $\CC_{\rm small}$ from boundary to boundary, so that the augmentation of $\Gamma_{\rm br}$ by these paths continues to produce a 1-skeleton of a cell decomposition of $\CC$. Then these generalized traversing paths together with $\pp_e,\pp_e'$ generate $\mb P$ (with the only relations being that the sum of $\pp_e,\pp_e'$ paths around the boundary of any component of $\CC_{\rm small}$ must vanish); while applying $h$ produces generators of $H_1^-(\Sigma)$ (with the only new relations of the form $h(\pp_e)=h(\pp_e')$). This completes the proof of Lemma \ref{lemma:gen} above.

\subsubsection{Example: the tetrahedron}
\label{sec:tetH}

\begin{wrapfigure}{r}{2in}
\centering
\includegraphics[width=1.5in]{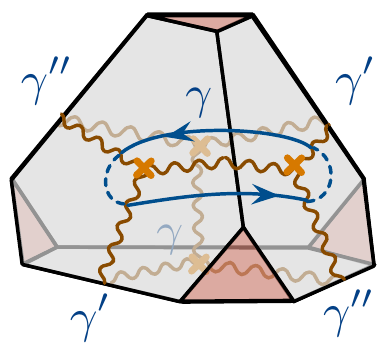}
\caption{Six edge cycles generating $H_1^-(\Sigma_\Delta)$.}
\label{fig:tetH}
\end{wrapfigure}

Let us illustrate Lemma \ref{lemma:gen} when $M=\Delta$ is a truncated tetrahedron, with its canonical triangulation $\mb t_{\rm 2d}$. The cover $\Sigma_\Delta\overset\pi\to \CC_\Delta=\pd M_\Delta$ is branched over four points (one in the center of each big hexagonal face of the tetrahedron), so
\be \CC_\Delta\simeq \pd \Delta \simeq S^2\,,\qquad \Sigma_\Delta \simeq T^2\,,\ee
Since $\pi_*:H_1(\Sigma_\Delta)\to H_1(\CC_\Delta)$ is trivial, all of $H_1^-(\Sigma_\Delta) = H_1(\Sigma_\Delta)=\Z^2$ is odd. Lemma \ref{lemma:gen} shows that $H_1^-(\Sigma_\Delta)$ is generated by six edge cycles $\gamma_e$ subject to the four relations that the sum around every vertex vanishes. This implies that cycles associated to opposite edges are equal. Denoting pairs of edge cycles $\gamma,\gamma',\gamma''$ as in Figure \ref{fig:tetH}, we see that
\be \label{tetH}
\begin{array}{c} H_1^-(\Sigma_\Delta) = \Z\langle \gamma,\gamma',\gamma''\,|\,\gamma+\gamma'+\gamma'' = 0\rangle\,, \\[.2cm]
 \langle \gamma,\gamma'\rangle = \langle \gamma',\gamma\rangle
 = \langle \gamma'',\gamma\rangle = 1\,,
\end{array}
\ee
with intersection product as described by Lemma \ref{lemma:intH}.

\subsection{$\wt G$, $\wt K$, and the gluing map}
\label{sec:KG}

Let's now revisit the setup of Section \ref{sec:glue-gen}, where we start with a framed, triangulated 3-manifold $M$, glue interiors of pairs of faces to form $M_0$ (with defects), and fill in the defects to form $M'$. Let us assume that neither $M$ nor $M'$ have defects. We can use paths to explicitly describe subgroups $\wt G\subset \wt K \subset \wt H_1^-(\Sigma_0)$ appearing in Propositions \ref{prop:glue}, \ref{prop:glue}', and their embeddings into $\wt H_1^-(\Sigma)$ via the injection $\wt g$. Let $\mb t_{\rm 2d}$ denote the induced triangulation of $\CC_{\rm big}$ and $\mb t_{\rm 2d}'$ denote the induced triangulation of $\CC_{\rm big}'=(\CC_0)_{\rm big}$.
 Let $\mb P,\,\mb P_0,\,\mb P'$ denote the algebras of paths on the respective small boundaries $\CC_{\rm small},(\CC_0)_{\rm small},\CC'_{\rm small}$.

\begin{wrapfigure}[19]{r}{1.7in}
\centering
\includegraphics[width=1.6in]{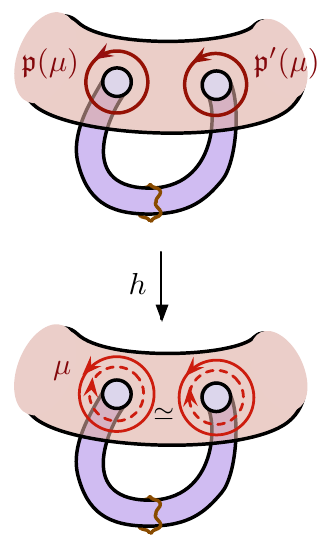}
\caption{Generators of $G$}
\label{fig:defmu}
\end{wrapfigure}

We begin with the subgroups of $\wt H_1^-(\Sigma_0)$. 
The small boundary $(\CC_0)_{\rm small}$ of $M_0$ has pairs of holes at the ends of each defect $I_j$, \ie\ at $(\CC_0)_{\rm small}\cap (\CC_0)_{\rm def}$. Let $\pp_{\mu_j}$ and $\pp'_{\mu_j}$ be paths on $(\CC_0)_{\rm small}$ circling counter-clockwise (from the perspective of the small boundary) around the holes at the ends of $I_j$, as in Figure \ref{fig:defmu}. Then the homology classes
\be \mu_j = \wt h(\pp_{\mu_j}) = \wt h(\pp'_{\mu_j}) \label{relmu} \ee
are precisely the cycles $\mu_j^+-\mu_j^-$ discussed in \eqref{defGtw}, and generate~$\wt G$. Letting $\mb P_G\subset \mb P_0$ denote the subgroup generated by $\pp_{\mu_j},\pp'_{mu_j}$, we have
\be \wt G = \im\,(\wt h: \mb P_G \to \wt H_1^-(\Sigma_0))\,.\ee
We then define a subgroup $\wt K \supset \wt G$ as
\be \wt K = \im\,(\wt h:(\mb P_0\oplus\Z_2) \to \wt H_1^-(\Sigma_0))\,, \ee
extending $\wt h$ as usual by $\Z_2=\Z\langle u\rangle/\Z\langle 2u\rangle$,  with $\wt h(u)=u$. Set $G = p(\wt G)=h(\mb P_G)$, $K=p(\wt K)=h(\mb P_0)$.

It is clear%
\footnote{To see $\tilde K\subset \tilde K'$ observe that $\tilde K \subset \ker\langle \tilde G,*\rangle|_{\tilde H_1^-(\Sigma_0)} = \ker\langle \tilde G',*\rangle|_{\tilde H_1^-(\Sigma_0)} = \tilde K'$.} %
that $\wt G\subset \wt K\subset \wt K'$, in the notation of Proposition \ref{prop:glue}'. To see that $\wt K/\wt G\simeq \wt H_1^-(\Sigma')$ so that $\wt K$ satisfies the hypotheses of part (b) of Proposition \ref{prop:glue}', we use:

\begin{lemma} \hspace{-.15cm}{\bf a} \label{lemma:SESP}
The maps $\wt h:(\mb P_0\oplus \Z_2)\to \wt H_1^-(\Sigma_0)$ and $\wt h:(\mb P'\oplus \Z_2)\to \wt H_1^-(\Sigma')$ fit into a commutative diagram
\be \label{PKG}
\begin{array}{ccccccccc} 0 & \to & \mb P_G & \overset i\hookrightarrow & \mb P_0\oplus \Z_2 & \overset {q_{\mb P}}\to & \mb P'\oplus \Z_2 &\to&  0 \\[.1cm]
 && \wt h\, \raisebox{.3cm}{\rotatebox{-90}{$\to\hspace{-.3cm}\to$}} 
 && \wt h\, \raisebox{.3cm}{\rotatebox{-90}{$\to\hspace{-.3cm}\to$}}
 && \wt h\, \raisebox{.3cm}{\rotatebox{-90}{$\to\hspace{-.3cm}\to$}} \\
 0 & \to & \wt G & \overset {\tilde i}\hookrightarrow & \wt K & \overset {\tilde q}\to & \wt H_1^-(\Sigma') &\to&  0\,.
\end{array}
\ee
with top and bottom rows exact. Thus $\wt K$ satisfies the hypotheses of Prop. \ref{prop:glue}'(b).
\end{lemma}

\begin{figure}[htb]
\centering
\includegraphics[width=4.8in]{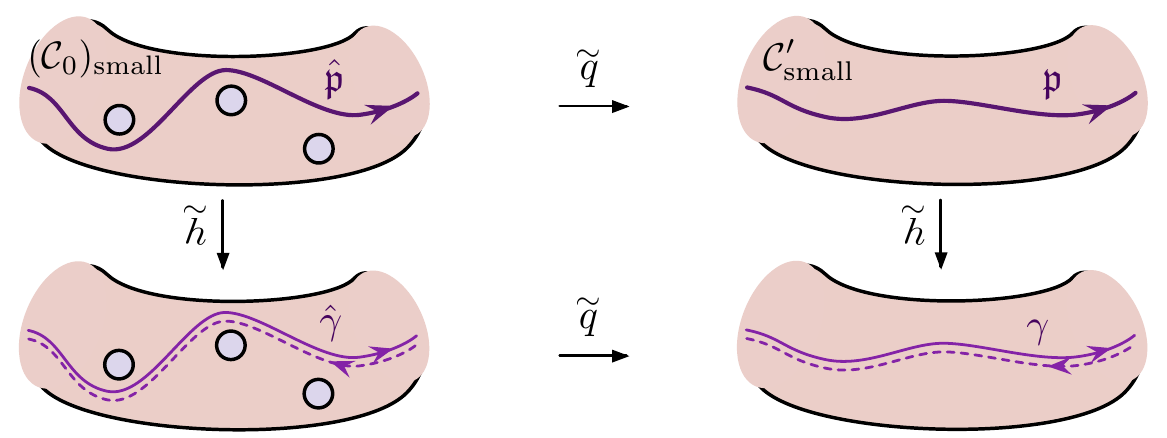}
\caption{Commutativity of the second square in \eqref{PKG}.}
\label{fig:PKG}
\end{figure}

\emph{Proof.}
The inclusions $\mb P_G \overset i\hookrightarrow \mb P_0$, $\wt G\overset {\tilde i}\hookrightarrow \wt K$ are obviously injective, and the first square  commutes due to the definitions of $\wt G,\wt K$.
On the top row, the map $q_{\mb P}$ includes a path $\pp$ on $(\CC_0)_{\rm small}$ as a path on $\CC_{\rm small}'$, where the holes from defects have been filled in; while $q_{\mb P}(u):=u$ for $u\in \Z_2$. We can invert $q_{\mb P}$ by lifting any path $\pp$ on $\CC_{\rm small}'$ to a path $\hat \pp$ on $(\CC_0)_{\rm small}$, choosing some way for $\hat\pp$ to pass around the defects; thus $\mb P' \simeq \mb P_0/\mb P_G$ and $\mb P' \oplus\Z_2\simeq (\mb P_0\oplus \Z_2)/\mb P_G$. 
(Alternatively, exactness of the top row follows from the long exact sequence for the pair $((\CC_0)_{\rm small}/E^\circ{}', \CC'_{\rm small}/E^\circ{}')$ where $E^\circ{}'$ is the union of interiors of small edges of $\mb t'_{2d}$.)
By definition of $\wt K$ and by Lemma \ref{lemma:gen} we know that all the downward maps are surjective. On the bottom row, we first observe that any cycle in $\wt K$ has zero intersection number with all defect cycles $\mu_j^\pm \in \wt G'$ (as in \eqref{defGtw}), so we may define $\wt q$ to act by including these cycles in $\wt H_1^-(\Sigma)$. In other words, $\wt q$ is the restriction of the map in \eqref{gluetw} to $\wt K\subset \wt K'$. With this definition, the second square also commutes (see Figure \ref{fig:PKG}), and surjectivity of $\wt q$ on the bottom row follows from surjectivity of $q_{\mb P}$ on the top. Thus $\wt K/\wt G \simeq \wt H_1^-(\Sigma')$.~\;$\square$ \medskip

\noindent{\bf Remark.} By applying $p$ to \eqref{PKG} and killing all the $\Z_2$'s and fiber classes, we see that $G=p(\wt G)=h(\mb P_G)$ and $K=p(\wt K)=h(\mb P_0)$ satisfy the hypotheses of Proposition \ref{prop:glue}(c).

Explicitly, the sequences \eqref{PKG} imply that to generate $K$ (or $\wt K=K\oplus \Z_2)$ we may start with generators $\gamma_e,\tau,...$ of $H_1^-(\Sigma')$ labelled by paths $\pp_e,\pp_\tau,...$, use $q_{\mb P}^{-1}$ to lift the paths to any $\hat \pp_e,\hat\pp_\tau,...\in \mb P_0$ (unique up to addition of $\pp_{\mu_j}$ and $\pp_{\mu_j}'$), adjoin the $\pp_{\mu_j}$'s, and then apply $h$ to obtain generators $\hat \gamma_e,\hat\tau,...$ and $\mu_j$ of $K$. The lifts $\hat \gamma_e,\hat\tau,...$ satisfy $q(\hat \gamma_e)=\gamma_e,\, q(\hat\tau)=\tau$, etc., and are unique up to addition of $\mu_j$'s.

\begin{figure}[htb]
\centering
\includegraphics[width=5.3in]{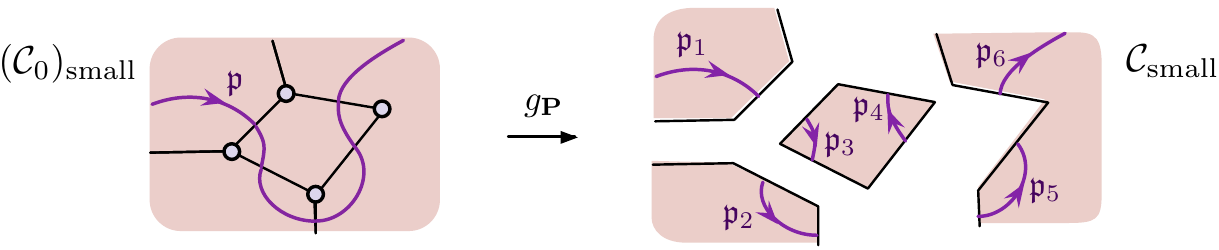}
\caption{The cutting map acting on paths.}
\label{fig:defgP}
\end{figure}

We may similarly use paths to find an algorithm for the action of the injection $\wt g:\wt H_1^-(\Sigma_0)\to \wt H_1^-(\Sigma)$ on $\wt K$.

Recall that $\mb P$ is a relative homology group $H_1(\CC_{\rm small},E^\circ)$ where $E^\circ$ denotes the disjoint union of interiors of small edges of $\mb t_{\rm 2d}$. Similarly, $\mb P_0 = H_1((\CC_0)_{\rm small},E^\circ{}')$, where $E^\circ{}'$ is the union of interiors of small edges of $\mb t_{\rm 2d}'$. The small edges $E^\circ$ have an image $\hat E^\circ\supset E^\circ{}'$ in $(\CC_0)_{\rm small}$ --- along which $(\CC_0)_{\rm small}$ is cut to obtain $\CC_{\rm small}$.
Since $E^\circ{}' \subset \hat E^\circ$, there exists a natural ``cutting map'' $g_{\mb P}:H_1((\CC_0)_{\rm small},E^\circ{}')\to H_1((\CC_0)_{\rm small},\hat E^\circ)\simeq H_1(\CC_{\rm small},E^\circ)$, or
\be \label{defgP} 
g_{\mb P} :\, \begin{array}{ccc} \mb P_0 &\to & \mb P \\
  \pp &\mapsto& \sum_i \pp_i\,. \end{array}
\ee
Explicitly, $g_{\mb P}$ takes a path on $(\CC_0)_{\rm small}$ to the sum of paths $\sqcup_i \pp_i = \pp\cap \CC_{\rm small}$ in its intersection with patches of $\CC_{\rm small}$, as in Figure \ref{fig:defgP}.
We also define the homomorphism
\be {\rm cut}: \begin{array}{ccc} \mb P_0 &\to& \Z_2 = \{0,1\} \\
 \pp &\mapsto & \#(\pp\cap (E^\circ\bs E^\circ{}'))
\end{array} \label{defcut}
\ee
to be the unoriented intersection number between a path (or sum of paths) $\pp$ and the small edges of $E^\circ{}'$, included as a subset of $(\CC_0)_{\rm small}$. In other words, ${\rm cut}(\pp)$ is the number of times a path $\pp$ is cut when the map $g_{\mb P}$ is applied. It is well defined modulo 2. (As an example, in Figure \ref{fig:defgP}, ${\rm cut}(\pp) = 5\equiv 1$.) Then we construct
\be \wt g_{\mb P}: \begin{array}{ccc} \mb P_0\oplus \Z_2 &\to & \mb P\oplus \Z_2 \\
\pp + n u &\mapsto & g_{\mb P}(\pp) + ({\rm cut}(\pp)+n)u\,.
\end{array}
\ee
Note that $\wt g_{\mb P}$ acts trivially on $u\in \Z_2$, but sends a path $\pp\in\mb P_0$ to $g_{\mb P}(\pp)+{\rm cut}(\pp)u\in \mb P\oplus \Z_2$. \medskip

\noindent \textbf{Lemma \ref{lemma:SESP}\hspace{.03cm}b} \; \textit{The maps $\wt h:\mb P_0\to H_1^-(\Sigma_0)$ and $\wt h:\mb P\to H_1^-(\Sigma)$ and the injection $\wt g$ fit into a commuting square
\be \label{gPh}
\begin{array}{ccc} \mb P_0\oplus \Z_2 & \overset{\tilde g_{\mb P}}\to & \mb P\oplus \Z_2 \\
   \wt h\downarrow && \wt h\downarrow \\
   \wt H_1^-(\Sigma_0) & \overset g \hookrightarrow & \wt H_1^-(\Sigma)\,,
\end{array}
\ee
so that $\wt g\wt h = \wt h \wt g_{\mb P}$.} \medskip

\begin{figure}[htb]
\centering
\includegraphics[width=5.7in]{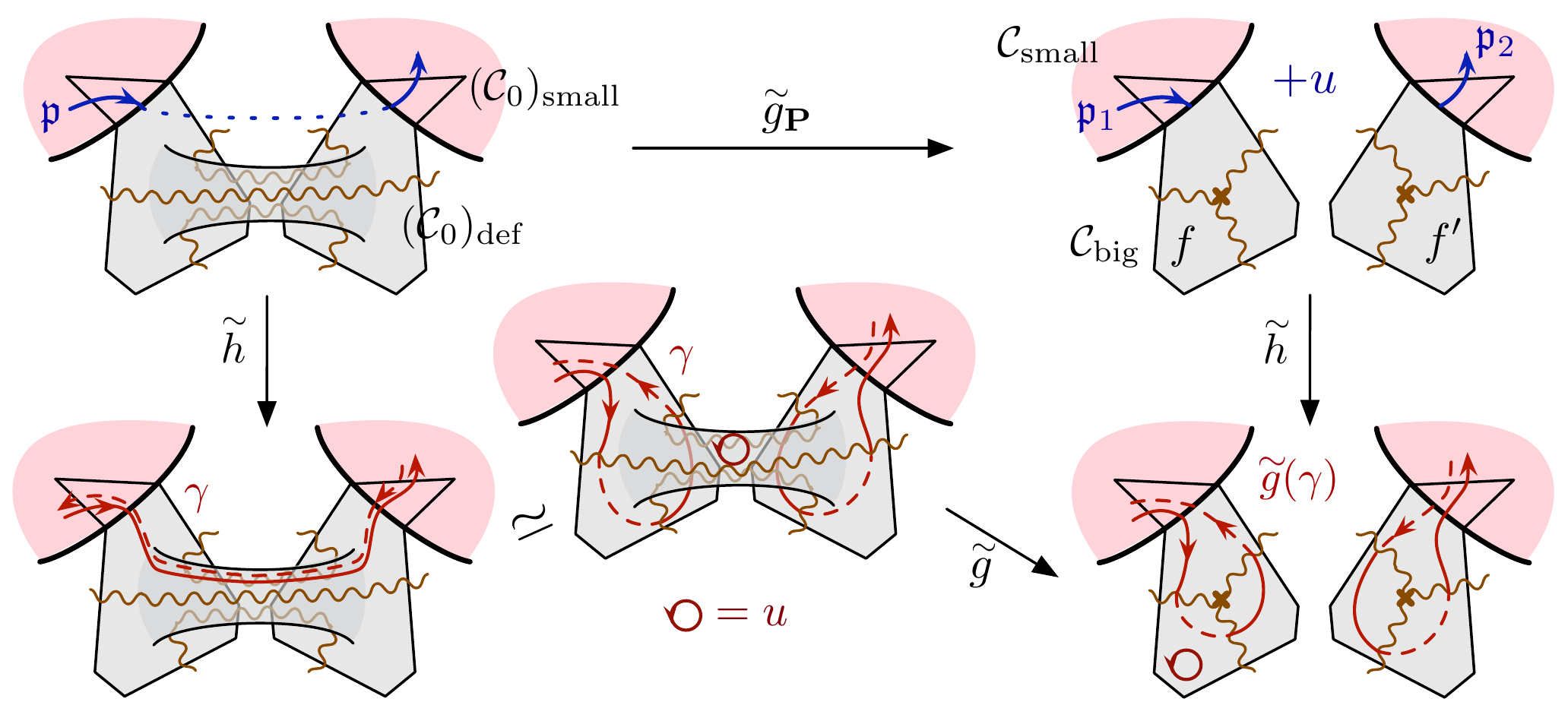}
\caption{Commutative diagram demonstrating $\wt g\circ \wt h(\pp) = \wt h\circ \wt g_{\mb P}(\pp)$. On the bottom row, the curve $\gamma=\wt h(\pp)$ in $T_1\Sigma_0$, which crosses the gluing region, is homologous to a sum of two smooth curves on either side of the gluing region along with a curve wrapping the fiber $u$, represented as a small loop. This extra ``twisted'' correction to $\wt g(\gamma)$ is accounted for by the contribution $u$ to $\wt g_{\mb P}(\pp) =g_{\mb P}(\pp) + {\rm cut}(\pp)u$ in the top row (here ${\rm cut}(\pp)=1$).}
\label{fig:pathhom-tw}
\end{figure}

\noindent \emph{Proof.} The twisted map $\wt g_{\mb P}$ was constructed precisely to ensure commutativity. The proof follows from a local argument, illustrated in  Figure \ref{fig:pathhom-tw}. \; $\square$ \medskip

\noindent {\bf Remark.} Lemma \ref{lemma:SESP}b lets us algorithmically calculate $\wt g(\gamma)\in \wt H_1^-(\Sigma)$ for any $\gamma\in \wt K=\im\,\big[\wt h:(\mb P_0\oplus \Z_2)\to \wt H_1^-(\Sigma_0)\big]$. By applying $p$, killing all $\Z_2$'s and fiber classes, we may also use the bare version of the cutting map $g_{\mb P}$ to calculate $g(\gamma)\in H_1^-(\Sigma)$ for any $\gamma\in K=\im\,\big[h:\mb P_0\to H_1^-(\Sigma_0)\big]$. In un-twisted homology, $gh=hg_{\mb P}$.

\subsection{Detour: Neumann's chain complex}
\label{sec:N}

We may use the formalism developed so far to provide a new topological interpretation for the chain complex that Neumann introduced in \cite{Neumann-combinatorics} to understand the symplectic properties of the hyperbolic gluing equations. (The actual connection with gluing equations will come in the next section.)

Suppose that $M'$ is a framed 3-manifold with no defects and no big boundary (the topology of $\CC'=\CC_{\rm small}'$ is unconstrained), which is glued from a collection $M=\sqcup_{i=1}^N \Delta_i$ of truncated tetrahedra. As usual we separate the gluing into two steps $M\leadsto M_0\leadsto M'$ where $M_0$ has defects $I_j$ along the edges of $M'$. In this section we work with untwisted homology, and untwisted maps, always projecting out the fiber class $u$.
We build a chain complex in several steps.

First, let the $\mu_j^\pm\subset \Sigma_0$ be the actual curves (as in Figure \ref{fig:SxSg} of Section \ref{sec:glue-gen}) used to represent the generators $\mu_j=\mu_j^+-\mu_j^-$ of $G\subset H_1^-(\Sigma_0)$, and let $N_\mu\subset \Sigma_0$ be a disjoint union of neighborhoods of these curves. (Alternatively $N_\mu$ is the lift of the defect boundary $(\CC_0)_{\rm def}$ to $\Sigma$.) Thus $N_\mu$ is a disjoint union of $2N$ annuli; and $H_1^-(N_\mu) \simeq \Z^{N}$ is generated by the $\mu_j$ as curves in $N_\mu$. The image of the inclusion
\be  H_1^-(N_\mu) \overset i\to H_1^-(\Sigma_0) \ee
is precisely $G$. We may compose this with the injection $g:H_1^-(\Sigma_0)\hookrightarrow H_1^-(\Sigma) = \oplus_{i=1}^N H_1^-(\Sigma_{\Delta i})$ (of finite cokernel) to get
\be H_1^-(N_\mu) \overset{g\,\circ\, i}\longrightarrow H_1^-(\Sigma)\,,\qquad \im\,(g\circ i) = g(G)\,.  \ee

Next, we dualize the map $g\circ i$ with respect to the canonical \emph{symmetric} pairing $(\mu_i,\mu_j)=\delta_{ij}$ on $H_1^-(N_\mu)$ and the intersection product $\langle\;,\;\rangle$ on $H_1^-(\Sigma)$ to get a complex
\be H_1^-(N_\mu) \overset{g\,\circ\, i}\longrightarrow H_1^-(\Sigma) \overset{(g\,\circ\, i)^*}\longrightarrow H_1^-(N_\mu)\,. \label{complex2} \ee
By definition, the kernel of $(g\,\circ\, i)^*$ is $K'':=\ker\, \langle g(G),*\rangle|_{H_1^-(\Sigma)}$. This includes $g(K')$, where $K'=\ker\,\langle G,*\rangle|_{H_1^-(\Sigma_0)}$; indeed, by Proposition \ref{prop:glue}a $K''/g(K')$ is torsion. Moreover, setting $K=\im\,[h:\mb P_0\to H_1^-(\Sigma_0)]$ as in Section \ref{sec:KG}, we obtain from Lemma \ref{lemma:gen} that $g(K)\subset g(K')\subset K''$ and $K''/g(K)$ is torsion. Thus the homology of \eqref{complex2} equals $g(K)/g(G)\simeq H_1^-(\Sigma')$ modulo torsion.

Note that $H_1^-(\Sigma')$ is isomorphic to the standard homology $H_1(\CC')=H_1(\CC_{\rm small}')$ since $\Sigma'\to \CC'$ is a trivial double cover (\cf\ Lemma \ref{lemma:dis} of Appendix \ref{app:cell}). The isomorphism $\ell^-:H_1(\CC')\overset\sim\to H_1^-(\Sigma')$ acts as multiplication by two on the intersection form.

We may further resolve the complex on the left using an odd version of the long exact sequence in relative homology for the pair $(\Sigma_0,N_\mu)$, 
\be \label{Neumann-cx}
\begin{array}{ccccc}
   H_2^-(\Sigma_0) &\overset{j}\to& H_2^-(\Sigma_0,N_\mu) &\overset{\delta}\to& H_1^-(N_\mu) \overset{g\,\circ\, i}\longrightarrow H_1^-(\Sigma) \overset{(g\,\circ\, i)^*}\longrightarrow H_1^-(N_\mu) \\
   \rotatebox{90}{$=$} && \rotatebox{90}{$=$} \\
   0 && H_2^-(\Sigma')\,.
\end{array}
\ee
Note that $H_2^-(\Sigma_0)=0$ because $\Sigma_0$ is connected; and $H_2^-(\Sigma_0,N_\mu) \simeq H_2^-(\Sigma') \simeq \Z^{\# \CC'}$ has one generator for every component of $\CC'=\pd M'$ (the odd-double-lift of the fundamental class of that component to $\Sigma'$). By Lemma \ref{lemma:cx} (Appendix \ref{app:cx}), the odd version of the long exact is exact modulo 2-torsion. Thus the homology of \eqref{Neumann-cx} is zero, except at $H_1^-(\Sigma)$, where it equals $H_1^-(\Sigma')$ (modulo torsion).

The complex \eqref{Neumann-cx} is identical to the one presented in Theorem 4.1 of Neumann \cite{Neumann-combinatorics}, but now the various groups involved have topological meaning. The maps $\delta$ and $g\circ i$ are what Neumann calls $\alpha$ and $\beta$. The isomorphism $\ell^-:H_1(\CC')\overset\sim\to H_1^-(\Sigma')$ is called $\delta$.
We could also extend the complex further to the right, as in \cite{Neumann-combinatorics}, by dualizing the long exact sequence. By computing the homology of \eqref{Neumann-cx} modulo torsion, we have re-proven \cite[Thm. 4.1]{Neumann-combinatorics}. Neumann goes further and carefully identifies the torsion groups as well. It would be interesting to investigate the torsion further from the perspective of odd homology.

\section{Framed flat connections on boundaries}
\label{sec:PGL}

Having described gluing in terms of odd homology --- both abstractly as in Proposition \ref{prop:glue} and concretely in terms of path algebras $\mb P$ --- we now proceed to describe how gluing acts on spaces of framed flat $PGL(2):=PGL(2,\C)$ connections.

We begin by reviewing the construction of cluster-like $\C^*$ coordinates for moduli spaces of framed flat $PLG(2)$ connections on boundaries of framed, triangulated 3-manifolds, following \cite{FG-Teich, DGG-Kdec} and \cite[Appendix A]{DGV-hybrid}.%
\footnote{Some additional subtleties and restrictions arise, especially as compared to \cite{DGG-Kdec} and Appendix A of \cite{DGV-hybrid}, because we focus only on boundaries of framed 3-manifolds rather than their interiors.} %
We show directly that the coordinates are labelled by cycles in odd homology, with Poisson, symplectic, and $K_2$ structures induced by the intersection form in odd homology (Propositions \ref{prop:coords}--\ref{prop:PB}).

We then explain gluing of framed $PGL(2)$ connections works, and derive the standard combinatorial gluing equations that appeared in \cite{DGG-Kdec, DGV-hybrid}, a generalization of Thurston's gluing equations for hyperbolic 3-manifolds. Observing that these equations take the promised form \eqref{eq-intro}, we deduce that gluing is $K_2$ symplectic reduction (Theorem \ref{thm:H}).

%We assume throughout that the 3-manifolds are defect-free.%
%
%\footnote{It should be clear from Section \ref{sec:gluePGL} that appropriate moduli spaces in the presence of defects can be defined as well. They have appeared, indirectly, in work such as \cite{AK-new}. To keep things (relatively) simple, we will not consider defects here.}
%

\subsection{Moduli spaces}
\label{sec:PGLdef}

Let $\CC=\pd M$ be the boundary of a framed 3-manifold, as defined in Section~\ref{sec:framed}. In particular, $\CC = \CC_{\rm big}\cup_{\sqcup_i S^1_i} \CC_{\rm small} \cup_{\sqcup_j S^1_j} \CC_{\rm def}$ has a splitting into big, small, and defect parts. Let $\mb t_{\rm 2d}$ be a triangulation of the big boundary --- \ie\ a tiling of $\CC_{\rm big}$ by hexagons. Let $E\to \CC$ be a trivial $\C^2$ bundle and $\PP E$ its projectivization (a $\cp^1$ bundle). Let $\CC^*$ be $\CC$ with a point removed in the center of every small disc of $\CC_{\rm small}$ --- thus $\CC^*$ is a ``punctured'' boundary. We define%
\be \begin{array}{rl} \CX[\CC] &:=\, \{\text{framed flat $PGL(2)$ connections on $E\to \CC^*$} \\ & \hspace{.5in} \text{with unipotent monodromy at punctures}\}  \\[.2cm]
 &:=\, \{\text{flat $PGL(2)$ connections on $E\to \CC^*$} \\
  & \hspace{.5in} \text{with unipotent monodromy at punctures} \\
  & \hspace{.5in} \text{and a choice of invariant flat line on each component of $\CC_{\rm small}$} \} \\
 &=\, \{\text{flat connections on $\PP E\to \CC^*$} \\
  & \hspace{.5in} \text{with unipotent monodromy at punctures} \\
  & \hspace{.5in} \text{and a choice of global flat section on each component of $\CC_{\rm small}$} \}\,.
\end{array}
\ee
We emphasize that these moduli spaces are defined modulo $PGL(2)$ gauge isomorphism.
As stated, the framing data is a choice of invariant flat line in the $\C^2$ bundle along the small boundary, or a flat section of the projectivized bundle.

If $M$ has no defects, and $\pi_1(\CC_{\rm small})$ is abelian (\ie\ $\CC_{\rm small}$ contains only small discs, annuli, tori, and spheres), then $\CX[\CC]$ has some additional properties.
In this case a choice of framing always exists and is usually unique.%
\footnote{Conversely, if $\pi_1(\CC_{\rm small})$ is nonabelian, the existence of framing is highly restrictive!} %
Forgetting the framing, one gets a map to the standard character variety
\be \label{charvar}
\CX[\CC] \overset{\text{forget}}\to {\rm Hom}^{\rm un}\big(\pi_1(\CC^*),\, PGL(2)\big)/PGL(2)\,
\ee
(where ``un'' means we restrict puncture holonomies to be unipotent).
For particular flat connections, there may be an enhanced choice of framing, which resolves some singularities in the character variety.
If $\CC$ is connected, the complex dimension at smooth points is 
\be \label{dimX}
d(\CX[\CC]) = \left\{\begin{array}{l@{\quad}l}
 \;\;0 & \text{$\CC$ is a small sphere} \\%[.2cm]
 \;\;2 & \text{$\CC$ is a small torus} \\%[.2cm]
 \begin{array}{l}(2g-2){\rm dim}(G) + n({\rm dim}(G)-{\rm rank}(G)) \\
 \hspace{.5in}= 6g-6+2n \end{array} & \chi(\CC^*)<0\,,
\end{array}\right.
\ee
where $n$ denotes the number of small discs in $\CC_{\rm small}$, or the number of punctures of $\CC^*$.
(If $\CC = \sqcup_i\CC_i$ contains multiple connected components, the space $\CX[\CC] = \prod_i\CX[\CC_i]$ factorizes.)

At its smooth points, the standard character variety has a non-degenerate complex symplectic form given by the Atiyah-Bott formula%
\footnote{Technically, \eqref{AB} is a formula on the space of all connections; the symplectic form on the character variety obtained by a symplectic reduction of \eqref{AB} as in \cite{AtiyahBott-YM}.}
\be \omega  =  \int_{\C^*} {\rm Tr}(\delta \CA\wedge\delta \CA)\,, \label{AB} \ee
which pulls back to a symplectic form on an open subset of $\CX[\CC]$ (which we also call $\omega$). In some cases it is known how to extend $\omega$ to a symplectic form on all of $\CX[\CC]$. For example, if $\CC$ consists only of big boundary and small discs (no annuli or tori), $\CX[\CC]$ coincides with the space studied by Fock and Goncharov \cite{FG-Teich} using cluster coordinates, and the appropriate extension of $\omega$ was defined therein.

We will not try to extend $\omega$ to a symplectic form on all of $\CX[\CC]$ here, though it does seem possible to do so. Instead we will work on open patches $\CP[\CC;\mb t_{\rm 2d}] \subset \CX[\CC]$ of the form $(\C^*)^{d(\CX[\CC])}$, labelled by triangulations $\mb t_{\rm 2d}$ (as in \cite{FG-Teich}), on which it is fairly straightforward to calculate and extend $\omega$.

\subsection{Coordinate functions}
\label{sec:PGLcoords}

Fix a 2d triangulation $\mb t_{\rm 2d}$ of $\CC_{\rm big}$ as above.
Let $\mb P = \mb P[\mb t_{\rm 2d}]$ be the abelian group of paths on $\CC_{\rm small}$ defined in Section \ref{sec:path}.
For each path $\pp\in \mb P$ we can define a $\C^*$-valued function $x_\pp$ on $\CX[\CC]$, which turns out to depend only on the homology class $\gamma=\wt h(\pp)\in H_1^-(\Sigma) \subset \wt H_1^-(\Sigma)$. Throughout this section, we use the splitting induced by $\wt h$ (Lemma \ref{lemma:gen}) to view $H_1^-(\Sigma)$ as a subgroup of $\wt H_1^-(\Sigma)$. We proceed as follows.

\begin{figure}[htb]
\centering
\includegraphics[width=6in]{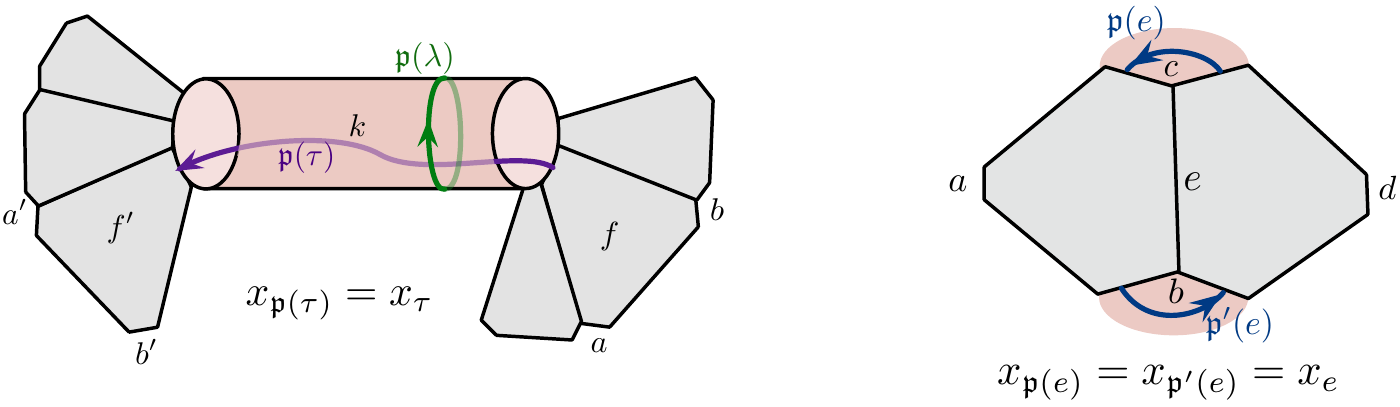}
\caption{Left: configuration of lines that defines the function $x_\pp$ for an open path (here for the ``twist'' of an annulus). Right: specialization of this configuration for paths that surround a big-boundary edge.}
\label{fig:twist-e-setup}
\end{figure}

For every $\pp\in \mb P$ represented by a closed path and $\CA \in \CX[\CC]$, set
\be x_\pp := \text{squared-eigenvalue of ${\rm Hol}_\pp(\CA)$ corresponding to the framing eigenline}\,. \ee
This makes sense since $\pp$ lies entirely on a single small boundary component, and the framing of $\CA$ specifies a distinguished eigenline on that component. For example, on a small torus $t$ the A- and B-cycle paths $\pp_\alpha^{(t)},\pp_\beta^{(t)}$ give rise to
\be x_\alpha^{(t)} :=  x_{\pp_\alpha^{(t)}} \,, \qquad  x_\beta^{(t)} := x_{\pp_\beta^{(t)}}\,, \label{xab} \ee
which are the squares of holonomy eigenvalues along respective cycles. On every small annulus $a$ there is a ``length'' function
\be  x_\lambda^{(a)} := x_{\pp_\lambda^{(a)}} \label{xlambda} \ee
measuring the holonomy eigenvalue around the girth of the annulus (Figure \ref{fig:twist-e-setup}, left). For a small disc, the path $\pp$ around the circumference of the disc gives $x_\pp = 1$ since the the holonomy must be unipotent there; this is consistent with the fact that $\pp$ is contractible.

If $\pp$ is represented by an open path, it starts and ends on the boundaries of some big hexagons $f,f'$, which are faces of $\mb t_{\rm 2d}$. For a given framed flat connection $\CA\in \CX[\CC]$, let $k$ denote the framing line of $\CA$ along $\pp$. Let $a,b$ (respectively, $a',b'$) denote the framing lines at the two small edges of $f$ (respectively, $f'$) disjoint from $\pp$, oriented as at the ends of $\pp(\tau)$ in Figure \ref{fig:twist-e-setup}. We can unambiguously parallel-transport $a,b,a',b'$ and $k$ to a single fiber of $E\to \CC$ above (say) the midpoint of $\pp$. For example, we first transport $a,b$ into the hexagon $f$, then transport them to the start of $\pp$, then along $\pp$; and similarly for $a',b'$. Working in a single fiber, we define the function $x_\pp$ to be the generalized cross-ratio
\be x_\pp = \frac{\langle a\wedge k\rangle \langle b\wedge k\rangle}{\langle a\wedge b\rangle}\cdot\frac{\langle a'\wedge b'\rangle}{\langle a'\wedge k\rangle \langle b'\wedge k\rangle}\, \in \C^*\,. \label{xt} \ee
Here $\langle *\wedge*\rangle$ is a skew-symmetric volume form in the fiber of $E$. To calculate \eqref{xt} one chooses any vectors in the lines $a,b,a',b',k$ and applies the volume form. The result is independent of the normalization of these vectors, independent of the volume form, and invariant under a $PGL(2)$ action on the fiber --- thus independent of the precise point along $\pp$ at which we compare lines.
(Equivalently, working in a projectivized $\cp^1$ bundle, one can think of \eqref{xt} as a generalized cross-ratio of five points in $\cp^1$: $(a-b)(a'-k)(b'-k)/[(a-k)(b-k)(a'-b')]$.)
Note that the cross-ratio \eqref{xt} only makes sense at points of $\CX[\CC]$ where all pairs of lines $a,b,a',b',k$ being compared are independent --- \ie\ the configuration of lines is in ``general position.''

\begin{defn} \label{def:CP} Let $\CP[\CC;\mb t_{\rm 2d}]\subset \CX[\CC]$ be the algebraically open subset on which functions $x_\pp\in \C^*$ are well defined for all $\pp\in \mb P[\mb t_{\rm 2d}]$ represented by open paths (meaning all pairs of lines being compared to form $x_\pp$ are in general position), with the additional restriction%
\footnote{This is largely a technical restriction, necessary for non-degeneracy to hold in Prop. \ref{prop:coords}d below. The issue is that unipotency does not uniquely specify the conjugacy class of a holonomy matrix --- it might be either trivial or parabolic. From a 3d perspective, this conjugacy class is fixed (and generically parabolic), so restricting to non-unipotent holonomy was not necessary in (\eg) Appendix A.3 of \cite{DGV-hybrid}, analogous constructions in \cite{DGG-Kdec}, or many years of hyperbolic geometry.} %
that for each connected component of $\CC_{\rm small}$ the holonomies around non-contractible cycles of that component are not simultaneously unipotent. (For example, $x_\lambda^{(a)}\neq 1$ for annuli and $(x_{\alpha}^{(t)},x_{\beta}^{(t)})\neq (1,1)$ for tori.) Often we simply write $\CP[\CC]$.
\end{defn}

\noindent {\bf Remark.} For components of $\CC$ consisting of big boundary with holes filled by discs, $\CP[\CC;\mb t_{\rm 2d}]$ coincides with a cluster-coordinate chart of \cite{FG-Teich}. \medskip

\begin{prop} \label{prop:coords} \hspace{.1cm}\\
\indent a) For both open and closed paths $x_{-\pp} = x_{\pp}^{-1}$, and under concatenation $x_{\pp\,\circ\, \pp'} = x_{\pp}x_{\pp'}$. Therefore, we may extend path functions to all $\pp\in \mb P$ by linearity, defining $ x_{\pp+\pp'} := x_\pp x_{\pp'}$, to obtain a map
\be x: \begin{array}{rccl} \CP[\CC]\,\times &\mb P &\to& \C^* \\
    & \pp &\mapsto & x_\pp\,,
    \end{array} \label{homP}
\ee
that's a homomorphism on the second factor.

b) The function $x_\pp$ depends only on the homology class $\gamma=\wt h(\pp)\in H_1^-(\Sigma)\subset \wt H_1^-(\Sigma)$.

c) If $M$ has no defects, the surjection $\wt h:\mb P\to\hspace{-.3cm}\to H_1^-(\Sigma)$ from \eqref{defh} provides a map
\be x: \begin{array}{rccl} \CP[\CC]\,\times &H_1^-(\Sigma) &\to& \C^* \\
    & \gamma &\mapsto & x_\gamma\,,
    \end{array}  \label{xXH}
\ee
that's a homomorphism on the second factor. By defining $x_u:=-1$, we may extend \eqref{xXH} to a map $x: \CP[\CC] \times \wt H_1^-(\Sigma) \to \C^*$.

d) If in addition $\pi_1(\CC_{\rm small})$ is abelian%
\footnote{This assumption may be lifted, by using a more general version of the reconstruction procedure (with unipotent modifications to remove extraneous punctures) from Appendix \ref{app:traffic}.}, %
the map \eqref{xXH} is non-degenerate in the sense that any basis $\{\gamma_i\}_{i=1}^d$ for $H_1^-(\Sigma)\subset \wt H_1^-(\Sigma)$ provides an injection
\be  \label{coord-inj}
(x_{\gamma_i},...,x_{\gamma_d}) : \CP[\CC] \hookrightarrow (\C^*)^{{\rm rank}\,H_1^-(\Sigma)}\,, \ee
with image the subset where $x_\lambda^{(a)} \neq 1,\, (x_{\alpha}^{(t)},x_{\beta}^{(t)})\neq (1,1)$.

e) With the assumptions in (d), ${\rm rank}\,H_1^-(\Sigma)=d(\CX[\CC])$, so $\CP[\CC]\subset \CX[\CC]$ is a subset of maximal dimension, with coordinates $\{x_\gamma\}_{\gamma\in H_1^-(\Sigma)}$, isomorphic to $(\C^*)^{d(\CX[\CC])}\big|_{x_\lambda^{(a)} \neq 1,\, (x_{\alpha}^{(t)},x_{\beta}^{(t)})\neq (1,1)}$.

\end{prop}

\noindent \emph{Proof.} For part (a), $x_{-\pp}=x_\pp^{-1}$ is clear from the definitions and $x_{\pp\circ\pp'}=x_\pp x_{\pp'}$ follows from from a straightforward local calculation. When concatenating $\pp\circ\pp'$ to produce another open path, numerators and denominators of \eqref{xt} cancel out to ensure $x_{\pp\circ\pp'}=x_{\pp}x_{\pp'}$. When producing a closed path, all numerators and denominators cancel up to an overall factor, which is precisely the square of the holonomy eigenvalue along the closed path.

For part (b), observe that of $\pp_e,\pp_e'$ are paths around the endpoint of a big edge $e$ of $\mb t_{\rm 2d}$ then $x_{\pp_e}=x_{\pp_e'}$ (see \eqref{xe} below). If $M$ has no defects, the claim follows from Lemma \ref{lemma:gen}. If $M$ has defects, we must do a bit more work. Let $M'$ be the framed 3-manifold obtained by filling in the defects of $M$, and consider the composition $q\circ  h:\mb P\to H_1^-(\Sigma')$ as in \eqref{PKG} (we have trivially factored out a $\Z_2$, using untwisted versions of these maps). The kernel of this map is $\mb P_G\oplus \mb P_E$, where $\mb P_E=\langle \pp_e-\pp_e'\rangle$ are the relations of Lemma \ref{lemma:gen}. Thus the kernel $\mb P_{\rm def}$ of $h$ is a subgroup $\mb P_E \subset \mb P_{\rm def}\subset \mb P_G\oplus \mb P_E$. Indeed, $\mb P_{\rm def}$ is generated by elements $\langle \pp_{\mu_j}-\pp_{\mu_j}'\rangle$ for all defects $I_j$, corresponding to $h(\pp_{\mu_j})-h(\pp_{\mu_j}')=0$ in \eqref{relmu}, in addition to the edge relations in $\mb P_E$. Thus, for the claim of part (b), we must show that $x_{\pp_{\mu_j}}=x_{\pp_{\mu_j}'}$ for all defects $I_j$. This holds because the closed paths $\pp_{\mu_j}$ and $\pp_{\mu_j}'$ are homotopic on $\CC$, so the eigenvalues of holonomies around these paths must be equal up to inversion; since the relative orientations of $\pp_{\mu_j}$ and $\pp_{\mu_j}'$ with respect to $\CC_{\rm small}$ are reversed and the sheets of the cover $\Sigma$ are \emph{also} reversed at opposite ends of a defect, the eigenvalues $x_{\pp_{\mu_j}}$ and $x_{\pp_{\mu_j}'}$ are exactly equal.

Part (c) follows immediately from part (b) and Lemma \ref{lemma:gen}.

Injectivity of \eqref{coord-inj} (and the identification of the image) in (d) follows from Appendix \ref{app:traffic}, where we recreate the unique framed flat connection $\CA\in \CP[\CC]$ associated to any point $x_\gamma\in (\C^*)^{d(\CX[\CC])}\big|_{x_\lambda^{(a)} \neq 1,\, (x_{\alpha}^{(t)},x_{\beta}^{(t)})\neq (1,1)}$.

The equality ${\rm rank}\,H_1^-(\Sigma)=d(\CX[\CC])$ in part (e) follows by comparing \eqref{rankH} and \eqref{dimX}.  \; $\square$ \medskip

When $M$ has no defects and $\pi_1(\CC_{\rm small})$ is abelian, the functions $x_\gamma$ are actually quite familiar. 
Above, we already produced functions $x_\alpha^{(t)},x_\beta^{(t)}$ for small tori and ``lengths'' $x_\lambda^{(a)}$ for small annuli.
Going further, we may choose a path $\pp_\tau^{(a)}$ traversing each annulus. Then
\be x_\tau^{(a)} = x_{\pp_\tau^{(a)}} \ee
is a complex generalization of the Fenchel-Nielsen ``twist'' function of Teichm\"uller theory.
Finally, for every big edge $e$, we let $\pp_e,\pp'_e$ be the two paths running counter-clockwise around endpoints of $e$ on the small boundary (Figure \ref{fig:twist-e-setup}, right). It is easy to see that the functions $x_{\pp_e} = x_{\pp_e'}$ both reduce to a standard cross-ratio of the four framing lines surrounding the edge $e$. With orientation as on the right of Figure \ref{fig:twist-e-setup}, we have
\be x_e:= x_{\gamma_e}= x_{\pp(e)} = x_{\pp'(e)} = - \frac{\langle a\wedge b\rangle\langle c\wedge d\rangle}{\langle a\wedge c\rangle\langle b\wedge d\rangle}\,. \label{xe} \ee
This is a standard Fock-Goncharov cross-ratio function%
\footnote{This is the negative of the edge coordinate discussed in Appendix B of \cite{DGV-hybrid}, and agrees with the positive coordinates of \cite{FG-Teich}.}%
, a generalization of Thurston's shear coordinates in Teichm\"uller theory.

\begin{figure}[htb]
\centering
\includegraphics[width=6in]{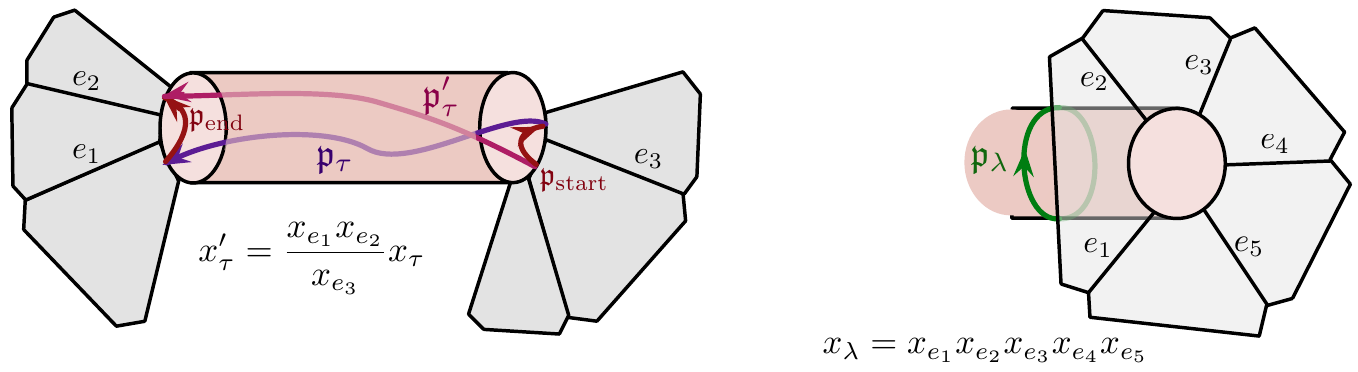}
\caption{Examples of the basic relations among $x_\lambda,x_\tau,x_e$: moving the endpoints of $\tau$ path (left); and expressing a holonomy as a product of edges (right).}
\label{fig:TLrels}
\end{figure}

The standard relations among functions $x_e,x_\tau,x_\lambda$, arising from $x_{\pp\circ\pp'}=x_\pp x_{\pp'}$, are illustrated in Figure \ref{fig:TLrels}. Moving the endpoints of a path $\pp_\tau$ results in multiplication of $x_\tau$ by $x_e$ functions; the product of $x_e$'s around any hole filled in by a small disc equals $1$ (consistent with the fact that holonomy around a small disc is trivial); and the product of $x_e$'s around the two ends of an annulus are inverses of each other, and equal to $x_\lambda^{\pm 1}$. From injectivity of \eqref{coord-inj} and the presentation \eqref{Hrels} of $H_1^-(\Sigma)$ in Section \ref{sec:path}, it follows that these are the only non-trivial relations. Indeed, we must have
\begin{align} \label{Pexplicit}
\CP[\CC] &=\Big\{(x_\alpha^{(t)},x_\beta^{(t)},x_\lambda^{(a)},x_\tau^{(a)},x_e) \,\in [(\C^*)^2\bs(1,1)]^{\# t}\times [(\C^*\bs (1))\times\C^*]^{\# a} \times (\C^*)^{\# e}\,\Big|  \\
 &\hspace{.5in}  \prod_{\text{$e$ on $d$}}x_e=1,\, \prod_{\text{$e$ on $\pd_1a$}} x_e = \prod_{\text{$e$ on $\pd_2a$}} x_e^{-1} = \lambda^{(a)} \,\Big\}\,.
\notag \end{align}

\subsection{Poisson brackets, symplectic form, and $K_2$ form}
\label{sec:Poisson}

\begin{figure}[htb]
%\centering
\hspace{-.25in}\includegraphics[width=6.5in]{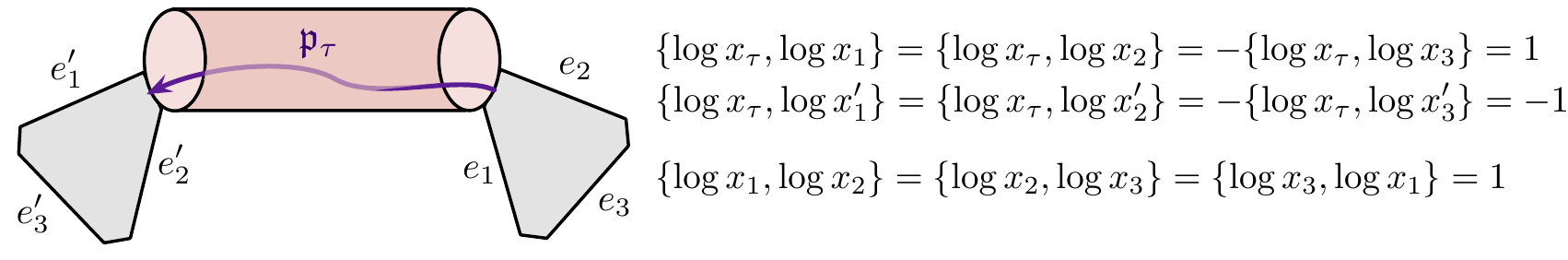}
\caption{Orientation convention for Poisson brackets of twist and edge coordinates. In the formulas, we abbreviate $x_{e_i} = x_i$, $x_{e_i'}=x_i'$, etc.}
\label{fig:Pet}
\end{figure}

When $M$ has no defects and $\pi_1(\CC_{\rm small})$ is abelian, we can calculate the Poisson bracket among functions $x_\gamma$ at generic points of $\CX[\CC]$ by pulling-back the Atiyah-Bott formula \eqref{AB}. In fact, all the relevant calculations have already been done in \cite{FG-Teich} and \cite[Appendix B]{GMNII}. (In both references, the idea was to calculate fundamental brackets involving the contractions of lines $\langle a\wedge b\rangle$ that appear in formulas such as \eqref{xt} and \eqref{xe}. Holonomy eigenvalues $x_\lambda$, etc., can also easily be written in terms of such contractions.)
We find
\be \label{PGL2PB}
\begin{array}{r@{\quad}c}
\text{for each small torus:} & \{\log x_\alpha,\log x_\beta\} = 2\,\langle \pp_\alpha,\pp_\beta\rangle =  2 \\
\text{for each small annulus:} & \{\log x_\tau,\log x_\lambda\} = 2\,\langle \pp_\tau,\pp_\lambda \rangle = \pm 2 \\
%\text{for edges $e$ adjacent to an annulus:}
 & \{\log x_\tau, \log x_e\} = \pm 1 \qquad \text{if $e$ on faces adjacent to $\pd\pp_\tau$}\\
\text{for big edges $e,e'$:} & \{\log x_e,\log x_{e'}\} = \text{\# faces shared by $e,e'$}\,,
\end{array}
\ee
with all other brackets vanishing.
Here $\langle \pp_\tau,\pp_\lambda \rangle$ denotes an intersection number of paths on an annulus (similarly on a torus), with orientation such that $\langle \pp_\tau,\pp_\lambda \rangle = +1$ on the left of Figure \ref{fig:twist-e-setup} (the paths there intersect on the back side of the annulus). For the $\{\log x_\tau,  \log x_e\}$ bracket, there are three edges bounding the face $f$ at the start of $\pp_\tau)$ and three bounding $f'$ the end of $\pp_\tau$, with signs of brackets given in Figure \ref{fig:Pet}. Finally, the number of faces in an edge-edge bracket is counted with orientation, such that $\{\log x_{e_1},\log x_{e_2}\}=+1$ in Figure \ref{fig:Pet}.

Although the brackets \eqref{PGL2PB} are calculated at points of $\CP[\CC] \subset\CX[\CC]$ that project to smooth points of the standard character variety, we will use \eqref{PGL2PB} to define them on the entire $(\C^*)^d$ patch $\CP[\CC]$. One then ought to verify that this extension is natural and consistent --- in particular that it agrees on overlaps of patches $\CP[\CC;\mb t_{\rm 2d}]$ labelled by \emph{different} big-boundary triangulations. The main ingredient of such a verification (as in \cite{FG-Teich}) is the fact that a flip of 2d triangulation acts as a symplectomorphism --- in fact, a $K_2$ symplectomorphism. We will not pursue the global nature of Poisson/symplectic structures on framed moduli spaces here. Strictly speaking we do not need it: we will confine ourselves to working with fixed patches and triangulations.

By direct comparison to \eqref{int-H} in Lemma \ref{lemma:intH}, we see that the Poisson brackets on $\CP[\CC]$ are all encoded in the intersection form in odd homology: $\{\log x_\gamma,\log x_{\gamma'}\} = \langle \gamma,\gamma'\rangle$ for any $\gamma,\gamma'\in H_1^-(\Sigma)$. The symplectic form on $\CP[\CC]$ can be expressed concretely by inverting the Poisson brackets,
\be \omega|_{\CP[\CC]} = \frac12\sum_{i,j} (\epsilon^{-1})^{ij}\, \frac{dx_{i}}{x_{i}}\wedge \frac{dx_{j}}{x_{j}}\,,\qquad \epsilon_{ij} := \langle \gamma_i,\gamma_j\rangle\,,\ee
where $\{\gamma_i\}$ is any basis for $H_1^-(\Sigma)$ and the $x_i:=x_{\gamma_i}$ are corresponding coordinates on $\CP[\CC]$. (Remember that the intersection form on odd homology is non-degenerate, so $\epsilon$ is invertible, though not necessarily over the integers.) Moreover, the symplectic form lifts to a $K_2$ avatar
\be \hat \omega := \frac12\sum_{ij}  (\epsilon^{-1})^{ij} x_i \wedge x_j \quad \in K_2(\CP[\CC])_\Q  \label{K2}\ee
in the second algebraic $K$-group of the field of functions on $\CP[\CC]$ (tensored%
\footnote{It follows from Lemma \ref{lemma:int} (Appendix \ref{app:basics}) that  $\epsilon^{-1}$ contains either integers or half-integers. Thus one could actually tensor with $[\tfrac12]$ (killing $2^p$ torsion) rather than with $\Q$. We don't keep careful track of this here.} with $\Q$). Such $K_2$ avatars are discussed from various perspectives in \cite{FG-Teich, DGG-Kdec, Dunfield-Mahler, Champ-hypA}. We summarize these observations as:

\begin{prop} \label{prop:PB}
If $M$ has no defects and $\pi_1(\CC)$ is abelian, so that the functions $x_\gamma$ ($\gamma\in H_1^-(\Sigma)\subset \wt H_1^-(\Sigma)$) contain complete coordinates on $\CP[\CC]$, the holomorphic symplectic form (agreeing with the Atiyah-Bott form), its $K_2$ avatar, and the holomorphic Poisson brackets are all encoded by the intersection product on $H_1^-(\Sigma)$. In particular,
\be \{x_\gamma,x_{\gamma'}\} = \langle \gamma,\gamma'\rangle\, x_\gamma x_{\gamma'}\qquad\text{or}\qquad \{\log x_\gamma,\log x_{\gamma'}\}=\langle \gamma,\gamma'\rangle\,. \label{PB}\ee
\end{prop}

\subsubsection{Example: tetrahedron}
\label{sec:PGLtet}

\begin{wrapfigure}{r}{2in}
\centering
\includegraphics[width=1.9in]{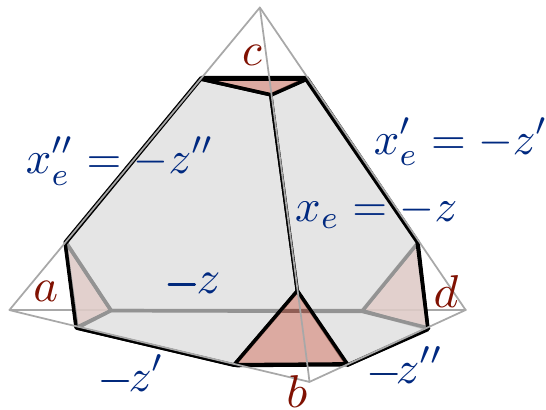}
\caption{Functions on $\CX[\CC_\Delta]$}
\label{fig:tetx}
\end{wrapfigure}

The punctured boundary of a tetrahedron $\CC_\Delta^*$ is a four-punctured sphere, and the space of framed flat $PGL(2)$ connections on $\CC_\Delta^*$ is two-dimensional. It has a canonical open subset $\CP[\CC_\Delta] = \CP[\CC_\Delta,\mb t_{\rm 2d}] \simeq \C^*\times\C^* \subset \CX[\Delta]$ corresponding to the canonical triangulation $\mb t_{2d}$ of the tetrahedron's big boundary \cite{DGG-Kdec}. 
Given the description of odd homology $H_1^-(\Sigma_\Delta)$ in Section \ref{sec:tetH}, Proposition \ref{prop:coords} says that $\CP[\CC_\Delta]$ is covered by six edge functions, subject to the relation that the product around every vertex is one. In terms of cross-ratios of framing lines, we have
\be \label{zcross}
x_e = - \frac{\langle a\wedge b\rangle\langle c\wedge d\rangle}{\langle a\wedge c\rangle\langle b\wedge d\rangle}\,,\quad
  x_e' = - \frac{\langle b\wedge d\rangle\langle c\wedge a\rangle}{\langle b\wedge c\rangle\langle d\wedge a\rangle}\,,\quad
  x_e''= - \frac{\langle d\wedge a\rangle\langle c\wedge b\rangle}{\langle d\wedge c\rangle\langle a\wedge b\rangle}\,,
\ee
with cross-ratios on opposite edges equal and $x_ex_e'x_e''=1$ (Figure \ref{fig:tetx}). Here $x_e=x_{\gamma}$, $x_e'=x_{\gamma'}$, etc. for cycles $\gamma,\gamma'\in \wt H_1^-(\Sigma_\Delta)$. The Poisson bracket and $K_2$ form in Prop. \ref{prop:PB} are
\be \{\log x_e,\log x_e'\} = \{\log x_e',\log x_e''\} = \{\log x_e'',\log x_e\}=1\,,\qquad \hat \omega = x_e\wedge x_e'\,.\ee
Setting
\be z = -x_e = x_{\gamma+u}\,,\quad z'=-x_e',\,\quad z'' = -x_e''\,, \label{zxapp}\ee
we recover the familiar phase space implicit in hyperbolic geometry \cite{Dimofte-QRS}.

In hyperbolic geometry there is usually an additional relation $z+z'{}^{-1}-1=0$. This arises when considering framed flat connections on $\CC_\Delta$ that can be extended to framed flat connections in the interior of the tetrahedron. Indeed, in this case the flat connection can be trivialized, and all three cross-ratios in \eqref{zcross} are computed in a single common fiber, leading to the additional relation. The submanifold $\CL_\Delta = \{z+z'{}^{-1}-1=0\} \subset \CP[\CC_\Delta]$ is a $K_2$ Lagrangian submanifold; it simply parameterizes the configuration space of four lines in $\C^2$. For further discussion, see \cite{Dimofte-QRS, DGG, DGG-Kdec}.

\subsection{Gluing $PGL(2)$ connections on boundaries}
\label{sec:gluePGL}

Our next task is to describe the gluing equations for framed $PGL(2)$ flat connections on boundaries. To do so, we first explain abstractly what it means to glue framed flat connections on boundaries.

Suppose that we glue framed 3-manifolds $M\leadsto M_0\leadsto M'$ in two steps, as in Section \ref{sec:glue-gen}. Assume that $M$ and $M'$ are defect-free, and that their respective small boundaries $\CC_{\rm small}$, $\CC_{\rm small}'$ have abelian fundamental groups.

There is always a map
\be g_{PGL(2)}^{(1)}\,:\; \CX[\CC] \to \CX[\CC_0^{**}]\,, \ee
Namely, given a framed flat connection $\CA$ on $\CC^*$, we may use a gauge transformation to trivialize it in the interiors of faces of $\mb t_{\rm 2d}$, and then (trivially) identify connection on pairs of faces that are glued together. We also identify framing flags on adjacent pieces of small boundary.
We must recall, however, that by definition $\CA$ is a framed flat connection on the punctured boundary $\CC^*=\CC\bs\{\text{point on each small disc}\}$ rather than on $\CC$ itself, so gluing produces a flat connection not on $\CC_0$ itself but on $\CC_0$ with some collection of points (punctures) removed from its small boundary, which we denote $\CC_0^{**}$. The connection is framed on $(\CC_0^{**})_{\rm small}$, and has unipotent holonomy around each puncture.

Next, in order to fill in the defects $I_j$ of $\CC_0^{**}$, it is necessary and sufficient to require that the holonomy $M_j$ around each defect is trivial. (If the holonomy is trivial, then connection can be trivialized along the defect, which allows the defect to be filled in.) Therefore, we get a map
\be g_{PGL(2)}^{(2)}\,:\; \CX[\CC_0^{**}]\big|_{M_j=I} \to \CX[\CC'{}^{**}]\,.\ee

\begin{figure}[htb]
\centering
\includegraphics[width=5.5in]{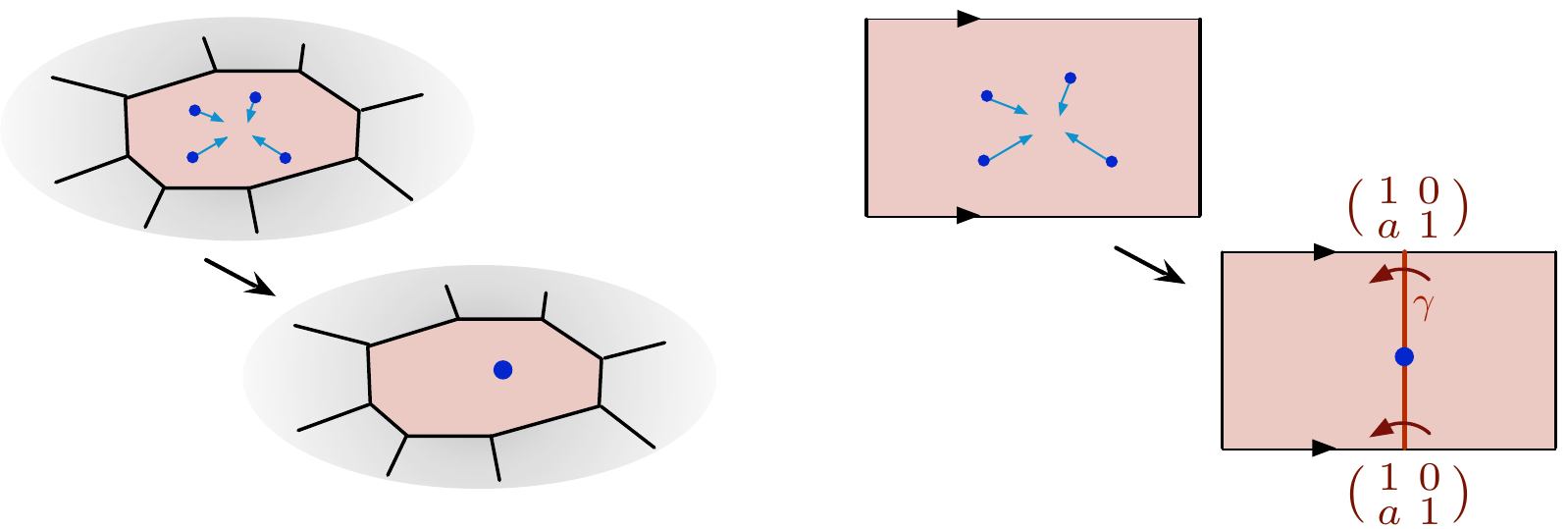}
\caption{Removing extraneous punctures on small boundaries. For a disc (left) we just collapse the punctures. For an annulus or torus (right) we collapse the punctures and modify the flat connection to trivialize the unipotent holonomy around the puncture. }
\label{fig:remp}
\end{figure}

Again, we observe that $\CC'{}^{**}$ may have extraneous punctures, arising from disjoint discs in $\CC_{\rm small}$ that have been connected in $\CC_{\rm small}'$. We can remove these punctures. For every small sphere component of $\CC_{\rm small}'$, we collapse all (potential) punctures to a single one; the holonomy around this single puncture is necessarily trivial, so we may remove it.
For every small disc, we similarly collapse all (potential) punctures to a single one, which still has unipotent holonomy (Figure \ref{fig:remp}, left). For every small annulus we modify the connection so that it extends over the punctures. There is a unique way to do this.
Abstractly, we note that since all puncture holonomies are unipotent and the connection is framed on the annulus (meaning, as usual, that it preserves a line) it makes sense to define an invariant holonomy eigenvalue $x_\lambda$ and a rescaling/twist coordinate $x_\tau$ from one end of the annulus to the other, just as in Section \ref{sec:PGLcoords}. In fact, we can calculate full (basepointed) holonomy matrices $M_\lambda,M_\lambda'$, both with eigenvalues $x_\lambda$, around the two ends of the annulus.
Then from Appendix \ref{app:traffic} we know that there's a unique framed flat connection on the smooth annulus with fixed $M_\lambda,M_\lambda',x_\tau$, as long as $x_\lambda\neq 1$.

To implement the modification concretely, we can collapse all punctures on an annulus to a single point, then cut the annulus along a curve $\gamma$ beginning and ending at this point as in Figure \ref{fig:remp} (right). We re-glue the connection on the annulus with an extra unipotent modification $\left(\begin{smallmatrix} 1 & 0 \\ a & 1\end{smallmatrix}\right)$ across $\gamma$. As long as $x_\lambda\neq 1$, there is a unique value of $a$ that trivializes the puncture holonomy. 

We may similarly remove extraneous punctures from small tori, so long as the invariant holonomy eigenvalues satisfy $(x_\alpha,x_\beta)\neq (1,1)$. Concretely, we collapse all punctures to a single one, then add a unipotent modification along a nontrivial curve $\gamma$ passing through the single puncture. The modification that trivializes the puncture holonomy is uniquely specified so long as the holonomy along $\gamma$ itself is not unipotent --- so, given $(x_\alpha,x_\beta)\neq (1,1)$, an appropriate $\gamma$ can always be chosen.

Altogether, we have built a removal-of-punctures map
\be g_{PGL(2)}^{(3)}\,:\; \CX[\CC'{}^{**}]\big|_{x_\lambda\neq 1,\,(x_\alpha,x_\beta)\neq (1,1)} \to \CX[\CC']\,,\ee
and composing it with the first two parts of the gluing procedure we obtain a gluing map
\be \label{PGLglue}
\boxed{g_{PGL(2)} = g_{PGL(2)}^{(3)}\circ g_{PGL(2)}^{(2)}\circ g_{PGL(2)}^{(1)}\,:\; \CX[\CC]\big|_{R} \to \CX[\CC']}\,.
\ee
Here the restrictions `$R$' are a lift to $\CX[\CC]$ of the defect conditions $M_j=1$ and puncture-removal conditions $x_\lambda\neq 1,\,(x_\alpha,x_\beta)\neq (1,1)$.

\subsection{Symplectic properties of the gluing equations}
\label{sec:glueX}

Let's now restrict to patches $\CP[\CC,\mb t_{\rm 2d}]$, $\CP[\CC',\mb t_{\rm 2d}']$, corresponding to the triangulation of $M$ and the induced triangulation of $M'$, and describe concretely how the gluing map \eqref{PGLglue} acts on $\C^*$ coordinates.

Consider the intermediate space $\CX[\CC_0^{**}]$, where $\CC_0^{**}$ is the boundary of $M_0$, possibly with some additional unipotent punctures. This space supports the same path-functions $x_\pp$, $\pp\in \mb P_0\oplus \Z_2$ as $\CX[\CC_0]$  (since unipotent punctures do not affect the definition or properties of the $x_\pp$), which depend only on the homology class $\gamma\in \wt K= \im\,\big[\wt h:(\mb P_0\oplus \Z_2)\to \wt H_1^-(\Sigma)\big]$. As usual, we work with the convention that $x_u=-1$ for the fiber class $u\in \Z_2$. The functions $x_\gamma$, $\gamma\in \wt K$ take well-defined $\C^*$ values on a subset $\CP[\CC_0^{**},\mb t_{\rm 2d}'] \subset \CX[\CC_0^{**}]$ that is defined the same way as in Def. \ref{def:CP}.

The image%
\footnote{Technically, we should relax the non-unipotent restriction in the definition of $\CP[\CC_0^{**}]$ in order for the image of $g_{PGL}^{(1)}:\CP[\CC] \to \CX[\CC_0^{**}]$ to be fully contained in $\CP[\CC_0^{**}]$. This does not affect the following argument.} %
of $g_{PGL}^{(1)}:\CP[\CC;\mb t_{\rm 2d}] \to \CX[\CC_0^{**}]$ lies in $\CP[\CC_0^{**};\mb t_{\rm 2d}']$. Moreover,
\be g_{PGL}^{(1)}{}^* (x_\gamma) = x_{\tilde g(\gamma)}\,,\qquad  \gamma\in \wt K\,, \label{xPGL1} \ee
with $\wt g:\wt H_1^-(\Sigma_0)\to \wt H_1^-(\Sigma)$ the map of \eqref{gluetw} and \eqref{gPh}. To understand this, suppose $\gamma=\wt h(\pp)$ for a path $\pp\in \mb P_0$. We can compute the path coordinate $x_\gamma=x_\pp$ by breaking $\pp$ up into segments $g_{\mb P}(\pp)=\sum_i\pp_i$, $\pp_i\in \mb P$ (as in Figure \ref{fig:defgP}), applying the definition \eqref{xt} to each segment to get functions $x_{\pp_i}$, and multiplying them together. Successive numerators and denominators of \eqref{xt} cancel out --- just as in the proof of the concatenation relation $x_{\pp\circ\pp'}=x_\pp x_{\pp'}$ --- up to a single sign $(-1)$. This extra sign arises due to the reversed relative orientation at the head of one segment and the tail of the next. Thus $x_\pp = ...x_{\pp_{i+1}}(-1)x_{\pp_i}(-1)x_{\pp_{i-1}}...$. The modifications by these signs are precisely encoded in the fiber-class corrections to the twisted cutting map $\wt g_{\mb P}$ \eqref{gPh}. Using $\wt h\wt g_{\mb P}=\wt g\wt h$, we obtain \eqref{xPGL1}.

The gluing conditions requiring trivial defect holonomies $M_j$ simply say that $x_{\mu_j}\equiv 1$ for all defect cycles $\mu_j \in \wt G \subset \wt K$. The lifts of these conditions to $\CP[\CC,\mb t_{\rm 2d}]$ are $x_{\tilde g(\mu_j)}\equiv 1$. If these conditions are satisfied, we get a map
\be g_{PGL(2)}^{(3)}\circ g_{PGL(2)}^{(2)}:\, \CP[\CC_0^{**},\mb t_{\rm 2d}']\big|_{x_\mu\equiv 1} \to \CP[\CC',\mb t_{\rm 2d}']\,.  \ee
Now the functions $x_\gamma$ ($\gamma\in \wt H_1^-(\Sigma')$) on $\CP[\CC',\mb t_{\rm 2d}']$ pull back%
\footnote{Again, note that the unipotent modifications made by $g_{PGL(2)}^{(3)}$ to remove extraneous punctures do not affect path coordinates. So to understand \eqref{xPGL2} it suffices to look at the action of $g_{PGL(2)}^{(2)}$, which just fills in defects.} %
to
\be \big(g_{PGL(2)}^{(3)}\circ g_{PGL(2)}^{(2)}\big)^*(x_\gamma) = x_{\hat\gamma}\,,\qquad \hat\gamma\in \wt K\,,\; \wt q(\hat\gamma)=\gamma\,, \label{xPGL2} \ee
where $\hat\gamma \in \wt K\subset \wt H_1^-(\Sigma_0)$ is any preimage of $\gamma$ under the map $\wt q$ from \eqref{PKG}. This preimage is unique modulo $\wt G$, so $x_{\hat\gamma}$ is well defined on $\CP[\CC_0^{**},\mb t_{\rm 2d}']\big|_{x_\mu\equiv 1}$.

By combining \eqref{xPGL1} and \eqref{xPGL2} we see that for any $\gamma\in \wt K\subset \wt H_1^-(\Sigma_0)$, the path-functions on $\CP[\CC,\mb t_{\rm 2d}]\big|_{M_j=1}$ and  $\CP[\CC',\mb t_{\rm 2d}']$, evaluated on a connection $\CA$ and its glued image $g_{PGL(2)}(\CA)$, are related by
\be \boxed{x_{\tilde g(\gamma)} = x_{\tilde q(\gamma)}\,,\qquad \forall\; \gamma\in \wt K\subset \wt H_1^-(\Sigma_0)\,.} \label{glue-simple} \ee
These are the gluing equations. They subsume the gluing conditions $x_{\tilde g(\mu)}\equiv 1$ for $\mu\in \wt G$, since $\wt q(\mu)\equiv 0$. Moreover, 
since $\wt q$ is surjective and the path-functions on $\CP[\CC',\mb t_{\rm 2d}']$ include a complete set of coordinates, the gluing map
\be g_{PGL(2)}:\, \CP[\CC,\mb t_{\rm 2d}]\big|\raisebox{-.1cm}{$(x_{\tilde g(\mu)}\equiv 1,R')$} \to\hspace{-.3cm}\to \CP[\CC',\mb t_{\rm 2d}'] \label{gsurj} \ee
must be surjective. (Here $R'$ denotes the additional non-unipotent restrictions lifted from $\CP[\CC',\mb t_{\rm 2d}']$, namely $x_{\tilde g(\hat \lambda)}\neq 1,\, (x_{\tilde g(\hat \alpha)},x_{\tilde g(\hat \beta)})\neq (1,1)$, for $\wt q(\hat\lambda,\hat\alpha,\hat\beta)=(\lambda,\alpha,\beta)$.) Indeed,

\begin{thm} \label{thm:H}

The $PGL(2)$ gluing map \eqref{gsurj} is a (holomorphic) $K_2$ symplectic reduction of a finite quotient,
\begin{align} g_{PGL(2)}&:\; (\CP[\CC;\mb t_{\rm 2d}]_{R'}/Z)\big/\!\!\big/ (\C^*)^{{\rm rank}\, G} := (\CP[\CC;\mb t_{\rm 2d}]_{R'}/Z)\big|_{x_{\tilde g(\mu)}\equiv 1}\big/(\C^*)^{{\rm rank}\, G} \notag \\
 &\hspace{3.5in}\;\overset\sim\to\; \CP[\CC';\mb t_{\rm 2d}']\,,
\label{red}
\end{align}
where $Z\simeq \wt H_1^-(\Sigma)/\wt H$ is the torsion group of Prop. \ref{prop:glue}'(b), and the $(\C^*)^{{\rm rank}\, G}$ action is generated by using the $x_{\tilde g(\mu_j)}$, $\mu_j\in \wt G$, as moment maps.

\end{thm}

Explicitly, if $\wt H\subset \wt H_1^-(\Sigma)$ is the finite-index subgroup containing $g(\wt K)$, such that $\wt H_1^-(\Sigma') \simeq \wt H/\!/\wt g(\wt G) = \wt g(\wt K)/\wt g(\wt G)$, we recall that $Z=\wt H_1^-(\Sigma)/\wt H$ contains at most 2-torsion and 4-torsion. The group $Z$ act naturally on $\CP[\CC;\mb t_{\rm 2d}]$ by multiplying functions $x_\gamma$ by 4-th roots of unity. Namely, letting $\{\sigma_i\}$ be any generators of the ${\rm Hom}\,\big(Z,\Z_4=\{1,i,-1,-i\}\big)$, the action on $\CP[\CC;\mb t_{\rm 2d}]$ is generated by $x_\gamma\mapsto \sigma_i[\gamma]x_\gamma$.

We may also explicitly describe the action of the moment maps as follows. Let $\{\mu_i\}_{i=1}^{{\rm rank}\,G}$ be a basis for $\wt G$ and let $(t_i)_{i=1}^{{\rm rank}\,G}\in (\C^*)^{{\rm rank}\,G}$ be some corresponding parameters. If $\mu = \sum a_i\mu_i\in \wt G$ define $t^\mu=\prod_i t_i^{a_i}$. Then, by virtue of \eqref{PB} in Proposition \ref{prop:PB}, the moment map action of $x_{\tilde g(\mu)}$ is
\be x_\gamma \mapsto (t^\mu)^{\langle \mu,\gamma\rangle} x_\gamma\,,\qquad \forall\; \gamma\in \wt H_1^-(\Sigma)\,,\; \mu\in \wt G\,. \label{mom} \ee

\medskip
\noindent{\it Proof of Theorem \ref{thm:H}}. We analyze the quotient on the LHS of \eqref{red}. By construction, the invariant functions on the ordinary quotient $\CP[\CC,\mb t_{\rm 2d}]/Z$ are precisely $x_\gamma$ for $\gamma\in \wt H\subset \wt H_1^-(\Sigma)$. These include a complete set of coordinates. Since $g(\wt K) = \ker\,\langle \wt g(\wt G),*\rangle|_{\tilde H}$ by Proposition \ref{prop:glue}' and Lemma \ref{lemma:SESP}, the  functions that are further invariant under the $(\C^*)^{{\rm rank}\,G}$ action \eqref{mom} are labeled by $\gamma\in \wt g(\wt K)\subset  \wt H_1^-(\Sigma)$. These functions coordinatize $(\CP[\CC;\mb t_{\rm 2d}]/Z)\big/(\C^*)^{{\rm rank}\, G}$. Subsequently restricting to $x_{\tilde g(\mu)}\equiv 1$ for all $\mu\in \wt G$ produces a space coordinatized by $x_{[\gamma]}$ for $[\gamma]\in \wt g(\wt K)/\wt g(\wt G)\simeq \wt H_1^-(\Sigma')$. The identification $\wt g(\wt K)/\wt g(\wt G)\simeq \wt H_1^-(\Sigma')$ sends $[\gamma]$ to $\wt q\circ \wt g^{-1}(\gamma)$ and preserves the intersection form. Therefore, $(\CP[\CC;\mb t_{\rm 2d}]_{R'}/Z)\big/(\C^*)^{{\rm rank}\, G}\big|_{x_{\tilde g(\mu)}\equiv 1}$ is canonically (1-1) symplectomorphic to $\CP[\CC';\mb t_{\rm 2d}']$, with path-coordinates related as $x_{\tilde g(\gamma)}=x_{\tilde q(\gamma)}$ for $\gamma\in \wt K$. These are precisely the gluing equations \eqref{gsurj}, so the symplectomorphism is the desired gluing map.

To see that this is a $K_2$ symplectomorphism, first note that $\CP[\CC;\mb t_{\rm 2d}]$ and $\CP[\CC;\mb t_{\rm 2d}]/Z$ have the same $K_2$ forms $\hat \omega \in K_2(\CP[\CC])_\Q$, given by \eqref{K2}. Let us choose a basis $\{\gamma_i,\mu_j\}$ for the untwisted subgroup $K\subset \wt K$ such that  $\{\wt q(\gamma_i)\}$ is a basis for $H_1^-(\Sigma')\subset \wt H_1^-(\Sigma')$. Complete this to a basis $\{\wt g(\gamma_i),\wt g(\mu_j),\eta_k\}$ of a torsion-free (untwisted) subgroup $H\subset \wt H$, so that the intersection form in this basis is block-diagonal,
\be \langle \wt g(\gamma_i),\wt g(\gamma_{i'})\rangle = \langle\wt q(\gamma_i),\wt q(\gamma_{i'})\rangle = \epsilon'_{ii'}\,,\qquad \langle \wt g(\mu_j),\eta_k\rangle = c_j\,\delta_{jk}\quad (c_j\in \Z)\,, \ee
where $\epsilon'$ is the intersection form on $H_1(\Sigma')$ and all other intersection products vanish. (The fact that we can find such a basis follows from the untwisted version of $\wt H/\!/\wt g(\wt G)= \wt g(\wt K)/\wt g(\wt G)\simeq \wt H_1^-(\Sigma')$.) Then the $K_2$ form on $\CP[\CC;\mb t_{\rm 2d}]/Z$ is
\be \hat\omega  = \frac12\,\sum_{i,i'} (\epsilon'{}^{-1})^{ii'} x_{\tilde g(\gamma_i)}\wedge x_{\tilde g_(\gamma_i')} + \sum_j \frac1{c_j} x_{\tilde g(\mu_j)}\wedge x_{\eta_j}\,. \ee
Setting $x_{\tilde g(\mu_j)}= 1$ and using the fact that $1\wedge *=0$ in K-theory, we find that $\hat\omega|_{(x_{\tilde g(\mu)}\equiv 1)} = \hat\omega'$ reproduces the $K_2$ form on $\CP[\CC',\mb t_{\rm 2d}']$. \; $\square$

\subsection{Example: Thurston's gluing equations, $K_2$ forms, and volumes}
\label{sec:Thurs}

In the Introduction (Section \ref{sec:intro-eg}), we claimed that Thurston's gluing equations are a special case of \eqref{glue-simple}, and therefore that Theorem \ref{thm:H} directly implies the symplectic properties found by Neumann and Zagier. We take a moment to spell out exactly how this works.

Suppose that $M'$ is an oriented hyperbolic 3-manifold with an ideal triangulation $M'=\cup_{i=1}^N\Delta_i$. We denote by $M=\sqcup_{i=1}^N \Delta_i$ the disjoint union of tetrahedra in the triangulation, viewing both $M$ and $M'$ as framed 3-manifolds. For example, $M'$ might be a closed hyperbolic manifold with a spun triangulation, in which case the topological ideal boundary $\CC'=\pd M'$ consists entirely of small spheres; or $M'$ could have $n_c$ cusps, corresponding to small torus components of $\CC'$.
One might also consider $M'$ with geodesic boundary, in which case $\CC'_{\rm big}$ is non-empty. This is another special case of our general framework, but we'll ignore it for the moment to keep this example simple. Thus, $\CC_{\rm big}'$ is empty, and $\CC_{\rm small}'$ contains $n_c\geq 0$ small tori and some number of small spheres. The number of edges $I_j$ of the triangulation is the same as the number $N$ of tetrahedra. We glue $M\leadsto M_0\leadsto M'$ as in Section \ref{sec:glue-gen}, letting $M_0$ be a framed 3-manifold with $N$ defects $I_j$.

For each tetrahedron $\Delta_i$ in $M$ we find a $K_2$ symplectic space $\CP[\CC_{\Delta i}] = \{x_i,x_i',x_i''\,|\,x_ix_i'x_i''=1\}\simeq (\C^*)^2$ as in Section \ref{sec:PGLtet}. Here $x_i=x_{\gamma_i},x_i'=x_{\gamma_i'}$, etc., with $\gamma_i,\gamma_i',\gamma_i''$ and $u$ generating $\wt H_1^-(\Sigma_\Delta)$ as in Section \ref{sec:tetH}. The triple of hyperbolic shapes for each tetrahedron are $(z_i,z_i',z_i'')=(-x_i,-x_i',-x_i'')$.
Altogether,
\be \CP[\CC]=\prod_{i=1}^N \CP[\CC_{\Delta i}] \simeq (\C^*)^{2N}\,, \qquad \wt H_1^-(\Sigma) = \oplus_{i=1}^N H_1^-(\Sigma_\Delta)\oplus \Z_2 \simeq \Z^{2N}\oplus \Z_2\,, \ee
where the splitting of $\wt H_1^-(\Sigma)$ comes naturally from Lemma \ref{lemma:gen}.
Similarly, for each small torus $T^2_t$ in $\CC_{\rm small}'$, there is an odd homology group $\wt H_1^-(\Sigma_{T^2})$ generated by $(\alpha^{(t)},\beta^{(t)},u)$ and a moduli space $\CP[\CC_{T^2}]=\{x_\alpha^{(t)},x_\beta^{(t)}\}\simeq (\C^*)^2\bs(1,1)$. Thus
\be \CP[\CC'] =\prod_{t=1}^{n_c} \CP[\CC_{T^2_t}] \simeq (\C^*)^{2n_c}\,,\qquad \wt H_1^-(\Sigma') =\oplus_{i=1}^N H_1^-(\Sigma_{T^2_t})\oplus\Z_2 \simeq \Z^{2n_c}\oplus\Z_2\,.\ee

Now consider the intermediate manifold $M_0$. The boundary $\CC_0=\pd M_0$ small and defect parts; $N$ defects end at $2N$ holes on $(\CC_0)_{\rm small}$. The holes on $(\CC_0)_{\rm small}$ lie at the vertices of the triangular tiling of $\CC_{\rm small}'$, as illustrated back in \ref{fig:intro-cycles}.  The subgroup $\wt G = \im\,[\wt h:\mb P_G \to \wt H_1^-(\Sigma_0)]$ is generated by cycles $\mu_j$ (Figure \ref{fig:defmu}), one for each defect $I_j$, coming from paths $\pp_{\mu_j}$ that surround the holes of $(\CC_0)_{\rm small}$. The subgroup $\wt K = \im\,[\wt h:(\mb P_0\oplus \Z_2)\to \wt H_1^-(\Sigma_0)]$ is generated by the $\mu_j$, together with lifts $\hat\alpha^{(t)},\hat\beta^{(t)}$ (Section \ref{sec:KG}) of $\alpha^{(t)},\beta^{(t)}\in H_1^-(\Sigma')$, and the fiber class $u$. Concretely, $\hat \alpha^{(t)}=\wt h(\hat\pp_\alpha^{(t)})$ and $\hat\beta^{(t)}=\wt h(\hat\pp_{\beta}^{(t)})$, where $\hat\pp_\alpha^{(t)},\hat\pp_{\beta}^{(t)}$ are paths on $\CC_{\rm small}$ representing A and B cycles, passing in some chosen way around the holes. These are exactly the types of boundary paths that appeared in \cite{NZ, Neumann-combinatorics}, said to be in normal position with respect to the tiling of $\CC_{\rm small}'$.

The first part of Thurston's gluing equations state that the product of tetrahedron shapes $z_i,z_i',z_i''$ around any edge $I_j$ must equal one. These are precisely our trivial-holonomy constraints, of the form
\be x_{\tilde g(\mu_j)} = 1\,,\qquad \forall\; \mu_j\in \wt G\,. \ee
To see this, note that under the cutting map $g_{\mb P}$ (Section \ref{sec:KG}) a path $\pp_{\mu_j}$ surrounding defect $I_j$ is cut into paths $\pp_e$ associated to all the edges $e$ of tetrahedra identified with $I_j$ in the gluing. The extra modification in the twisted map $\wt g_{\mb P}$ adds a fiber class $u$ to $\wt g_{\mb P}(\pp_{\mu_j})$ for every cut that is made; thus from $\wt h\wt g_{\mb P}=\wt g\wt h$ we find
\be \wt g(\mu_j) = \sum_{\text{around $I_j$}} (\gamma_i+u\;\text{or}\;\gamma_i'+u\;\;\text{or}\;\;\gamma_i''+u)\,, \notag \ee
\be  x_{\tilde g(\mu_j)} = \prod_{\text{around $I_j$}} (-x_i,-x_i',-x_i'') = \prod_{\text{around $I_j$}} (z_i,z_i',z_i'')\,,
\ee
matching the hyperbolic gluing equations.

\begin{figure}[htb]
\centering
\includegraphics[width=4in]{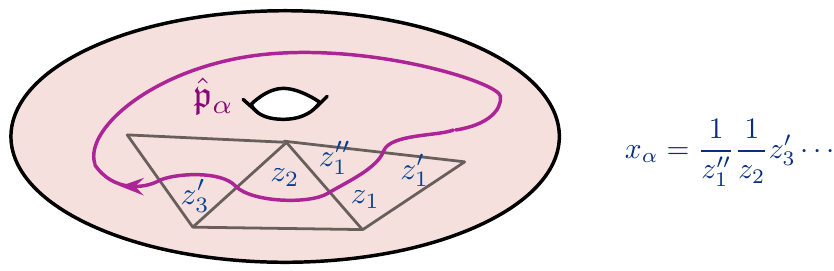}
\caption{Calculating the cusp equation for $\ell^2=x_\alpha^{(t)}$ on a small boundary that's tiled by small vertex triangles of truncated tetrahedra.}
\label{fig:LM}
\end{figure}

The second part of the gluing equations, sometimes called the ``cusp equations,'' states that the squares of A and B cycle eigenvalues $\ell_t^2, m_t^2$ (\ie\ eigenvalues of the hyperbolic holonomy) around each cusp equal the product of shapes $z_i^{\pm 1},z_i'{}^{\pm 1},z_i'{}^{\pm 1}$ at dihedral angles subtended by paths on $\CC_{\rm small}$, representing the respective cycles and in normal position with respect to the tiling of $\CC_{\rm small}$ --- in other words, our paths $\hat\pp_\alpha^{(t)}$ and $\hat\pp_\beta^{(t)}$. The exponents $\pm1$ correspond to whether an angle is subtended clockwise or counterclockwise, \cf\ Figure \ref{fig:LM}. These cusp equations again are just the remaining equations of the form $x_{\tilde q(\gamma)}=x_{\tilde g(\gamma)}$ in \eqref{glue-simple}, for $\gamma=\hat\alpha^{(t)}$ or $\gamma=\hat\alpha^{(t)}$. Indeed, on one hand,
\be \ell_t^2 = x_\alpha^{(t)} = x_{\tilde q(\hat \alpha^{(t)})}\,,\qquad m_t^2 = x_\beta^{(t)} = x_{\tilde q(\hat \beta^{(t)})}\,. \ee
On the other hand $\tilde g(\hat \alpha^{(t)}) = \wt h\circ \wt g_{\mb P}(\hat \pp_\alpha^{(t)})$ is a sum of cycles $\pm(\gamma_i+u)$, $\pm(\gamma_i'+u)$, $\pm(\gamma_i''+u)$ corresponding to angles subtended by $\hat \pp_\alpha^{(t)}$ (similarly for $\hat \beta^{(t)}$); so $x_{\tilde g(\hat \alpha^{(t)})}$ and $x_{\tilde g(\hat \beta^{(t)})}$ are precisely the desired products of shapes.

Now Theorem \ref{thm:H} amounts to the statement that the defect functions $x_{\tilde g(\mu_j)}$ are a set of mutually commuting moment maps on $\CP[\CC]$, which also commute with the A and B cycle functions $x_{\tilde g(\hat \alpha^{(t)})}$ and $x_{\tilde g(\hat \beta^{(t)})}$. Thus $\CP[\CC']$ is the symplectic reduction of (a finite quotient of) $\CP[\CC]$. Obviously this implies that there are $\tfrac12(\dim_\C\,\CP[\CC]-\dim_\C\,\CP[\CC'])=N-n_c$ independent moment maps. Homologically, this follows from the easy fact that
\be {\rm rank}\,\wt G = {\rm rank}\, G = N-n_c\,. \ee
Indeed, the sum of generators $\mu_j$ for all defects that end on a given cusp (counted with multiplicity) is null-homologous, so each cusp produces one relation among the $N$ generators.%
\footnote{It is a short exercise to show that these relations are all independent. See \cite{Neumann-combinatorics} or the review \cite{GF-review}.} %

Alternatively, one could write the gluing equations in a fixed basis as
\be \begin{array}{rcl}  \ell_t^2 &=& \pm\prod_i z_i^{A_{ti}} z_i'{}^{A_{ti}'} \\
m_t^2 &=& \pm\prod_i z_i^{B_{ti}} z_i'{}^{B_{ti}'} \\
1=x_{\tilde q(\mu_j)} &=& \pm\prod_i z_i^{C_{ji}} z_i'{}^{C_{ji}'} \end{array}\,,\qquad
 g\begin{pmatrix} \alpha \\ \beta\\ \mu \end{pmatrix} = \begin{pmatrix} A & A'\\ B & B' \\ C & C' \end{pmatrix} \begin{pmatrix} \gamma \\ \gamma' \end{pmatrix}\,.
\ee
Then the matrix of the untwisted map $g|_K=p\circ \wt g_K$ shown here has rank ${\rm rank}\,K=N+n_c$ and preserves the intersection form, \ie\ $(g|_K)J_{2N}(g|_K)^T = 2J_{2n_c}\oplus 0_{N\times N}$, as in the Introduction.

The reduction of $K_2$ forms is expressed as
\be \label{K2knot}
\hat\omega' \;=\; \sum_{t=1}^{n_c} \frac12\, \ell^2_t\wedge m^2_t \;=\; \sum_{i=1}^N z_i\wedge z_i'\,\Big|\raisebox{-.1cm}{$(x_{\tilde g(\mu)}\equiv 0)$}
 \;=\; \hat\omega\big|\raisebox{-.1cm}{$(x_{\tilde g(\mu)}\equiv 0)$}\,. \ee
There is a standard map $\eta: K_2(\CP[\CC'])\to \Omega^1(\CP[\CC'])$\, defined by $\eta(a\wedge b)= \log|a|\,d\arg b-\log|b|\,d\arg a$ (see \cite{DGG-Kdec} and references therin), which provides a canonical anti-derivative of the symplectic form on $\CP[\CC']$, since $d\eta(\hat\omega') = {\rm Im}\,\omega'$. By computing $\eta(\hat \omega')$ using both sides of \eqref{K2knot} and further restricting to the $K_2$ Lagrangian $\CL \subset \CP[\CC']$ defined by $z_i+z_i^{-1}-1=0$ together with the gluing equations,%
\footnote{This ``gluing variety'' generalizes the geometric component of the A-polynomial for a knot \cite{NZ, Dunfield-Mahler, Champ-hypA}.} %
we directly recover the Neumann-Zagier formula for variation of the volume,
\be \label{dVol} d\,{\rm Vol}(M') = \eta(\hat\omega')\big|_\CL\,.\ee

\section{Non-abelianization}
\label{sec:NA}

We showed in Section \ref{sec:PGL} that coordinates for framed flat $PGL(2)$ connections on the boundaries of framed 3-manifolds can be labelled by elements of the homology of double covers, with Poisson brackets of the former matching the intersection product of the latter. Moreover, we found that gluing equations \eqref{glue-simple} have a homological interpretation that makes their symplectic properties manifest. However, we did not explain \emph{why} the relation between $PGL(2)$ connections and homology of double covers existed, or why it was particularly natural. We now aim to fill this gap.

We will first consider the moduli space of flat $GL(1):=GL(1,\C)$ connections on a double cover $\Sigma$. This very simple space naturally has holonomy coordinates labelled by the first homology of $\Sigma$ (viewed as an abelianization of $\pi_1(\Sigma)$), with Atiyah-Bott Poisson bracket induced by the intersection form. The gluing equations for flat $GL(1)$ connections manifestly take the form \eqref{glue-simple}.

Then we borrow and extend the non-abelianization construction of \cite{GMN-spectral} to build a non-abelianization map $\Phi$, a (nontrivial) symplectomorphism between flat $GL(1)$ flat connections on a double cover $\Sigma$ and framed flat $PGL(2)$ connections on the base $\CC$. While the non-abelianization map of \cite{GMN-spectral} was mainly discussed in a rich geometric context --- involving a choice of complex structure on $\CC$ and an interpretation of $\Sigma$ as a spectral cover --- we will simply used the topological structure of the boundary of a framed 3-manifold (and a choice of big-boundary triangulation) to define $\Phi$.%
\footnote{\label{foot:HN}%
Topological descriptions of non-abelianization were also discussed in \cite{GMN-spectral, GMN-snakes} and \cite{HN-FN}, which overlap with our constructions. In our language, \cite{GMN-spectral, GMN-snakes} considered moduli spaces of $PGL(K)$ (not just $PGL(2)$) connections on surfaces $\CC$ consisting of an arbitrary big part $\CC_{\rm big}$, with 2d ideal triangulation $\mb t_{\rm 2d}$, and all holes of $\CC_{\rm big}$ filled by small discs. It was found for $PGL(2)$ that coordinates induced by non-abelianization coincide with Fock-Goncharov cluster coordinates. \cite{HN-FN} generalized the $PGL(2)$ non-abelianization map to construct complex Fenchel-Nielsen coordinates for a surface $\CC$ together with a pants decomposition. Such coordinates arise for us when $\CC$ (as a framed boundary) contains big 3-holed spheres connected by small annuli, \cf\ \cite{DGV-hybrid}. Our definition of $\Phi$ in this case differs from that of \cite{HN-FN}, but is ultimately equivalent.}

\subsection{Abelian flat connections}
\label{sec:ab}

We begin by defining a space of abelian flat connections whose coordinates are manifestly labelled by elements of twisted homology $\wt H_1^-(\Sigma)$. As observed in \cite{GMN-spectral} (and hinted in \cite{FG-Teich}), the non-abelianization construction requires such twisting.

For a closed oriented surface $\Sigma$, define
\begin{align} \label{defPab}
\wt \CX_{\rm ab}[\Sigma] &= \{\text{twisted $GL(1)$ flat connections on $\Sigma$}\} \\
 &:= \{\text{$GL(1)$ flat connections on $T_1\Sigma$ with fiber holonomy $-1$}\}\,. 
\notag
\end{align}
These are flat connections on a (necessarily trivial) complex line bundle $L\to T_1\Sigma$, or, equivalently, $GL(1):=GL(1,\C)$ local systems on $T_1\Sigma$.
The space $\wt \CX_{\rm ab}[\Sigma]$ is isomorphic to $(\C^*)^{{\rm rank}\, \tilde H_1(\Sigma)} = (\C^*)^{{\rm rank}\,H_1(\Sigma)}$. Indeed, a flat $GL(1)=GL(1,\C)$ connection on $T_1\Sigma$ is uniquely parametrized by its $GL(1)\simeq \C^*$ -valued holonomies $x_\gamma$ for $\gamma\in H_1(T_1\Sigma)$, with $x_{\gamma+\gamma'}=x_\gamma x_{\gamma'}$.
We are requiring that the fiber holonomy is $x_u=-1$, whence the holonomies $x_\gamma$ naturally become labelled by elements of twisted homology $\wt H_1(\Sigma)$.

Letting $\CX_{\rm ab}[\Sigma]$ denote the space of standard (untwisted) $GL(1)$ flat connections on $\Sigma$, we note that there is an isomorphism $\wt \CX_{\rm ab}[\Sigma]\simeq \CX_{\rm ab}[\Sigma]$ induced by any splitting $\wt H_1(\Sigma)\simeq H_1(\Sigma)\oplus \Z_2$. Such a splitting is given, for example, by the structure of a framed 3-manifold, with $\Sigma$ the canonical cover of the boundary (Lemma \ref{lemma:gen}). Equivalently, a choice of spin structure on $\Sigma$ induces an isomorphism $\wt \CX_{\rm ab}[\Sigma]\simeq \CX_{\rm ab}[\Sigma]$. Namely, a spin structure defines a (fiber-wise) 2-fold cover of $T_1\Sigma$; the pull-back of a twisted flat connection to the cover gives a connection with fiber holonomy $+1$, which may subsequently be pushed forward to $\Sigma$ itself, providing the desired isomorphism (\cf\ \cite[Sec. 10]{GMN-spectral}).

Now suppose $\Sigma\overset\pi\to\CC$ is the canonical double cover of the boundary of a framed 3-manifold. We want an odd version of $\wt \CX_{\rm ab}[\Sigma]$ whose coordinates are labelled by $\gamma\in \wt H_1^-(\Sigma)$. It can be defined as a projectivization, in the following sense. 
Let $\{\alpha_i\}_{i=1}^r$ be a basis for $H_1(\CC)$, and for an $r$-tuple of parameters $t = (t_1,..,t_r)\in (\C^*)^r$ and $\alpha = \sum_i a_i\alpha_i$ in $H_1(\CC)$ let $t^\alpha := t_1^{\alpha_1}\cdots t_r^{\alpha_r}$. We set
\be \label{defPab-}
\wt \CX_{\rm ab}^-[\Sigma] := \wt \CX_{\rm ab}[\Sigma]/(\C^*)^r\,,\qquad \text{with action}\;\; x_\gamma \mapsto t^{\pi_*\circ\, p(\gamma)}x_\gamma\,,\quad \gamma\in \wt H_1(\Sigma)\,,\; t\in (\C^*)^r\,.\ee
Here $p(\gamma)\in H_1(\Sigma)$ is the projection to standard homology, and $\pi_*\circ p(\gamma)\in H_1(\CC)$ is the subsequent projection to the base. The coordinates invariant under the $(\C^*)^r$ action are precisely those $x_\gamma$ with $\gamma\in \wt H_1^-(\Sigma)$. Thus the holonomies provide a map
\be \label{defpairab} x:\,\wt \CX_{\rm ab}^-[\Sigma] \times \wt H_1^-(\Sigma) \to \C^*\,, \ee
which is a homomorphism on the second factor, and non-degenerate in the sense that any basis $\{\gamma_i\}$ for $\wt H_1^-(\Sigma)$ induces an isomorphism
\be \label{abiso} (x_{\gamma_i}) : \, \wt \CX_{\rm ab}^-[\Sigma] \,\overset\sim\to\, (\C^*)^{{\rm rank}\,\tilde H_1^-(\Sigma)} = (\C^*)^{{\rm rank}\, H_1^-(\Sigma)}\,.\ee

\noindent\textbf{Remark.} One could also consider a space $\wt \CX_{\rm ab}'[\Sigma]$ defined as the slice of $\wt \CX_{\rm ab}[\Sigma]$ on which $x_{P_+\gamma}=1$ for all $\gamma\in \wt H_1(\Sigma)$. Our space $\wt \CX_{\rm ab}^-[\Sigma]$ is a finite quotient of $\wt \CX_{\rm ab}'[\Sigma]$ by $(\Z_2)^s$ (acting as multiplication by $-1$ on some $\C^*$ factors), where $s$ is the number of $\Z_2$ factors in the torsion group $H_1(\Sigma)/(\im\,P_+\oplus\ker\,P_+)$. Under the non-abelianization map of Section \ref{sec:NAdef}, $\wt\CX_{\rm ab}^-[\Sigma]$ maps to framed flat $PGL(2)$ connections on a punctured base $\CC^*$, whereas $\wt\CX_{\rm ab}'[\Sigma]$ would map to twisted framed flat $SL(2)$ connections on $\C^*$. \medskip

The spaces $\wt \CX_{\rm ab}[\Sigma]$ and $\wt \CX_{\rm ab}^-[\Sigma]$ have natural holomorphic Poisson and symplectic structures given by the Atiyah-Bott formula $\omega_{\rm ab} = \int_\Sigma \delta \aa\wedge\delta \aa$.%
\footnote{In this formula, $\aa$ is a standard flat $GL(1)$ connection on $\Sigma$, represented locally as a 1-form. Implicitly, $\aa$ is obtained by using \emph{any} isomorphism between twisted and untwisted moduli spaces $\tilde \CX_{\rm ab}[\Sigma]\simeq  \CX_{\rm ab}[\Sigma]$. The symplectic form is independent of the choice of isomorphism. This is manifest in expressions \eqref{GL1PB}, \eqref{GL1w}, which are invariant under sign changes $x_\gamma\mapsto x_{\gamma+u}= -x_\gamma$.} %
In this abelian setting, it is trivial to compute the Poisson brackets (on $\wt\CX_{\rm ab}^-[\Sigma]$, say)
\be \label{GL1PB} \{\log x_\gamma,\log x_{\gamma'}\} = \langle \gamma,\gamma'\rangle\,,\qquad \gamma,\gamma'\in \wt H_1^-(\Sigma)\,. \ee
To write the symplectic form in coordinates, we must choose some $\{\gamma_i\}_{i=1}^{{\rm rank}\,H_1^-(\Sigma)}$, $\gamma_i\in \wt H_1^-(\Sigma)$, lifting a basis $\{\ol \gamma_i\}$ of $H_1^-(\Sigma)$. Letting $\epsilon_{ij} = \langle \gamma_i,\gamma_j\rangle$ denote the non-degenerate intersection pairing and $x_i:=x_{\gamma_i}$, we have
\be \label{GL1w} \omega_{\rm ab} = \frac12 \sum_{i,j}(\epsilon^{-1})^{ij}\frac{dx_i}{x_i}\wedge \frac{dx_j}{x_j}\,.\ee
More elegantly, the symplectic form is induced by the Poincar\'e dual of the cup product (composed with a projection to odd homology)
\be \label{GL1cup}
\cup^*:\begin{array}{ccc}
  H_2(\Sigma)&\to& H_1^-(\Sigma)\wedge H_1^-(\Sigma) \\
  {}[\Sigma] &\mapsto& \frac12\sum_{ij} (\epsilon^{-1})^{ij} \ol\gamma_i\wedge\ol\gamma_j\,.
\end{array}
\ee
From \eqref{GL1cup} we also obtain the canonical lift to K-theory (modulo torsion)
\be \label{GL1K2} \hat \omega_{\rm ab} := \frac12 \sum_{i,j}(\epsilon^{-1})^{ij}{x_i}\wedge {x_j}\;\in\, K_2(\wt \CX_{\rm ab}^-[\Sigma])_\Q \,.\ee

\subsection{Gluing abelian connections}
\label{sec:glueab}

Let $M$ be any triangulated framed 3-manifold, and $M'$ a framed 3-manifold glued by identifying faces of $M$. As usual we split the gluing into two parts $M\leadsto M_0\leadsto M'$. The gluing gluing maps for spaces of flat connections $\wt\CX_{\rm ab}^-$ on the canonical covers $\Sigma,\Sigma_0,\Sigma'$ of the respective boundaries. Due to the  pairing \eqref{defpairab}, these gluing maps are automatically dual to the gluing maps in homology.

Concretely, first consider the gluing $M\leadsto M_0$, in which only interiors of faces of a big-boundary triangulation $\mb t_{\rm 2d}$ are identified. The pre-images of each face in $T_1\Sigma$ retracts to a single $S^1$ fiber. Thus, any $GL(1)$ flat connection $\aa\in \wt\CX_{\rm ab}^-[\Sigma]$ can be trivialized along the gluing region (aside from the universal fiber holonomy $-1$), and the gluing gives a map
\be \label{gab1} \hspace{.3in} g_{GL(1)}^{(1)} : \begin{array}{cccc} \wt\CX_{\rm ab}^-[\Sigma] &\to\hspace{-.3cm}\to& \wt\CX_{\rm ab}^-[\Sigma_0] \\
 x_{\tilde g(\gamma)} & \mapsto & x_\gamma &\quad (\gamma\in \wt H_1^-(\Sigma_0))\,.
\end{array}
\ee
The holonomies $x_\gamma$ of $g_{GL(1)}^{(1)}(\aa)$ must be equal to the holonomies $x_{\tilde g(\gamma)}$ of $\aa$ itself, with $\wt g$ as in Section \ref{sec:glue-gen}. Since $\wt g$ is an injection, $g_{GL(1)}^{(1)}$ is a surjection. (More precisely, there is an isomorphism $g_{GL(1)}^{(1)}:\,\wt\CX_{\rm ab}^-[\Sigma]/{\rm coker}\,\wt g \overset\sim\to \wt\CX_{\rm ab}^-[\Sigma_0]$, with generators of the 2-torsion group ${\rm coker}\,\wt g\simeq{\rm coker}\,g$ acting as multiplication by $-1$ on $\C^*$ factors.)
 Moreover, since $\wt g$ preserves the intersection product, $g_{GL(1)}^{(1)}$ is a symplectomorphism. 

Second, when filling in defects $M_0\leadsto M'$, a flat connection on $\Sigma_0$ induces a flat connection on $\Sigma'$ if and only if its holonomies along defect cycles $\nu\in \wt G'$ are trivial (again using the notation of Section \ref{sec:glue-gen}). We denote this condition as ``$x_{\nu\in\tilde G'}= 1$.'' When it is satisfied, we may trivialize the flat connection in a neighborhood of the defects, then fill them in, obtaining a map
\be \label{gab2} \hspace{.3in} g_{GL(1)}^{(2)} : \begin{array}{cccc}
 \wt \CX_{\rm ab}^-[\Sigma_0]\big|\raisebox{-.15cm}{($x_{\nu\in\tilde G'}= 1$)} &
  \to\hspace{-.3cm}\to & \wt\CX_{\rm ab}^-[\Sigma'] \\
 x_\gamma &\mapsto & x_{\tilde q(\gamma)} & \quad (\gamma\in \wt K'\subset \wt H_1^-(\Sigma_0))\,.
\end{array}
\ee
Now since $\wt H_1^-(\Sigma',\Z)=\wt H_1^-(\Sigma_0)/\!/\wt G'=\wt K'/\wt G'$ (Prop. \ref{prop:glue}') we see that \eqref{gab2} is symplectic reduction,
$\wt\CX_{\rm ab}^-[\Sigma'] \simeq \wt\CX_{\rm ab}^-[\Sigma_0]/\!/(\C^*)^{{\rm rank}\,\tilde G'}= \wt\CX_{\rm ab}^-[\Sigma_0]\big|_{x_{\nu\in\tilde G'}=1}\big/(\C^*)^{{\rm rank}\,\tilde G'}$,
where the $(\C^*)^{{\rm rank}\,\tilde G'}$ action is generated by using the $x_{\nu\in \tilde G'}$ as moment maps (analogous to the action described in \eqref{mom} following Theorem \ref{thm:H}).

Combining the descriptions of \eqref{gab1}--\eqref{gab2}, we deduce that the combined gluing map $g_{GL(1)}:=g_{GL(1)}^{(2)}\circ g_{GL(1)}^{(1)}$ is a symplectic reduction of a finite quotient, governed by gluing equations
\be x_{\tilde g(\gamma)} = x_{\tilde q(\gamma)}\,,\qquad \forall\; \gamma\in \wt K'\,.\ee
When $M$ and $M'$ are defect-free, we may use Lemma \ref{lemma:SESP} to replace $\wt G',\wt K'$ by the finite-index subgroups $\wt G:=\wt h(\mb P_G)$ and $\wt K := \wt h(\mb P_0\oplus \Z_2)$ corresponding to the path algebra on $(\CC_0)_{\rm small}$. Then 
\be \label{ab-symp} g_{GL(1)}:\, \CX_{\rm ab}^-[\Sigma]\big|\raisebox{-.15cm}{($x_{\tilde g(\tilde G)}\equiv1$)} \to\hspace{-.3cm}\to \CX_{\rm ab}[\Sigma'] \simeq (\CX_{\rm ab}^-[\Sigma]/Z)/\!/(\C^*)^{{\rm rank}\, \tilde G}\,,
\ee
governed by gluing equations identical to \eqref{glue-simple}: $x_{\tilde g(\gamma)} = x_{\tilde q(\gamma)}$ for all $\gamma\in \wt K$. (The torsion group $Z\simeq \wt H_1^-(\Sigma)/\wt H$ acts on $\CX_{\rm ab}^-[\Sigma]$ the same way as in Theorem \ref{thm:H}.)

\subsection{Spectral networks and non-abelianization}
\label{sec:NAdef}

Suppose that $\CC$ is the boundary of a framed 3-manifold without defects, with canonical cover $\Sigma\overset\pi\to \CC$, and that $\pi_1(\CC_{\rm small})$ is abelian. Fix a triangulation $\mb t_{\rm 2d}$ of $\CC_{\rm big}$. By the results of Section \ref{sec:ab}, the symplectic moduli space of twisted abelian flat connections on $\Sigma$, $\wt \CX_{\rm ab}[\Sigma]\simeq (\C^*)^{{\rm rank}\,H_1^-(\Sigma)}$ is coordinatized by $x_\gamma$, $\gamma\in\wt H_1^-(\Sigma)$.
Define the open subset
\be \CP_{\rm ab}[\Sigma] := \big\{\aa\in \wt\CX_{\rm ab}^-[\Sigma]\,\big|\, x_\lambda^{(a)}(\aa)\neq 1,\;(x_\alpha^{(t)}(\aa),x_\beta^{(t)}(\aa))\neq (1,1) \big\} \subset  \wt\CX_{\rm ab}^-[\Sigma]\,, \label{def:CPab} \ee
on which holonomies around tori and annuli are never totally unipotent. This is analogous to the definition of the $PGL(2)$ moduli space $\CP[\CC,\mb t]\subset \CX[\CC]$ on p. \pageref{def:CP}. In this section we construct a non-abelianization map
\be \label{NA}
\Phi[\mb t_{\rm 2d}]:\; \CP_{\rm ab}[\Sigma] \to \CP[\CC;\mb t_{\rm 2d}]\,,
\ee
defined via a topological ``spectral network'' that (slightly) extends that of \cite{GMN-spectral} (\cf\ Footnote \ref{foot:HN}). We show that

\begin{prop} \label{prop:NA} With notation as above,
the map $\Phi[\mb t_{\rm 2d}]$ is 1-1 and a (holomorphic) $K_2$ symplectomorphism. It identifies twisted $GL(1)$ holonomies $x_\gamma$ with the  path-functions $x_\gamma$ on $\CP[\CC;\mb t_{\rm 2d}]$ labelled by elements $\gamma\in \wt H_1^-(\Sigma)$; in other words $\Phi[\mb t_{\rm 2d}]^*(x_\gamma) = x_\gamma$\,.

\end{prop}

We construct $\Phi$ in two steps. First, given a twisted abelian flat connection $\aa$ on $\Sigma$, we can push it forward to a framed flat $PGL(2)$ connection $\pi_*\aa$ on $\CC\bs \mathfrak b$, \ie\ in the complement of the branch locus,
\be \label{pfA} \pi_*: \wt \CX_{\rm ab}^-[\Sigma] \to \CX[\CC\bs\mathfrak b]\,.\ee
To see this, note that the projection $\pi:\Sigma\to \CC$ can be extended uniquely to a bundle map $\pi:T_1(\Sigma\bs\mathfrak b)\to T_1(\CC\bs\mathfrak b)$ that is globally two-to-one and an isomorphism on the unit-tangent fibers. Then a flat line bundle $L\to T_1\Sigma$ induces a flat rank-two bundle $E'\to T_1(\CC\bs\mathfrak b)$, locally of the form $L^+\oplus L^-$, where $L^\pm$ denote (locally) the bundles over the two sheets of the cover $\pi:T_1(\Sigma\bs\mathfrak b)\to T_1(\CC\bs\mathfrak b)$. The induced flat connection $\aa'$ on $E'$ (locally of the form $\aa^+\oplus \aa^-$) has holonomy valued in $GL(2)$, and equal to $-\id$ around unit-tangent fibers. Taking its projectivization, we obtain a connection $\aa''$ of $PGL(2)$ holonomy. Since $-\id\simeq \id$ in $PGL(2)$, the holonomy of $\aa''$ around unit-tangent fibers is trivial, so $\aa''$ descends to a flat $PGL(2)$ connection on a rank-two bundle $E_*\to(\CC\bs\mathfrak b)$, which we call $\pi_*\aa$.

A priori, this procedure describes a map $\pi_*$ from $\wt \CX_{\rm ab}[\Sigma]$ to $\CX[\CC\bs \mathfrak b]$. Our definition of the ``odd'' space $\wt \CX_{\rm ab}^-[\Sigma]$ in \eqref{defPab-}, however, was precisely engineered so that the map $\pi_*$ would factor through to $\wt \CX_{\rm ab}^-[\Sigma]$. Indeed, the action $x_\gamma\mapsto t^{\pi_*\circ\,p(\gamma)}x_\gamma$ in \eqref{defPab-} is just a lift of the projectivization action on $E'$.

The connection $\pi_*\aa$ is also naturally framed. Recall from Section \ref{sec:cc} that the canonical cover $\Sigma\overset\pi\to \CC$ can be glued together from two sheets $\Sigma^+$ and $\Sigma^-$, each a copy of $\CC\bs\Gamma_{\rm br}$. We correlate the local decomposition $E' = L^+\oplus L^-$ above with the labeling of these sheets.
Then induces a canonical diagonal decomposition $E_* = L^+\oplus L^-$ (as a projectively flat bundle) over each connected component of $\CC\bs\Gamma_{\rm br}$, and in particular over $\CC_{\rm small}$. We take the flag $L^+\subset E_*$ to be the framing data on each component of $\CC_{\rm small}$.

\begin{figure}[htb]
%\centering
\hspace{-.25in}
\includegraphics[width=6.5in]{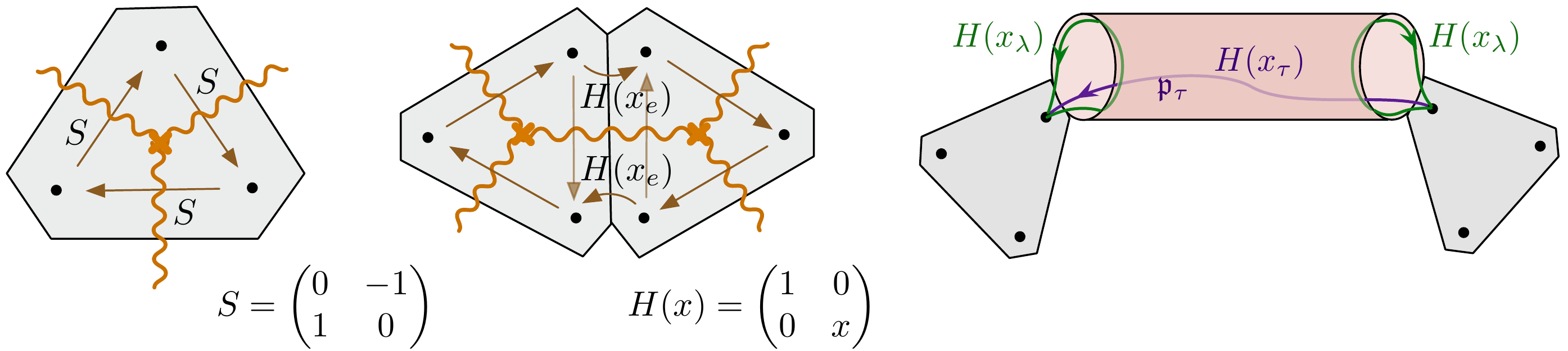}
\caption{Local frames and parallel transport for the flat connection $\pi_*\aa$}
\label{fig:ab-basis}
\end{figure}

The second step is to modify the push-forward connection $\pi_*\aa$, defined over $\CC\bs\mathfrak b$, to obtain a new (and more interesting) connection that extends over the branch points but may have unipotent singularities at punctures on the small boundary. We proceed as follows.

For small-sphere and small-torus components of $\CC$, no modification is required. The remaining connected components of $\CC$ consist of triangulated big boundary with holes filled in by small discs and annuli. We then choose local frames for the projectively flat bundle $E_*\to \CC\bs\mathfrak b$ over three points $p$ in each face $t$ of $\mb t_{\rm 2d}$ so that the parallel transport of $\pi_*\aa$ is given by the transformations in Figure \ref{fig:ab-basis}. Locally, each frame is a choice of vectors in the lines $L^\pm$. Thus the parallel transport is diagonal over connected components of $\CC\bs\Gamma_{\rm br}$, given by matrices $H(x_\pp)={\rm diag}\,(1,x_\pp)$ for appropriate functions $x_\pp=(x_e,x_\tau,x_\lambda)$ on $\wt X_{\rm ab}^-[\Sigma]$. Across branch cuts, the lines are exchanged and the parallel transport $S=\left(\begin{smallmatrix} 0 & -1 \\ 1 & 0 \end{smallmatrix}\right)$ is anti-diagonal.%
\footnote{A similar choice of local frames was used in the non-abelianization constructions of \cite{GMN-spectral, GMN-snakes}.
Note that the sign in $S=\left(\begin{smallmatrix} 0 & -1 \\ 1 & 0 \end{smallmatrix}\right)$ is particularly natural given the push-forward construction of $\pi_*\aa$. Namely, since $\aa$ has holonomy $-1$ around unit-tangent fibers of $T_1\Sigma$, it follows that the holonomy $h_{\rm br}$ of $\aa'$ around any cycle surrounding a branch point (and wrapping the unit-tangent fiber of $T_1\CC$ any number of times) must satisfy $h_{\rm br}^2 = -\id \in GL(2)$. To be compatible with the local decomposition $E'=L^+\oplus L^-$, it must also be anti-diagonal, and can be chosen (modulo gauge equivalence) as $h_{\rm br}=\pm\left(\begin{smallmatrix} 0 & -1 \\ 1 & 0 \end{smallmatrix}\right)$. After projectivizing, we obtain holonomy $S^3 = S = \left(\begin{smallmatrix} 0 & -1 \\ 1 & 0 \end{smallmatrix}\right)\in PGL(2)$ around branch points.}
The nontrivial holonomy $S^3=S$ around branch points is what prevents $E_*$ from extending over the branch points.

\begin{figure}[htb]
\hspace{-.25in}
\includegraphics[width=6.5in]{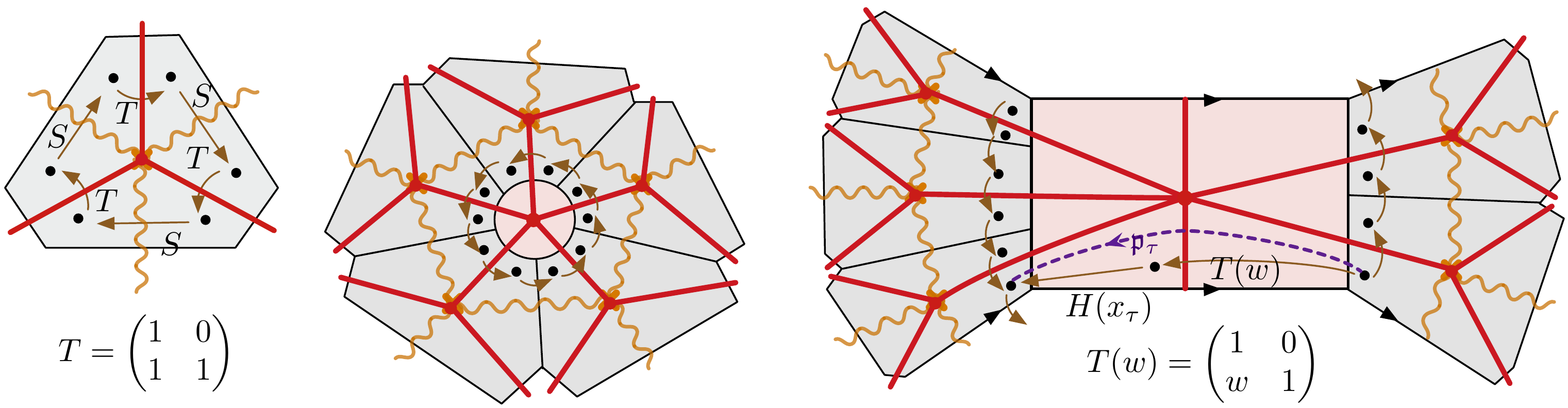}
\caption{Structure of the spectral network $\CW$ (in red) on a face of $\mb t_{\rm 2d}$ (left), in the neighborhood of a small disc (center), and in the neighborhood of a small annulus (right). We indicate the unipotent modifications of $\pi_*\aa$ corresponding to the edges of the spectral network on faces (left) and wrapping around annuli (right).}
\label{fig:W}
\end{figure}

To fix this we introduce a (topological) \emph{spectral network} $\CW=\CW[\mb t_{\rm 2d}, \{\pp_\tau\}]$ on $\CC$. It depends on $\mb t_{\rm 2d}$ and has a mild dependence on a choice of paths $\pp_\tau$ traversing small annuli.
Roughly, $\CW$ is an unoriented graph dual to small edges of $\mb t_{\rm 2d}$ and to the paths $\pp_\tau$, as in Figure~\ref{fig:W}. Formally, $\CW$ has vertices at all branch points and at one point in the interior of each small disc and annulus. Its edges connect the vertex in each small disc $d$ to every branch point in a face of $\mb t_{\rm 2d}$ surrounding $d$; connect the vertex in each small annulus $a$ to the branch points in the faces surrounding $a$, without crossing $\pp_\tau^{(a)}$; and connect the vertex in $a$ to itself via a circular path in the homology class of~$\pp_\lambda^{(a)}$.

We perform unipotent modifications of $E_*$ across walls of the spectral network, in two rounds. First, note that each edge $s$ of $\CW$ ending at a branch point $b\in\mathfrak b$  passes close to a single point $p$ labeling a frame for $E_*$. We add a unipotent modification by $T=\left(\begin{smallmatrix} 1 & 0 \\ 1 & 1 \end{smallmatrix}\right)$ when moving across $s$ clockwise from the perspective of $b$ --- meaning explicitly that we split the one frame over $p$ into two, with new parallel transport $T$ between them. This creates a new projectively flat bundle $E_*'$. Its holonomy around any branch point is trivial by virtue of the relation $(ST)^3 = \id$, hence $E_*'$ extends over the branch locus. However, $E_*'$ may have nontrivial unipotent holonomy at the vertices of $\CW$ on small discs and annuli.

Now, for each small annulus $a$, we add a second unipotent modification to $E_*'$ by $T(w_a) := \left(\begin{smallmatrix} 1 & 0 \\ w_a & 1 \end{smallmatrix}\right)$ across the edge of $\CW$ wrapped around $a$ (homologous to $\pp_\lambda^{(a)}$). For concreteness, we take this modification to split either of the two frames of $E_*'$ at the tail of $\pp_\tau^{(a)}$. (These frames were already split by the first modification; since $[H(w_a),H(1)]=0$, it does not matter which we take.) Then, as long as $x_\lambda^{(a)}\neq 1$, there is a unique $w_a\in\C$ that trivializes the holonomy around the vertex of $\CW$ on $a$. (The calculation is identical to solving for $t$ in \eqref{det-t}, p. \pageref{det-t}.)

Following these two rounds of unipotent modifications, we arrive at a projectively flat bundle $E$ that extends over all of $\CC$ except the vertices of $\CW$ on small discs --- \ie\ over the ``punctured boundary'' $\CC^*$ of Section \ref{sec:PGLdef}. Let $\Phi(\aa)$ denote the flat $PGL(2)$ connection on $E$. Since unipotent modifications preserve the flag $L^+\subset E$ over each component of $C_{\rm small}$, $\Phi(\aa)$ is a naturally a framed flat connection. Thus, we've defined a map
\be \label{defNA}
\Phi:\, \wt \CX_{\rm ab}^-[\Sigma]\big|_{x_\lambda\neq 1} \to \CX[\CC]\,.
\ee
Its definition depends both on a triangulation $\mb t_{\rm 2d}$ and on a choice of paths $\pp_\tau$. We will prove momentarily that the dependence on $\pp_\tau$'s is trivial and that when $\Phi$ is restricted to both $x_\lambda\neq 1$ (for annuli) and $(x_\alpha,x_\beta)\neq (1,1)$ (for tori) its image is precisely $\CP[\CC,\mb t_{\rm 2d}]$. \medskip

\noindent \textit{Proof of Prop. \ref{prop:NA}}\;   Given $\aa\in \wt\CX_{\rm ab}^-[\Sigma]$ with $x_\lambda\neq 1$ and $(x_\alpha,x_\beta)\neq 1$, observe that $\Phi(\aa)$ is a framed flat $PGL(2)$ connection on $\CC^*$ with exactly the same set of local frames and parallel transports as are used in Appendix \ref{app:traffic} to uniquely reconstruct a connection $\CA\in \CP[\CC,\mb t_{\rm 2d}]$ from its coordinates. Thus $\Phi(\aa)\in \CP[\CC,\mb t_{\rm 2d}]$. Indeed, the path-coordinates $x_\pp$ of $\Phi(\aa)\in \CP[\CC]$, showing up as diagonal entries of parallel transport matrices, are manifestly equal to the corresponding path-coordiantes $x_\pp$ labeling the abelian holonomies of $\aa$. (In each case, $x_\pp$ depends only on the twisted homology class $\wt h(\pp)\in \wt H_1^-(\Sigma)$.) From \eqref{abiso} and Prop. \ref{prop:coords} it follows that $\Phi:\CP_{\rm ab}[\Sigma]\to \CP[\CC,\mb t_{\rm 2d}]$ is one-to-one.

Since path-functions $x_\pp$ uniquely determine $\Phi(\aa)$, independent of any particular choice of traversing paths $\pp_\tau$ for the annuli, the map $\Phi$ cannot depend on the choice of $\pp_\tau$'s used in unipotent modification. Thus it depends at most on a big-boundary triangulation $\mb t_{\rm 2d}$.

The fact that $\Phi$ is a $K_2$ symplectomorphism follows immediately by comparing the $K_2$ forms \eqref{K2}, \eqref{GL1K2} (\cf\ the Poisson brackets \eqref{PB}, \eqref{GL1PB}), which look identical in $x_\gamma$ coordiantes, and are both controlled by the intersection form on $H_1^-(\Sigma)$.
There is also an alternative, coordinate-free, proof of the fact that $\Phi$ is a symplectomorphism. Namely, we decompose $\Phi$ as a composition of push-forward ($\pi_*$) and unipotent modifications. The push-forward map is obviously a symplectomorphism for the Atiyah-Bott symplectic/Poisson structures. It was then shown in \cite[Section 10.4]{GMN-spectral} that unipotent modification preserves the holomorphic symplectic structure. \; $\square$

\subsection{Non-abelianization commutes with gluing}
\label{sec:NAglue}

We combine the results of the previous sections in a final theorem about gluing. Suppose that we glue framed 3-manifolds $M\leadsto M_0\leadsto M'$ as in Section \ref{sec:glue-gen}, where $M$ and $M'$ have no defects, and their small boundaries $\CC_{\rm small}$ and $\CC'_{\rm small}$ have only discs, annuli, and tori. Let us fix compatible big-boundary triangulations $\mb t_{\rm 2d}$ for $\CC_{\rm big}$ and $\mb t_{\rm 2d}'$ for $(\CC_0)_{\rm big}=\CC_{\rm big}'$.

\begin{thm} \label{thm:NA}

Gluing and non-abelianization maps fit into a commutative diagram
\be \begin{array}{cccl} \label{NAglue}
  \CP_{\rm ab}[\Sigma]\big|_{x_{\tilde g(\tilde G)}\equiv 1,\,R'} &\overset{g_{GL(1)}}{\to\hspace{-.3cm}\to} &  \CP_{\rm ab}[\Sigma'] &=  (\CP_{\rm ab}[\Sigma]_{R'}/Z) \big/\!\!\big/ (\C^*)^{{\rm rank}\,\tilde G} \\[.2cm]
 \hspace{-.38in}\Phi[\mb t_{\rm 2d}]\downarrow & & 
  \hspace{-.38in}\Phi[\mb t_{\rm 2d}']\downarrow \\[.2cm]
 \CP[\CC;\mb t_{\rm 2d}]\big|_{x_{\tilde g(\tilde G)}\equiv 1,\,R'} & \overset{g_{PGL(2)}}{\to\hspace{-.3cm}\to} &  \CP[\CC';\mb t_{\rm 2d}'] &
  = (\CP[\CC;\mb t_{\rm 2d}]_{R'}/Z)\big/\!\!\big/ (\C^*)^{{\rm rank}\,\tilde G}\,,
\end{array}\ee
where $R'$ is the technical restriction $x_{\tilde g\tilde q^{-1}(\lambda)}\neq 1$, $(x_{\tilde g\tilde q^{-1}(\alpha)},x_{\tilde g\tilde q^{-1}(\beta)})\neq (1,1)$ and $Z\simeq \wt H_1^-(\Sigma)/\tilde H$ (as in Theorem \ref{thm:H}). The vertical maps are 1-1 $K_2$ symplectomorphisms and the horizontal maps are $K_2$ symplectic reduction of finite quotients.
Thus the trivial symplectic reduction in the gluing of $GL(1)$ moduli spaces on $\Sigma$ \eqref{ab-symp} induces the non-trivial symplectic reduction (Theorem \ref{thm:H}) in the gluing of $PGL(2)$ moduli spaces on $\CC$. Both reductions are governed by the same gluing equations
\be  x_{\tilde g(\gamma)} \,=\, x_{\tilde q(\gamma)}\,,\qquad \forall\;\gamma\in \tilde K \subset \tilde H_1^-(\Sigma_0)\,, \label{NAeqs} \ee
ultimately arising from the isomorphism $\wt H_1^-(\Sigma') = \wt H/\!/\wt g(\wt G) = \wt g(\wt K)/\wt g(\wt G)$ with $\wt H\subset \wt H_1^-(\Sigma)$ of finite index (Lemma \ref{lemma:SESP}, Prop. \ref{prop:glue}').
\end{thm}

\noindent\emph{Proof.} Commutativity of the square \eqref{NAglue} follows easily by comparing the description of abelian gluing in Section \ref{sec:glueab} to the description of $PGL(2)$ gluing in Section \ref{sec:gluePGL}. Alternatively (and more explicitly), we have already seen in \eqref{glue-simple} and \eqref{ab-symp} that both the $GL(1)$ and $PGL(2)$ gluing maps are governed by the same gluing equations \eqref{NAeqs}, labelled by elements of twisted odd homology; since the non-abelianization maps $\Phi$ preserve the $x_\gamma$ functions (by Prop. \ref{prop:NA}), it follows that the square must commute. The remaining claims follow immediately from the result of Prop. \ref{prop:NA} that the non-abelianization maps are 1-1 symplectomorphisms. \; $\square$

\subsection*{Acknowledgements}

We especially wish to thank D. Futer, L. Hollands, A.B. Goncharov, G. Moore, A. Neitzke, and C. Zickert for helpful conversations and advice. We are also indebted to Center for Quantum Geometry of Moduli Spaces for their hospitality during the workshop \emph{Pressure Metrics and Higgs Bundles}, Aug 19-22, 2013 (partial funding through the GEAR network), where the ideas for this project were conceived. T.D. 
is supported by DOE grant DE-SC0009988, and in part by the European Research
Council under the European Union's Seventh Framework Programme (FP7/2007-2013), ERC Grant agreement 335739.
R.V. is supported by a VENI grant from the Netherlands Organization for Scientific Research.

\appendix

\section{Odd results}
\label{app:odd}

In this section, we collect some basic results about odd homology, and review their proofs. Notation is as in Section \ref{sec:oddH}. In particular, $\Sigma\overset\pi\to\CC$ is an oriented double cover of a closed, oriented surface, branched over a locus of points $\mathfrak b$.

\subsection{Basics}
\label{app:basics}

Let $\sigma :\Sigma\to \Sigma$ denote the deck-transformation homeomorphism. It preserves orientation and its fixed-point locus is precisely $\mathfrak b(\pi)$. It induces a push-forward automorphism on homology groups $\sigma_*: H_\bullet(\Sigma)\overset{\sim}{\to} H_\bullet(\Sigma)$.  Letting $P_\pm := \id \pm \sigma_*$, we define
\be \label{defH-} H_\bullet^-(\Sigma) := \ker\,P_+\,,\qquad H_\bullet^+(\Sigma) := \ker\,P_-\,. \ee
Notice that
\be \label{Pprops} P_\pm^2 = 2P_\pm\,,\qquad P_\pm P_\mp = 0\,,\qquad P_++P_-=2\,\id\,.\ee
Thus, $P_\pm$ are close to being orthogonal projection operators. They fail to be proper projections due to the factors of 2 in \eqref{Pprops}.

\begin{lemma} \label{lemma:P}
The quasi-projections $P_\pm$ obey the following properties:

a) $\im\,P_\pm \subset \ker\, P_\mp$;

b) $\ker\,P_+\cap \ker\,P_- = \{0\}$;

c) $\ker\, P_\mp/\im\,P_\pm \simeq (\Z_2)^{r_\pm}$ (for some $r_\pm$) is finite, containing only 2-torsion;

d) $H_\bullet/(\ker\, P_+\oplus \ker\, P_-)$ is also finite, containing only 2-torsion.

\end{lemma}

\noindent\emph{Proof.} Parts (a-b) follows trivially from \eqref{Pprops}. Parts (c-d) follow because $P_\pm$ can be promoted to honest orthogonal projections after tensoring with $\Q$. More concretely, given any $\alpha \in H_\bullet(\Sigma)$ such that $\alpha = 2\beta$ for some $\beta\in H_\bullet(\Sigma)$ we can uniquely decompose $\alpha = \alpha_++\alpha_-$ with $\alpha_\pm \in \im\,P_\pm$ (namely, $\alpha_\pm = P_\pm\beta$), proving (d); and given any $\beta\in \ker\, P_\mp$ we have $2\beta=P_\pm \beta$, proving (c).  \;$\square$ \medskip

The first homology group $H_1(\Sigma)$ has a non-degenerate skew-symmetric intersection form $\langle\;,\;\rangle: \bigwedge^2 H_1(\Sigma,\Z)\to \Z$, which is preserved by $\sigma_*$ (since $\sigma$ is an orientation-preserving homeomorphsim). Thus it is also preserved by $P_\pm$. We find that

\begin{lemma} \label{lemma:int}
Letting $\ker\,P_\pm,\im\,P_\pm$ denote kernels and images in $H_1(\Sigma)$ now, we have

a) $\langle \ker P_\pm,\ker P_\mp \rangle = 0$\,;

b) the intersection form is non-degenerate on $\ker\,P_+$ and on $\ker\, P_-$\,;

c) the intersection form is even on $\im\,P_+$ and $\im\,P_-$\,;

d) if $\langle \alpha,\beta\rangle=0$ for all $\beta\in \im\, P_+$ then $\alpha \in \ker\,P_+$ (and similarly with $+\leftrightarrow -$).

\end{lemma}

\noindent \emph{Proof.} These are all simple consequences of invariance of the intersection form under $\sigma_*,P_\pm$. For example, if $\alpha_\pm\in\ker\,P_\pm$ then $\langle \alpha_+,\alpha_-\rangle = \langle\sigma\alpha_+,\sigma\alpha_-\rangle = -\langle \alpha_+,\alpha_-\rangle$, hence $\langle \alpha_+,\alpha_-\rangle=0$, proving (a). Nondegeneracy on $H_1(\Sigma)$ then implies (b) and (d). For (c) note that if $\alpha,\beta\in \im\,P_+$ (say) then $\langle \alpha,\beta\rangle = \langle P_+\gamma,\beta\rangle = \langle \gamma,P_+\beta\rangle=2\langle\gamma,\beta\rangle\in 2\,\Z$\,. \;$\square$ \medskip

\begin{figure}[htb]
\centering
\includegraphics[width=4in]{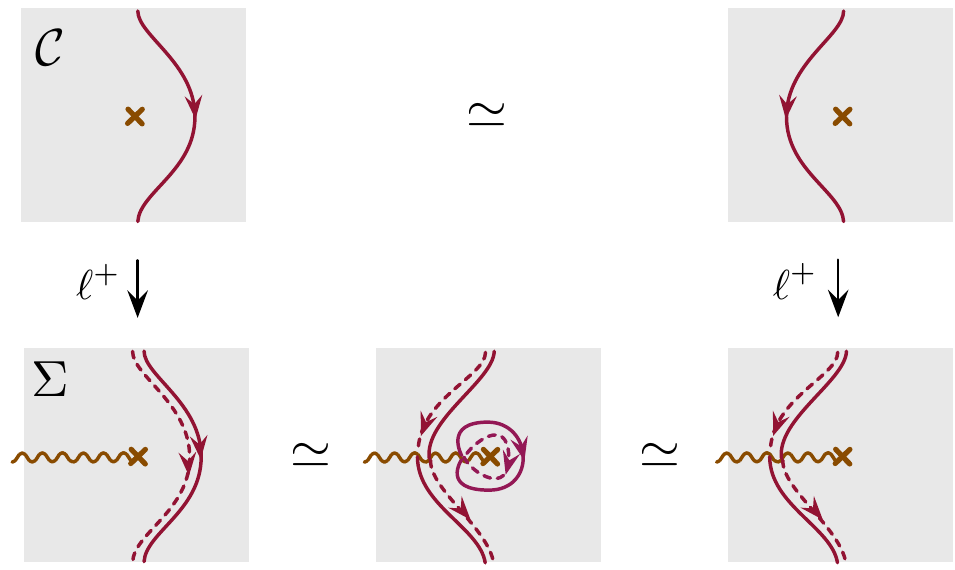}
\caption{Demonstrating the invariance of the even lift $\ell^+(\gamma) = [\pi^{-1}(\gamma)]$ under a homotopy of $\gamma$ through a branch point on $\CC$. We represent $\Sigma$ locally by drawing a branch cut emanating from the branch point, and draw the curves $\pi^{-1}(\gamma) \subset \Sigma$ using solid and dashed lines for the two sheets of the cover. On the bottom row of the figure, ``$\simeq$'' denotes equivalences in $H_1(\Sigma)$.}
\label{fig:ell}
\end{figure}

Let $\pi_*:H_1(\Sigma)\to H_1(\CC)$ denote the induced action of the projection $\Sigma\overset\pi\to\CC$\,. For any oriented curve $\gamma\subset (\CC\bs\mathfrak b)$ (in the complement of the branching locus), the preimage $\pi^{-1}(\gamma)$ consists of one or two oriented curves on $\Sigma$. Let $\ell^+(\gamma)$ denote the homology class $[\pi^{-1}(\gamma)]\in H_1(\Sigma)$, and notice that it only depends on the homology class $[\gamma]\in H_1(\CC)$. (It is clear that $\ell^+(\gamma)$ is invariant under local homotopies of $\gamma$ that do not cross the branching locus $\mathfrak b$. If a homotopy crosses $\mathfrak b$, invariance of $\ell^+(\gamma)$ is illustrated in Figure \ref{fig:ell}.)
Thus the ``even lift'' $\ell^+$ can be extended by linearity to a map $\ell^+:H_1(\CC)\to H_1(\Sigma)$. Moreover, we have
\be \pi_*\circ \ell^+ = 2\id\,, \qquad \ell^+\circ \pi_*=P_+\,,\qquad P_-\circ \ell^+ = 0\,.\ee

\begin{lemma} \label{lemma:+-}

a) $H_1^-(\Sigma) = \ker\, \pi_*|_{H_1(\Sigma)}$\,;

b) $\ell^+:H_1(\CC) \hookrightarrow H_1^+(\Sigma)$ is an injection with finite (2-torsion) cokernel\,;

c) thus $H_1(\Sigma,\Q) = H_1^+(\Sigma,\Q)\oplus H_1^-(\Sigma,\Q) \simeq H_1(\CC,\Q) \oplus H_1^-(\Sigma,\Q)$\,, and in particular
\be {\rm rank}\,H_1^-(\Sigma) = {\rm rank}\,H_1(\Sigma) - {\rm rank}\,H_1(\CC)\,. \notag \ee
\end{lemma}

\emph{Proof.} From $\pi_*\circ \ell^+ = 2\id$ we see that $\ell^+$ is injective, so from $\ell^+\circ \pi_*=P_+$ it follows that $\ker\,\pi_*|_{H_1(\Sigma)} = \ker\,P_+|_{H_1(\Sigma)}$, proving (a). From $P_-\circ\ell^+=0$ (or simply from the definition of $\ell^+$) we see that $\im\,\ell_+ \subset H_1^+(\Sigma)$. Moreover, given any $\alpha\in H_1^+(\Sigma)$ we have $2\alpha = P_+\alpha = \ell^+(\pi_*\alpha)$, hence $2\alpha\in \im\,\ell_+$, proving (b). Part (c) follows from (b) and Lemma \ref{lemma:P}d. \;$\square$

\subsection{Chain complexes}
\label{app:cx}

To analyze the effect of gluing on odd homology, it is convenient to have odd versions of standard exact sequences. The following result shows that we can restrict exact sequences to odd homology, modulo 2-torsion.

\begin{lemma} \label{lemma:cx}
Let \;$(A_\bullet,\delta_\bullet)=\big[\to A_i \overset{\delta_i}\to A_{i-1} \overset{\delta_{i-1}}\to A_{i-2}\to\big]$\; be a chain complex of abelian groups. Let $\sigma$ be an involution of $A_\bullet$ that preserves grading and commutes with $\delta_\bullet$. Set $P_+ := \id + \sigma$ and $A_\bullet^-:= \ker\,P_+|_{A_\bullet}$, and let $\delta_\bullet^-$ be the restriction of $\delta_\bullet$ to $A_\bullet^-$. Then

a) $(A_\bullet^-,\delta_\bullet^-)$ is also a chain complex;

b) $\sigma$ induces an involution on homology $H_\bullet(A_\bullet,\delta_\bullet)$, and letting $H_\bullet^-(A_\bullet,\delta_\bullet):= \ker\, P_+|_{H_\bullet(A_\bullet,\delta_\bullet)}$, there is an isomorphism $H_\bullet(A_\bullet^-,\delta_\bullet^-) \simeq H_\bullet^-(A_\bullet,\delta_\bullet)$ modulo 2-torsion;

c) If $(A_\bullet,\delta_\bullet)$ is exact and injective on the left, say \;$0\to A_d \overset{\delta_d}\hookrightarrow A_{d-1} \overset{\delta_{d-1}}\to A_{d-2}\to...$\; for some $d$, then $(A_\bullet^-,\delta_\bullet^-)$ is exact in the first two places (at $A_d^-$ and $A_{d-1}^-$) and its homology is 2-torsion thereafter.

\end{lemma}

\emph{Proof.} For (a), observe that if $\alpha\in A_i^-$ then $P_+(\delta_i\alpha) = \delta_i(P_+\alpha) = 0$. Therefore, $\im\,\delta_i^-\subset \im\,\delta_i\cap A_{i-1}^- \subset \ker\,\delta_{i-1}^-$, showing that $(A^-_\bullet,\delta^-_\bullet)$ is a complex.

Next, observing that $\sigma$ fixes both $\ker\,\delta_\bullet$ and $\im\, \delta_\bullet$, since it commutes with $\delta_\bullet$. Thus there is an induced involution on classes in the $i$-th homology group $H_i = \ker\,\delta_i/\im\,\delta_{i+1}$, given by $\sigma[\alpha] :=[\sigma\alpha]$, and it makes sense to consider
\be H_i^-:= \ker\,P_+|_{H_i} = \langle\, \alpha \in \ker\,\delta_i\,|\, P_+\alpha \in \im\,\delta_{i+1}\,\rangle\big/\im\,\delta_{i+1}\,.\ee
Any class $[\alpha]\in H_i^-$ can be represented by $\alpha\in \ker\,\delta_i$ with $P_+\alpha = \delta_{i+1}\beta$. We apply $P_+$ to this relation to find $2P_+\alpha = P_+(\delta_{i+1}\beta)$. Letting $\alpha'=2\alpha-\delta_{i+1}\beta$ we see that $2[\alpha]\simeq [\alpha']$ and $P_+\alpha'=0$. Thus two times every class in $H_i^-$ has a representative annihilated by $P_+$. This means that modulo 2-torsion
\be H_i^-\simeq  (\ker\,\delta_i\cap A_i^-) \big/(\im\,\delta_{i+1}\cap A_i^-)  = \ker\,\delta_i^- \big/(\im\,\delta_{i+1}\cap A_i^-)\,. \ee
(More precisely, the RHS injects into the LHS with finite cokernel.)
Moreover, if $\alpha\in \im\,\delta_{i+1}\cap A_i^-$, \ie\ $\alpha = \delta_{i+1}\beta$ and $P_+\alpha = 0$, then $2\alpha = P_-\alpha = \delta_{i+1}(P_-\beta) \in \im\,\delta_{i+1}^-$. Therefore, the quotient $(\im\,\delta_{i+1}\cap A_i^-)/\im\,\delta_{i+1}^-$ is 2-torsion, and
\be  H_i^-\simeq \ker\, \delta_i^-\big/\im\,\delta_{i+1}^- = H_i(A_\bullet^-,\delta_\bullet^-) \ee
modulo 2-torsion, as needed for part (b).

Finally, suppose $(A_\bullet,\delta_\bullet)$ is an exact sequence as in (c). Part (b) implies that $H_\bullet(A_\bullet^-,\delta_\bullet^-)$ is 2-torsion; but we can do better. Since $\delta_d$ is injective, $\delta_d^-$ is also injective, so $H_d(A_\bullet^-,\delta_\bullet^-)=0$. Moreover, $\ker\,\delta_{d-1}^- = \ker\,\delta_{d-1} \cap A_i^- = \im\,\delta_{d}\cap A_i^- = \im\,\delta_d^-$ (the last equality again follows because $\delta_d$ is injective: if $\alpha=\delta_d\beta$ and $P_+\alpha=0$ then $\delta_d(P_+\beta)=0$, hence $P_+\beta=0$), so $H_{d-1}(A_\bullet^-,\delta_\bullet^-)=0$. \;$\square$

\subsection{Cellular description}
\label{app:cell}

For a branched cover $\Sigma\overset\pi\to\CC$, we can describe a set of generators for $H_1^-(\Sigma)$ and the intersection form on them very explicitly. We do so by applying Lemma \ref{lemma:cx} to a cell complex for~$\Sigma$.

Choose a finite cell decomposition of $\CC$ with 2-cells $\{f_i\}$, 1-cells $\{e_i\}$ and 0-cells $\{p_i\}$, such that every branch point of the cover $\Sigma\overset\pi\to\CC$ is a 0-cell. By lifting to the two sheets of $\Sigma$, this induces a cell decomposition of $\Sigma$ with 2-cells $D_2=\{f_i^+,f_i^-\}$, 1-cells $D_1=\{e_i^+,e_i^-\}$, and 0-cells $D_0=\{p_i^+,p_i^-|\,p_i\notin \mathfrak b\}\cup \mathfrak b$ (where superscripts $\pm$ indicate local choices of lifts). We immediately recover from this the standard Riemann-Hurwitz formula
\be \chi(\Sigma) = 2\chi(\CC) - \#(\mathfrak b)\,, \label{RH} \ee
and in combination with Lemma \ref{lemma:+-}(c) we get
\begin{lemma} \label{lemma:ranks}
Suppose that $\CC$ and $\Sigma$ are both connected. Then their genera are related by $g(\Sigma)=2g(\CC)+\#(\mathfrak b)/2-1$ and
\be \label{ranks}
\begin{array}{rcl}
{\rm rank}\,H_1(\Sigma) &=& 2\,{\rm rank}\,H_1(\CC)+\#(\mathfrak b)-2\,,\\
  {\rm rank}\,H_1^-(\Sigma) &=& {\rm rank}\,H_1(\CC)+\#(\mathfrak b)-2\, \\
   &=&  -\chi(\CC) + \#(\mathfrak b)\,.
 \end{array}
\ee
\end{lemma}

However, we can do better than describing ranks. Let
\be (C_\bullet,\pd_\bullet) = \big[0\to C_2\overset{\pd_2}\to C_1\overset{\pd_1}\to C_0\to 0\big] \ee
denote the cellular chain complex corresponding to $D_\bullet$ (so that $C_i = \Z\langle D_i\rangle$), whose homology is $H_\bullet(C_\bullet,\pd_\bullet) = H_\bullet(\Sigma)$. Let
\be (C_\bullet^-,\pd_\bullet^-) = \big[0\to C_2^-\overset{\pd_2^-}\to C_1^-\overset{\pd_1^-}\to C_0^-\to 0\big] \ee
be the chain complex obtained by applying $\ker\, P_+$ to the groups $C_i$ and specializing the boundary maps, as in Lemma \ref{lemma:cx}.
(This makes sense because the deck transformation $\sigma$ commutes with $\pd_\bullet$.)
Notice that $C_2^-$ and $C_1^-$ are generated by ``odd lifts'' of 2-cells and 1-cells
\be C_2^- = \Z\langle\ell^-(f_i)\rangle\,, \qquad C_1^- = \Z\langle\ell^-(e_i)\rangle\,,\ee
where $\ell^-(f_i) := f_i^+-f_i^-$ and $\ell^-(e_i):=e_i^+-e_i^-$.
It follows from Lemma \ref{lemma:cx} that $H_1^-(\Sigma) \simeq \ker\,\pd_1^-/\im\,\pd_2^-$ (modulo 2-torsion). More precisely, we have

\begin{lemma} \label{lemma:cell}
Suppose $\Sigma$ and $\CC$ are both connected. With a cell decomposition and notation as above:

a) $H_1^-(\Sigma) = \ker\, \pd_1^-\big/(\im\, \pd_2\cap C_1^-)$ is generated by 1-cycles $\sum_i a_i \ell^-(e_i)$ formed from odd lifts.

b) If the 0-cells consist entirely of branch points ($D_0=\mathfrak b$) then every odd lift $\ell^-(e_i)$ is automatically closed, so $H_1^-(\Sigma)$ is generated by the $\ell^-(e_i)$.

c) In (b), the intersection product $\langle \ell^-(e_i),\ell^-(e_j)\rangle$ equals the number of common endpoints of $e_i$ and $e_j$, counted with orientation (which is determined by the choice of lifts made in defining $\ell^-$).

\end{lemma}

\emph{Proof.} First observe that $\pd_2^-$ is injective. Indeed, the kernel of $\pd_2$ is generated by the fundamental class $[\Sigma] = \sum_i (f_i^++f_i^-)$, which is even; so if $\pd_2\beta = 0$, then $\beta = c[\Sigma]$, and $P_-\beta = c\,P_-[\Sigma]=0$; and if $\beta\in C_2^-$ as well then $P_-\beta=P_+\beta=0$ imply (Lemma \ref{lemma:P}(b)) $\beta=0$. It is also useful to note that $C_2^\pm=\ker\,P_\mp|_{C_2} = \im\,P_\pm|_{C_2}$, and similarly for $C_1$.

Now, as in the proof of Lemma \ref{lemma:cx}, we generally have that $H_1^-(\Sigma) = \ker\,P_+\big|(\ker\,\pd_1/\im\,\pd_2)$ is generated by classes $[\alpha]$ for $\alpha\in C_1$ such that $\pd_1\alpha=0$ and $P_+\alpha = \pd_2\beta$ for some $\beta\in C_2$. Then $0=P_-P_+\alpha = P_-\pd_2\beta = \pd_2(P_-\beta)$. By injectivity of $\pd_2^-$ this implies $P_-\beta=0$, and since $\ker\,P_-|_{C_2}=\im\,P_+|_{C_2}$ we get $\beta = P_+\gamma$ for some $\gamma \in C_2$. Setting $\alpha' := \alpha-\pd_2\gamma$, we find $[\alpha]=[\alpha']$ and $P_+\alpha'=0$. Therefore, $H_1^-(\Sigma)$ is equally well generated by classes $[\alpha]$ for $\alpha\in \ker\,\pd_1$ such that $P_+\alpha=0$, and we obtain the identity in part (a).

For part (b), notice that if $D_0=\mathfrak b$ consists entirely of branch points then $\pd_1(e_i^+) = \pd_1(e_i^-)$ for all 1-cells $e_i$. Therefore $\pd_1\,\ell^-(e_i) = 0$. Part (c) follows by noting that the only intersections of $\ell^-(e_i)$ and $\ell^-(e_j)$ can occur at branch points, and are simple. \;$\square$ \medskip

There is one situation not strictly covered by Lemma \ref{lemma:cell}.
If $\CC$ is connected but $\Sigma$ is not, then $\Sigma \simeq \CC^+\sqcup \CC^-$ is a trivial, disconnected double cover (here the sheets $\CC^\pm$ are identical copies of $\CC$). There can be no branch points. The homology $H_1(\Sigma)=H_1(\CC^+)\oplus H_1(\CC^-)$ is generated by lifts $\gamma^\pm$ of cycles $\gamma\in H_1(\CC)$ to the two sheets. Let $\ell^-(\gamma):=\gamma^+- \gamma^-$ denote the odd lift. Than it is easy to see from the decomposition $H_1(\Sigma)=H_1(\CC^+)\oplus H_1(\CC^-)$ that

\begin{lemma} \label{lemma:dis}

If $\Sigma\simeq \CC^+\sqcup \CC^-$ is a disconnected double cover, then

a) $\ell^-:H_1(\CC) \to H_1^-(\Sigma)$ is an isomorphism\,; and

b) $\langle \ell^-(\gamma),\ell^-(\gamma') \rangle = 2\langle \gamma,\gamma'\rangle\,.$
\end{lemma}

One of the main lessons of Lemmas \ref{lemma:cell}--\ref{lemma:dis} is that for any cover $\Sigma\overset\pi\to \CC$ the homology $H_1^-(\Sigma)$ can be represented by (sums of) curves $\gamma\subset \Sigma$ that are fixed set-wise by $\sigma$, such that $\sigma(\gamma)=\ol \gamma$ merely flips orientation. These curves are odd lifts $\ell^-(*)$ of either closed curves on $\CC$ or edges that connect branch points in a cell decomposition of $\CC$. Such generators of $H_1^-(\Sigma)$ are ``manifestly'' odd. We apply this lesson momentarily.

\subsection{Cutting and gluing}
\label{app:cutglue}

A final introductory observation concerns a basic cut-and-glue operation and its odd analogue.

\begin{figure}[htb]
\centering
\includegraphics[width=5in]{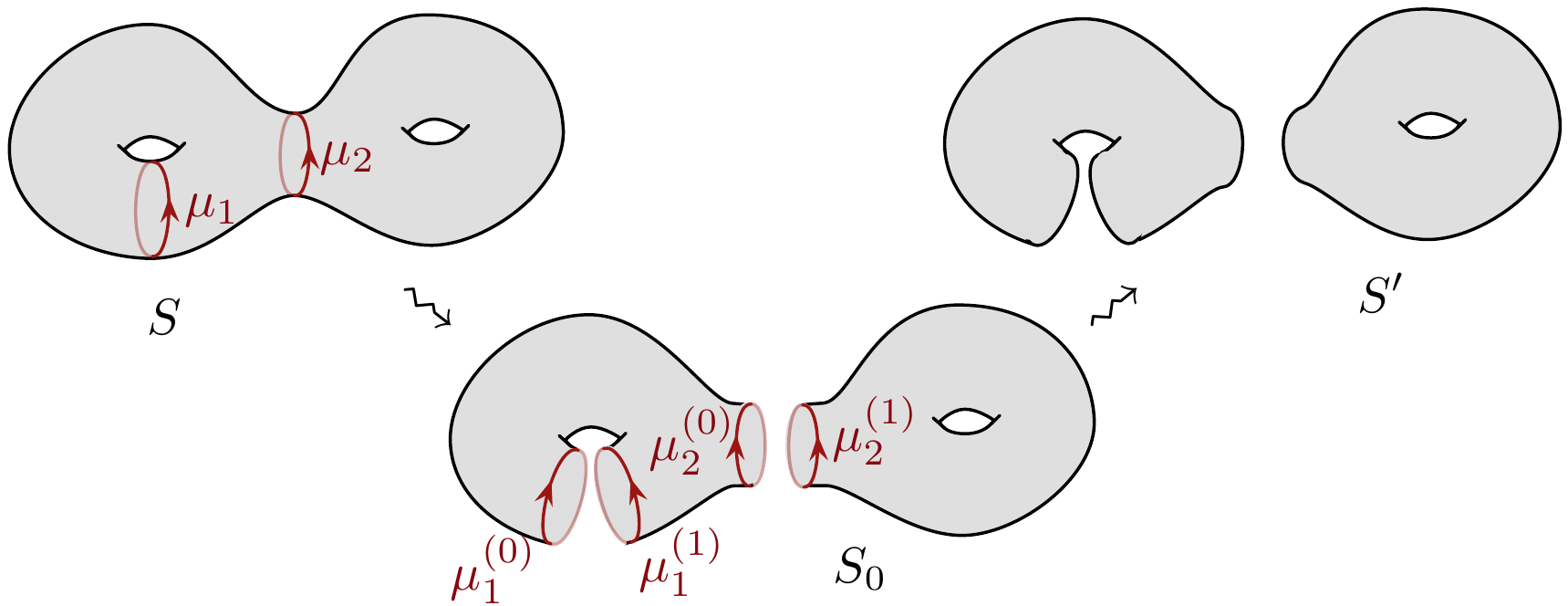}
\caption{Cutting $S$ along curves $\mu_i$ to form $S_0$, then filling in the boundaries of $S_0$ with discs to form the closed surface $S'$. In this case $G=\Z\langle\mu_1\rangle \subset H_1(S)\simeq \Z^4$ and $H_1(S')=H_1(S)/\!/G \simeq \Z^2$.}
\label{fig:cutglue}
\end{figure}

Let $S$ be a closed oriented surface and $\{\mu_i\}_{i=1}^r$ a collection of closed non-intersecting curves on $S$. Form a surface $S_0 = S\bs(\sqcup_i\mu_i)$ (possibly disconnected) by cutting along the $\mu_i$. Now $S_0$ has $2r$ circular boundary components, with images $\mu_i^{(0)}, \mu_i^{(1)}$ of the curves $\mu_i$ running along the boundaries. Cap off each boundary component of $S_0$ with a disc to form another closed, oriented surface $S'$, possibly with multiple components (Figure \ref{fig:cutglue}).

\begin{lemma} \label{lemma:S}

Let $G\subset H_1(S)$ be the subgroup of $H_1(S)$ generated by the classes of the $\mu_i$'s. Let $K = \ker\,\langle G,*\rangle|_{H_1(S)}$ be the subgroup of elements $H_1(S)$ that have zero intersection with all the $\mu_i$. Then there is a short exact sequence
\be 0 \to G \overset{i}\to K \overset{q}\to H_1(S') \to 0 \label{SES} \ee
that identifies $H_1(S')\simeq K/G =: H_1(S)/\!/G$ as a lattice symplectic quotient. In particular, the intersection form on $H_1(S')$ is induced by the quotient from the form on $H_1(S)$.

\end{lemma}

\emph{Proof.} 
First observe that since the $\mu_i$'s are non-intersecting, the intersection form vanishes on $G$, $\langle G,G\rangle = 0$. This it makes sense to define $K = \ker\,\langle G,*\rangle|_{H_1(S)}$, and $i:G\hookrightarrow K$ is just the inclusion.
Exactness of the rest of \eqref{SES} can be derived by comparing long exact sequences in relative homology for the pairs $(S,S_0)$ and $(S',S_0)$. However, it is useful to take a more concrete approach.

The map $q$ is defined as follows. Any homology class $[\gamma]\in K\subset H_1(S)$ can be represented by a curve (or sum of curves) $\gamma$ that does not intersect the $\mu_i$'s. Thus $\gamma$ can be naturally included in $S_0$ and (viewing $S_0$ as a subspace of $S'$) as a curve $\gamma' \subset S'$. We set $q([\gamma]):=[\gamma']$.
The only potential ambiguity in this definition comes from the $\mu_i$ themselves, which may be included either as $\mu_i^{(0)}$ or $\mu_i^{(1)}$ in $S_0$ and then in $S'$. But in $S'$ the curves $\mu_i^{(0)}, \mu_i^{(1)}$ bound discs, so they are contractible, and $q([\mu_i])=0$ unambiguously. Thus $q$ is well defined and its kernel includes all of $G$. Moreover, it preserves the intersection form.

Conversely, any $[\gamma']\in H_1(S')$ can be represented by a curve $\gamma'\subset S_0\subset S'$ (\ie\ a curve in the complement of the discs that filled in the holes of $S_0$). Viewing $S_0$ as a subspace of $S$, we can thus include $\gamma'$ as a curve $\gamma\subset S$, and try to define $q^{-1}([\gamma']) = [\gamma] \in K\subset H_1(S)$. There is ambiguity in this definition coming from how the representative $\gamma'$ is chosen to wind around the discs in $S'\bs S_0$; this ambiguity is precisely $G$. Therefore, we get a map $q^{-1}:H_1(S') \to K/G$ that provides an inverse to $q$, and demonstrates that $q:K/G \overset{\sim}\to H_1(S')$ is an isomorphism. \;$\square$ \medskip

We actually want a generalization of \eqref{SES} to odd homology. To this end, suppose that we perform a basic cut-and-glue operation $\CC\leadsto \CC_0 \leadsto \CC'$ along some non-intersecting closed curves $\ol\mu_i \subset \CC$. Suppose that the $\ol\mu_i$ are chosen to lie in the complement of the branching locus $\mathfrak b$ of a double cover $\Sigma\overset\pi\to \CC$. Then the cut-and-glue operation lifts to $\Sigma$. We first take all lifts of the $\ol \mu_i$ to $\Sigma$, noting that a given $\ol\mu_i$ may have a single lift $\mu_i$ if $\pi^{-1}(\mu_i)$ is connected and two lifts $\mu_i^\pm$ otherwise. We then cut $\Sigma$ along the $\mu_i^\pm,\mu_i$ to form $\Sigma_0$, and we fill in the boundaries of $\Sigma_0$ with discs to form $\Sigma'$. There are induced covering maps $\Sigma_0\overset\pi\to \CC_0$ and $\Sigma'\overset\pi\to \CC'$, such that each disc in $\CC'\bs\CC_0$ contains a new branch point and is covered by a single disc in $\Sigma'\bs\Sigma_0$ if and only if the corresponding cutting curve $\ol\mu_i$ has a single connected lift $\mu_i$ 
(otherwise a disc in $\CC'\bs\CC_0$ is covered by two disconnected discs in $\Sigma'\bs\Sigma_0$). We now have

\begin{lemma} \label{lemma:S-}

Let $G\subset H_1(\Sigma)$ denote the subgroup generated by the lifts $\mu_i^\pm,\mu_i$, and $K:=\ker\,\langle G,*\rangle|_{H_1(\Sigma)}$. Then $G$ and $K$ are fixed by the deck transformation $\sigma_*$, and we can take odd parts $G^-=\ker\,P_+|_G = G\cap H_1^-(\Sigma)$ and $K^-=\ker\,P_+|_K = K\cap H_1^-(\Sigma) = \ker\langle G^-,*\rangle|_{H_1^-(\Sigma)}$.
There is a complex
\be  0 \to G^- \overset{i^-}\to K^- \overset{q^-}\to H_1^-(\Sigma') \to 0 \label{SES-} \ee
with vanishing homology at all but the last spot, where the homology is 2-torsion. The maps $i^-,q^-$ preserve the intersection form. Therefore there is an injection $H_1^-(\Sigma)/\!/G^-=K^-/G^- \overset{q^-}\hookrightarrow H_1^-(\Sigma')$ with finite (2-torsion) cokernel, which preserves the intersection form.

If every cutting curve $\ol\mu_i\subset \CC$ has two distinct lifts $\mu_i^\pm$ to $\Sigma$, then \eqref{SES-} is exact, and $q^-:K^-/G^- \hookrightarrow H_1^-(\Sigma')$ is an isomorphism.

\end{lemma}

\emph{Proof.} First observe that $\sigma_*$ preserves $G$ because $\sigma_*\mu_i^\pm = \mu_i^\mp$ and $\sigma_*\mu_i=\mu_i$; and $\sigma_*$ preserves $K$ because it preserves the intersection form. So $G^-:=\ker\, P_+|_{G}$ and $K^-:=\ker\,P_+|_K$ make sense. 
Also, $K^-$ coincides with $\ker\,\langle G^-,*\rangle|_{H_1^-(\Sigma)}$: inclusion $\ker\,P_+|_K \subset \ker\,\langle G^-,*\rangle|_{H_1^-(\Sigma)}$ is obvious; conversely if $\alpha\in \ker\,\langle G^-,*\rangle|_{H_1^-(\Sigma)}$ then $\langle \mu_i^+-\mu_i^-,\alpha\rangle=0$ $\Rightarrow$ $\langle \mu_i^+,\alpha\rangle = \langle \mu_i^-,\alpha\rangle = \langle \sigma\mu_i^-,\sigma\alpha\rangle = -\langle \mu_i^+,\alpha\rangle$ $\Rightarrow$ $\langle \mu_i^\pm,\alpha\rangle=0$, and $\langle \mu_i,\alpha\rangle=0$ because $\mu_i$ is even and $\alpha$ is odd; so $\alpha\in \ker \langle G,*\rangle|_{H_1(\Sigma)}=K$, hence $\ker\,\langle G^-,*\rangle|_{H_1^-(\Sigma)}\subset \ker\,P_+|_K$.

Consider the exact sequence $0 \to G \overset{i}\to K \overset{q}\to H_1(\Sigma') \to 0$ from Lemma \ref{lemma:S}. Since $\sigma_*$ commutes with $i$ and $q$, we can take $\ker\,P_+$ of all groups and apply Lemma \ref{lemma:cx}(c) to obtain the complex \eqref{SES-}, with vanishing homology except perhaps at the last spot, where homology must be 2-torsion. The restricted maps $i^-$ and $q^-$ must preserve the intersection form, because $i$ and $q$ do.

In the case that every $\ol\mu_i\subset\CC$ has two distinct lifts, the discs in $\Sigma'\bs\Sigma_0$ contain no branch points. Thus (as at the end of Appendix \ref{app:cell}) every element of $H_1^-(\Sigma')$ can be represented by a sum of curves $\gamma' \subset \Sigma_0\subset \Sigma'$ that is manifestly odd, \ie\ $\sigma$ fixed $\gamma'$ setwise while reversing orientation. Then the inverse map $q^{-1}$ of Lemma \ref{lemma:S} sends $[\gamma']$ to an element $[\gamma]\in K/G$ that is \emph{also} manifestly odd, and its restriction $q^{-1}|_{H_1^-(\Sigma')}:H_1^-(\Sigma')\to K^-/G^-$ provides an inverse to $q^-$ in \eqref{SES-}, guaranteeing that $q^-:K^-/G^- \hookrightarrow H_1^-(\Sigma')$ is actually an isomorphism.

(If there are branch points on the discs in $\Sigma'\bs\Sigma_0$ this can fail. In order to represent all $[\gamma']\in H_1^-(\Sigma')$ by curves $\gamma'$ that are manifestly odd, the curves may have to go through the branch points --- and cannot be deformed into $\Sigma_0\subset \Sigma$ while keeping them manifestly odd. Thus while the inverse $q^{-1}|_{H_1^-(\Sigma')}:H_1^-(\Sigma')\to K/G$ still exists, its image may not be completely odd. This is the source of the 2-torsion in \eqref{SES-}.) \;$\square$

\section{Reconstructing framed flat connections}
\label{app:traffic}

Let $M$ by a framed, triangulated 3-manifold with $\pi_1(\CC_{\rm small})$  abelian.
We briefly explain how to uniquely reconstruct a framed flat connection $\CA\in \CP[\CC]\subset \CX[\CC]$ given functions $x_\pp$ satisfying the conditions on the RHS of \eqref{Pexplicit}, proving part (d) of Proposition \ref{prop:coords} in the process.

The basic idea (following \cite{FG-Teich} and precursors, \eg\ \cite{Fock-Teich}) is to use the functions to construct a distinguished set of projective bases (\ie\ frames) for the fiber of the bundle $E\to \CC$ at various points of $\CC$, together with maps between these bases corresponding to the parallel transport of $\CA$. There are several steps.

\begin{figure}[htb]
\centering
\includegraphics[width=5.8in]{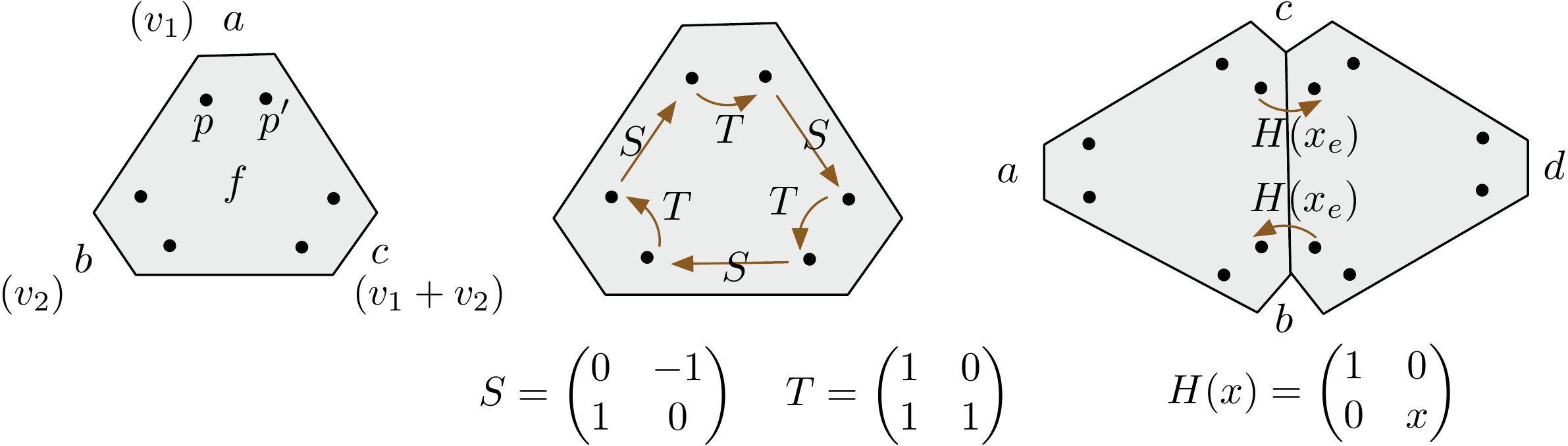}
\caption{Left: projective bases assigned to six points in a hexagon; for example at $p$ the basis is $(v_1,v_2)$ with $v_1\in a,\,v_2\in b,\,v_1+v_2\in c$; while at $p'$ the basis is $(v_1,w_2)$ with $v_1\in a,\,w_2\in b,\,v_1-w_2\in c$. Middle: $PGL(2)$ transformations among bases in any hexagon. Right: $PGL(2)$ transformations across big edges.}
\label{fig:bases}
\end{figure}

Start with a connected component of the big boundary $\CC_{\rm big}$, and a hexagon $f$ in the triangulation $\mb t_{\rm 2d}$. We assume that the flat connection is trivialized in the interior $f^\circ$ of $f$. Choose three arbitrary distinct lines $a,b,c$ in the fiber of $E$ over $f^\circ$, and make them framing lines at the three small edges of $f$ (this choice fixes gauge redundancy). Construct projective bases $b(p)$ for the fiber of $E$ over six points $p\in f^\circ$ as in Figure \ref{fig:bases}, as follows. For the element in $b(p)$ take any vector $v_1$ in the framing line on the small edge of $f$ closest to $p$. For the second element take a vector $v_2$ in the framing line on the small edge of $f$ second-closed to $p$, normalized that $v_1\pm v_2$ lies in the third framing line. The sign is specified by orientation, as in Figure \ref{fig:bases}.
The $PGL(2)$ transformations among the six bases in $f$ are given by matrices $S = \left(\begin{smallmatrix} 0 & -1 \\ 1 & 0 \end{smallmatrix}\right)$ and $T = \left(\begin{smallmatrix} 1 & 0 \\ 1 & 1 \end{smallmatrix}\right)$, as in the Figure, obeying $(ST)^3=1$ in $PGL(2)$.

Proceed to assign three framing lines and six projective bases to the rest of the hexagons in the connected component of $\CC_{\rm big}$, subject to the following rules: \\
1) framing lines at adjacent small edges (bounding a hole in $\CC_{\rm big}$) must agree; \\
2) the connection is trivialized in the interior of every hexagon, and the six bases there are constructed the same way they were for $f$, with the same $PGL(2)$ relations; \\
3) parallel transport across a big edge $e$ is given by $H(x_e) = \left(\begin{smallmatrix} 1 & 0 \\ 0 & x_e \end{smallmatrix}\right) \in PGL(2)$, where $x_e$ is the edge function\,. \\
Following these rules, there are no more arbitrary choices to be made, and we reconstruct a framed flat connection on $\CC_{\rm big}$ --- with framing lines at the holes of $\CC_{\rm big}$. We repeat for every connected component of $\CC_{\rm big}$.

Next, we extend the framed flat connection over the small boundary. Notice that the holonomy on a clockwise path around any hole on the big boundary takes the form
\be {\rm Hol}(hole) = \begin{pmatrix} 1 & 0 \\ * & \quad \prod_{\text{$e$ at hole}} x_e^{-1} \end{pmatrix}\,.\ee
This preserves the framing line on the boundary of the hole, rescaling it with (squared) eigenvalue $\prod_e x_e$. If the hole is to be filled in with a small disc, then $\prod_e x_e=1$ by \eqref{Pexplicit}, so we can uniquely extend the framed flat connection over a \emph{punctured} disc with unipotent holonomy at the puncture (as required for $\CX[\CC]$).

\begin{figure}[htb]
\centering
%\vspace{-.3cm}
\includegraphics[width=3.5in]{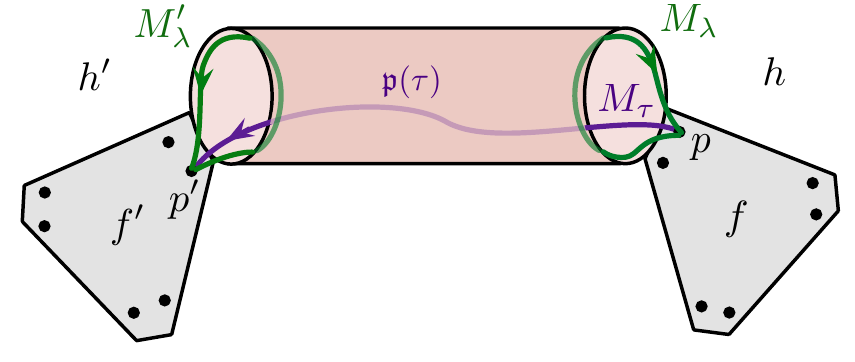}
%\vspace{-.3cm}
\caption{Reconstructing a framed flat connection on an annulus: the parallel transports must satisfy $M_\lambda = M_\tau^{-1}M_\lambda'M_\tau$.}
\label{fig:PGLannulus}
\end{figure}

If instead a pair of holes $h,h'$ are connected by a small annulus, we choose two hexagons $f,f'$ adjacent to the holes and a path $\pp_\tau$ running from $f$ to $f'$. Let $p\in f,\,p'\in f'$ be points with projective bases on the two sides of the annulus, as in Figure \ref{fig:PGLannulus}. Let $M_\lambda$, $M_\lambda'$ be the $PGL(2)$ holonomies around the ends of the annulus, with basepoints at $p,p'$, running clockwise (resp., counter-clockwise) from the viewpoint of $h$ (resp., $h'$). The holonomies are fixed to be
\be M_\lambda = \begin{pmatrix} 1 & 0 \\ a & x_\lambda^{-1} \end{pmatrix}\,,\quad M_\lambda' = \begin{pmatrix} 1 & 0 \\ a' & x_\lambda^{-1} \end{pmatrix}\,,\qquad x_\lambda = \prod_{\text{$e$ at $h$}}x_e = \prod_{\text{$e'$ at $h'$}}x_{e'}^{-1}\,\ee
for some (determined) $a,a'$. The function $x_\tau = x_{\pp(\tau)}$ fixes the eigenvalue of the $PGL(2)$ transformation $M_\tau$ between the bases at $p$ and $p'$, given by parallel transport along $\pp_\tau$,
\be M_\tau = \begin{pmatrix} 1 & 0 \\ t & x_\tau^{-1} \end{pmatrix} \label{det-t} \ee
for some undetermined $t$. In order to identify the framing lines at the two ends of the annulus and extend the framed flat connection over the annulus, it suffices to require $M_\lambda = M_\tau^{-1}M_\lambda' M_\tau$, which uniquely determines $t = (a'-ax_\tau^{-1})/(1-x_\lambda^{-1})$ as long as $x_\lambda\neq 1$.%
\footnote{Alternatively, one can think of $M_\tau$ as the $PGL(2)$ gauge transformation required to identify the holonomies $M_\lambda$ and $M_\lambda'$ at the two ends of an annulus. Determining $M_\tau$ fixes a relative gauge ambiguity between components of $\CC_{\rm big}$ that are connected by annuli.} 

This completes the reconstruction of $\CA$ over any connected component of $\CC$ containing both big and small parts. In particular, the $PGL(2)$ holonomy along any closed path may be obtained by combining the ``traffic rules'' of Figure \ref{fig:bases} for the big boundary with the transport $M_\tau$ along annuli. Triviality of the holonomy along any contractible path is ensured by the local identities $(ST)^3=H(x)SH(x)S = M_\lambda^{-1}M_\tau^{-1}M_\lambda' M_\tau=1$.

Finally, on an isolated small sphere the framed flat connection is trivial. (The connection itself is trivial, and its residual gauge symmetry can be used to trivialize the framing data.)
For a small torus, if we are given a choice of A and B cycles with (squared) holonomy eigenvalues $(x_\alpha,x_\beta)\neq (1,1)$, we may uniquely reconstruct the commuting holonomy matrices (modulo $PGL(2)$ gauge equivalence) as $M_\alpha = H(x_\alpha^{-1})$, $M_\beta = H(x_\beta^{-1})$.
We take the framing to be the unique eigenline with (squared) eigenvalues $(x_\alpha,x_\beta)$.

%\vspace{-.4cm}
\newpage

\bibliographystyle{JHEP_TD}
\bibliography{toolbox}

\end{document}